\documentclass[reqno]{amsart}

\usepackage[all]{xy}
\usepackage{graphicx}

\usepackage{amssymb}
\usepackage{amsmath}
\usepackage{mathrsfs}
\usepackage{epsfig}
\usepackage{amscd}
\usepackage{graphicx, color}

\usepackage{graphicx}
\usepackage{amssymb}
\usepackage{amsmath}
\usepackage{enumerate}
\usepackage{amsfonts}
\usepackage[mathcal]{euscript}
\usepackage{amsthm}
\usepackage{xcolor}
\usepackage{amscd}
\usepackage{mathrsfs}
\numberwithin{equation}{section}

\usepackage[colorlinks=true, pdfstartview=FitV, linkcolor=black,
			   citecolor=black, urlcolor=black]{hyperref}
\usepackage{bookmark}
\usepackage[titletoc]{appendix} 
\hypersetup{hypertexnames=false} 

          \theoremstyle{definition}
          \newtheorem{theorem}{Theorem}[section]
          
          \newtheorem{prop}[theorem]{Proposition}
          \newtheorem{lemma}[theorem]{Lemma}

          \newtheorem{cor}[theorem]{Corollary}

          \newtheorem{defn}[theorem]{Definition}
          \newtheorem{notation}[theorem]{Notation}
          \newtheorem{example}[theorem]{Example}
          \newtheorem{choice}[theorem]{Choice}

          \newtheorem{remark}[theorem]{Remark}
          \newtheorem{convention}[theorem]{Convention}
          
          \newtheorem{cond}[theorem]{Condition}

\numberwithin{equation}{section}

          \newcommand{\nc}{\newcommand}
          \nc{\DMO}{\DeclareMathOperator}	
          \nc{\commentout}[1]{}

          \nc{\newnotation}{\nomenclature}
          \nc{\wrap}{\cW}
          \nc{\Tw}{\mathsf{Tw}}
          \nc{\loc}{\mathsf{Loc}}
          \nc{\Top}{Top}
          \nc{\emb}{\mathsf{emb}}
          \nc{\ind}{\mathsf{Ind}}
          \nc{\Ind}{\mathsf{Ind}}
          \nc{\Loc}{\mathsf{Loc}}
          \nc{\Cob}{\mathsf{Cob}}
          \nc{\mul}{\mathsf{Mul}}
          \nc{\fat}{\mathsf{fat}}
          \nc{\cob}{\mathsf{Cob}}
          \nc{\coh}{\mathsf{Coh}}
          \nc{\Liouaut}{\Aut_{\mathsf{Liou}}}
          \nc{\idem}{\mathsf{Idem}}
          \nc{\sets}{\mathsf{Sets}}
          \nc{\near}{\mathsf{near}}
          \nc{\sing}{\mathsf{Sing}}
          \nc{\Sing}{\mathsf{Sing}}
          \nc{\perf}{\mathsf{Perf}}
          \nc{\block}{\mathsf{block}}
          \nc{\ssets}{\mathsf{sSets}}
          \nc{\cmpct}{\mathsf{cmpct}}
          \nc{\compact}{\mathsf{cmpct}}
          \nc{\pwrap}{\mathsf{PWrap}}
          \nc{\coder}{\mathsf{Coder}}
          \nc{\bimod}{\mathsf{Bimod}}
          \nc{\grmod}{\mathsf{GrMod}}
          \nc{\Morita}{\mathsf{Morita}}
          \nc{\morita}{\mathsf{Morita}}
          \nc{\spaces}{\mathsf{Spaces}}
          \nc{\pwrms}{\mathsf{PWrFuk}_{M,S}}
          \nc{\pwrmf}{\mathsf{PWrFuk}_{M,F}}
          \nc{\pwrapmf}{\mathsf{PWrFuk}_{M,F}}
          \nc{\fuk}{\mathsf{Fukaya}}
          \nc{\infwr}{\mathsf{InfWr}}
          \nc{\fukaya}{\mathsf{Fukaya}}
          \nc{\autml}{\mathsf{Aut}_{M,\Lambda}}
          \nc{\fukml}{\mathsf{Fukaya}_{M,\Lambda}}
          \nc{\fukmle}{\mathsf{Fukaya}_{M,\Lambda,\epsilon}}
          \nc{\fukmod}{\wrfukcompact(M)\modules}
          \nc{\lag}{\mathsf{Lag}}
          \nc{\lagm}{\lag_M}
          \nc{\lago}{\lag^o}
          \nc{\lagml}{\lag_{M,\Lambda}} 
          \nc{\lagmle}{\lag_{M,\Lambda,\epsilon}}
          \nc{\Fun}{\mathsf{Fun}}
          \nc{\fun}{\mathsf{Fun}}
          \nc{\vect}{\mathsf{Vect}}
          \nc{\chain}{\mathsf{Chain}}
          \nc{\chainn}{Chain}
          \nc{\wrfuk}{\mathsf{WrFukaya}}
          \nc{\wrfukcompact}{\mathsf{WrFukaya}_{\mathsf{cmpct}}}
          \nc{\pwrfuk}{\mathsf{PWrFukaya}}
          \nc{\inffuk}{\mathsf{InfFuk}}
          \nc{\pwrfukml}{\mathsf{PWrFukaya}_{M,\Lambda}}
          \nc{\inffukml}{\mathsf{InfFuk}_{M,\Lambda}}
          \nc{\nattrans}{\mathsf{NatTrans}}
          \nc{\corres}{\mathsf{Corres}}
          \nc{\fukep}{\fukaya_\Lambda(M,\epsilon)}
          \nc{\fukepop}{\fukaya_\Lambda(M,\epsilon)^{\op}}
          \nc{\lagep}{\lag_\Lambda(M,\epsilon)}
          \DMO{\cyl}{cyl} 
          \nc{\dbcoh}{D^b\mathsf{Coh}}
          \nc{\corr}{\mathsf{Corr}}
          \nc{\Liouauto}{{\Aut^o}}
          \nc{\Liouautb}{\Aut^{b}}
          \nc{\Liouautgr}{{\Aut^{gr}}}
          \nc{\Liouautgrb}{\Aut^{gr,b}}
          \nc{\Fuk}{\mathsf{Fuk}}

          \DMO{\im}{im}
          \DMO{\ev}{ev}
          \DMO{\stable}{Ex}
          \DMO{\inj}{inj}
          \DMO{\fib}{fib}
          \DMO{\conf}{Conf}
          \DMO{\chains}{Chains}
          \DMO{\cochains}{Cochains}
          \DMO{\cone}{Cone}
          \DMO{\Map}{Map}
          \DMO{\ran}{Ran}
          \DMO{\rot}{Rot}
          \DMO{\leg}{Leg}
          \DMO{\imm}{imm}
          \DMO{\adj}{adj}
          \DMO{\symp}{Symp}
          \DMO{\tree}{Tree}
          \DMO{\cube}{Cube}
          \DMO{\deep}{deep}
          \DMO{\back}{back}
          \DMO{\Hoch}{Hoch}
          \DMO{\front}{front}
          \DMO{\flow}{Flow}
          \DMO{\floer}{Floer}
          \DMO{\Maps}{Maps}
          \DMO{\exact}{exact}
          \DMO{\excess}{Excess}
          \DMO{\Decomp}{Decomp}
          \DMO{\decomp}{Decomp}
          \DMO{\collar}{collar}
          \DMO{\yoneda}{Yoneda}
          \DMO{\hamspace}{Ham}
          \DMO{\sympspace}{Symp}
          \DMO{\holomaps}{Holomaps}
          \DMO{\comp}{Comp}
          \DMO{\crit}{Crit}
          \DMO{\test}{{test}}
          \DMO{\sign}{sign}
          \DMO{\topp}{top}
          \DMO{\indx}{Index}
          \DMO{\Break}{Break} 
          \DMO{\zero}{zero} 
          \DMO{\ob}{Ob}
          \DMO{\gr}{Gr} 
          \DMO{\Gr}{Gr} 
          \DMO{\cl}{Cl} 
          \DMO{\grlag}{GrLag}
          \DMO{\Pin}{Pin}
          \DMO{\Graph}{Graph}
          \DMO{\pin}{Pin}
          \DMO{\gap}{Gap}
          \DMO{\Ex}{Ex}
          \DMO{\id}{id}
          \DMO{\End}{End}
          \DMO{\sym}{Sym}
          \DMO{\aut}{Aut}
          \DMO{\Aut}{Aut}
          \DMO{\haut}{hAut}
          \DMO{\hAut}{hAut}
          \DMO{\DK}{DK} 
          \DMO{\poly}{poly} 
          \DMO{\diff}{Diff}
          \DMO{\coll}{coll}
          \DMO{\dist}{dist} 
          \DMO{\coker}{coker} 
          \nc{\kernel}{\ker} 
          \DMO{\sspan}{span}
          \DMO{\hocolim}{hocolim}	
          \DMO{\holim}{holim}
          \DMO{\sk}{sk}

          \DMO{\ho}{ho}
          \DMO{\fin}{fin}
          \DMO{\tor}{Tor}
          \DMO{\ext}{Ext}
          \DMO{\ret}{Ret}
          \DMO{\ham}{Ham}
          \DMO{\con}{con}
          \DMO{\leaf}{leaf}
          \DMO{\supp}{supp}
          \DMO{\edge}{edge}
          \DMO{\colim}{colim}
          \DMO{\edges}{edges}
          \DMO{\Image}{image}
          \DMO{\roots}{roots}
          \DMO{\height}{height}
          \DMO{\finmod}{FinMod}
          \DMO{\leaves}{leaves}
          \DMO{\planar}{planar}
          \DMO{\vertices}{vertices}

\nc{\norm}[2]{{ \ensuremath{\|} #1 \ensuremath{\|}}_{#2}}
\nc{\Dbar}[1]{\ensuremath{{\bar{\partial}}_{#1}}}
\nc{\Ce}{\ensuremath{\mathbb{C}}}
\nc{\B}{\ensuremath{\mathbb{B}}}
\nc{\osc}{\operatorname{osc}}
\nc{\leng}{\operatorname{leng}}

          \nc{\lagg}{\lag^{\cG}}
          \nc{\iso}{\mathsf{Iso}}
          \nc{\Set}{\mathsf{Set}}
          \nc{\Ass}{\mathsf{ \bf Ass}}
          \nc{\Mod}{\mathsf{Mod}}
          \nc{\modules}{\mathsf{Mod}}
          \nc{\sset}{\mathsf{sSet}}
          \nc{\liou}{\mathsf{Liou}}
          \nc{\poset}{\mathsf{Poset}}
          \nc{\trno}{T^*\RR^n_{\geq 0}}
          \nc{\spectra}{\mathsf{Spectra}}
          \nc{\tensorfin}{\tensor^{\fin}}
          \nc{\lagptg}{\lag_{pt,pt}^{\cG}}
          \nc{\Fin}{\mathcal{F}\mathsf{in}}
          \nc{\lagnl}{\lag_{N,\Lambda}}
          \nc{\lagmlg}{\lag_{M,\Lambda}^{\cG}}
          \nc{\lagsplit}{\lag^{\mathsf{split}}}
          \nc{\lagktimes}{(\lag^{\dd k})^\times}
          \nc{\lagplanar}{\lag^{\times,\planar}}
          \nc{\Cont}{\text{\rm Cont}}
          \nc{\Ham}{\text{\rm Ham}}
          \nc{\Dev}{\text{\rm Dev}}
          \nc{\Lin}{\text{\rm Lin}}
          \nc{\Int}{\text{\rm Int}}
          \nc{\Hom}{\text{\rm Hom}}
          \nc{\Chord}{\text{\rm Chord}}
          \nc{\nbhd}{\mathcal{N}\text{\rm{bhd}}}

          \nc{\onef}{1_{\fukaya}}
          \nc{\smsh}{\wedge}
          \nc{\un}{\underline}
          \nc{\xto}{\xrightarrow}
          \nc{\xra}{\xto}
          \nc{\tensor}{\otimes}
          \nc{\del}{\partial}
          \nc{\dd}{\diamond}
          \nc{\tri}{\triangle}
          \nc{\bb}{\Box}
          \nc{\into}{\hookrightarrow}
          \nc{\onto}{\twoheadrightarrow}
          \nc{\contains}{\supset}
          \nc{\transverse}{\pitchfork}
          \nc{\uncirc}{\underline{\circ}}
          \nc{\thetacontact}{\theta} 

          \nc{\Jbar}{\overline{J}}
          \nc{\Fbar}{\overline{F}}
          \nc{\delbar}{\overline{\del}}
          \nc{\thetabar}{\overline{\theta}}
          \nc{\omegabar}{\overline{\omega}}
          \nc{\Liou}{\text{\rm Liou}}

          \nc{\Yhat}{\widehat{Y}}

          \nc{\Mliou}{M}

          \nc{\vece}{ {\vec \epsilon}}	
          \nc{\vecd}{ {\vec \delta}}
          \nc{\ov}{\overline}
          \DMO{\op}{op}
          \nc{\opp}{ ^{\op}}

          \nc{\hiro}{\textcolor{blue}}
          \nc{\YG}{\textcolor{orange}}

          \nc{\eqn}{\begin{equation}}
          \nc{\eqnn}{\begin{equation}\nonumber}
          \nc{\eqnd}{\end{equation}}
          \nc{\enum}{\begin{enumerate}}
          \nc{\enumd}{\end{enumerate}}
          \nc{\beastar}{\begin{eqnarray*}}
          \nc{\eeastar}{\end{eqnarray*}}

\numberwithin{equation}{section}

\def\R{{\mathbb R}}
\def\osc{{\hbox{\rm osc }}}

\def\C{{\mathbb C}}
\def\R{{\mathbb R}}

\def\11{{\mathbb I}}

\def\Jbar{{\widetilde J}}

\def\delbar{{\overline \partial}}
\def\dudtau{{\frac{\del u}{\del \tau}}}
\def\dudt{{\frac{\del u}{\del t}}}

          \def\dudtau{\frac{\del u}{d\tau}}
          \def\dudt{\frac{\del u}{dt}}

          \def\cC{\mathcal C}\def\cD{\mathcal D}
          \def\cF{\mathcal F}\def\cG{\mathcal G}
          \def\cJ{\mathcal J}
          \def\cN{\mathcal N}\def\cO{\mathcal O}
          
          \def\cW{\mathcal W}\def\c\Mliou{\mathcal \Mliou}

          \def\RR{\mathbb R}
          \def\\Mliou\Mliou{\mathbb \Mliou}
          \def\ZZ{\mathbb Z}

          \def\s\Mliou{\mathsf \Mliou}

          \def\b\Mliou{\mathbf \Mliou}

          \def\f\Mliou{\mathfrak \Mliou}

\def\Jbar{{\widetilde J}}

\def\delbar{{\overline \partial}}

\def\b{\beta}
\def\c{\chi}
  
\def\e{\varepsilon} 
\def\f{\phi}

\def\s{\sigma}

\def\CH{{\mathcal H}}

\def\CJ{{\mathcal J}}

\def\CW{{\mathcal W}}

\def\darr#1{\raise1.5ex\hbox{$\leftrightarrow$}
\mkern-16.5mu #1}

\def\roughly#1{\raise.3ex\hbox{$#1$\kern-.75em
\lower1ex\hbox{$\sim$}}}

\def\opname#1{\mathop{\kern0pt{\rm #1}}\nolimits}
\def\Re{\opname{Re}}

\def\End{\opname{End}}
\def\dim{\opname{dim}}

\def\dist{\opname{dist}}

\def\supp{\operatorname{supp}}

\def\Dev{\operatorname{Dev}}

\def\leng{\operatorname{leng}}

\def\End{\operatorname{End}}

\def\Aut{\operatorname{Aut}}
\def\coker{\operatorname{Coker}}
\def\span{\operatorname{span}}

\def\Cont{\operatorname{Cont}}

\def\Sing{\operatorname{Sing}}

\def\Image{\operatorname{Image}}
\def\ev{\operatorname{ev}}
\def\Int{\operatorname{Int}}

\def\ben{\begin{enumerate}}
\def\een{\end{enumerate}}
\def\be{\begin{equation}}
\def\ee{\end{equation}}
\def\bea{\begin{eqnarray}}
\def\eea{\end{eqnarray}}
\def\beastar{\begin{eqnarray*}}
\def\eeastar{\end{eqnarray*}}
\def\bc{\begin{center}}
\def\ec{\end{center}}

\renewcommand{\b}{\beta}

\def\Hoch{{\tt Hoch}}

\def\Cont{\operatorname{Cont}}

\def\Sing{\operatorname{Sing}}

\def\Ham{\operatorname{Ham}}

\def\Graph{\operatorname{Graph}}
\def\id{\text\rm{id}}

\begin{document}

\title{Geometry of Liouville sectors and the maximum principle}

\author{Yong-Geun Oh}
\address{Center for Geometry and Physics, Institute for Basic Science (IBS),
77 Cheongam-ro, Nam-gu, Pohang-si, Gyeongsangbuk-do, Korea 790-784
\& POSTECH, Gyeongsangbuk-do, Korea}
\email{yongoh1@postech.ac.kr}
\thanks{This work is supported by the IBS project \# IBS-R003-D1}


\begin{abstract}
We introduce a new package of Floer data of \emph{$\lambda$-sectorial almost complex structures}
$J$ and \emph{sectorial Hamiltonians} $H$ on the Liouville sectors introduced in \cite{gps, gps-2}
the pairs of which are amenable to the maximum principle
for the analysis of pseudoholomorphic curves relevant to the studies
of wrapped Fukaya categories and of symplectic cohomology. It is also amenable to the strong maximum principle
in addition when paired with cylindrical Lagrangian boundary conditions.
The present work answers to a question raised by Ganatra-Pardon-Shende in \cite{gps} concerning
a characterization of almost complex structures and Hamiltonians in that
all the relevant confinement results in the studies of wrapped Fukaya categories,
symplectic cohomology and closed-open (and open-closed) maps between them can be uniformly
established via the maximum principle through \emph{tensorial calculations, Hamiltonian calculus and sign considerations}
without doing any estimates. Along the way, we prove the existence of a pseudoconvex pair $(\psi,J)$ such that
\emph{$J$ is
$d\lambda$-tame and $\psi$ is an exhaustion function of $\nbhd(\del_\infty M \cup \del M)$ that also satisfies
the equation
$-d\psi \circ J = \lambda$ thereon for any Liouville sector with corners $(M,\lambda)$}.
\end{abstract}


\maketitle

\tableofcontents

\section{Introduction}

In \cite{oh-tanaka:liouville-bundles}, Tanaka and the present author constructed an unwrapped Floer theory
for bundles of Liouville manifolds. The output was a collection of unwrapped Fukaya categories associated to fibers of a Liouville bundle, along with a compatibility between two natural constructions of continuation maps. This set-up allowed
for them to make the construction of Floer-theoretic invariants
of smooth group actions on Liouville manifolds, and they exploited these constructions in~\cite{oh-tanaka:actions} to construct homotopically coherent actions of Lie groups on wrapped Fukaya categories, thereby proving a version of a conjecture from Teleman's 2014 ICM address.

In \cite{oh-tanaka:liouville-bundles}, the authors claimed that the same results hold for
the case of Liouville sectors~\cite{gps}. An anonymous referee pointed out a gap involving the existence of almost-complex structures suitable for maximum-principle-type arguments (as opposed to monotonicity arguments). The present paper fills this gap.

We introduce a new framework of \emph{sectorial Floer data} for Liouville sectors with corners
which is first amenable to the (interior) maximum principle for the relevant pseudoholomorphic curves.  It is also amenable
to the strong maximum principle
with $Z$-invariant-at-infinity Lagrangian boundary condition. Thirdly
it also pairs well with the \emph{sectorial Hamiltonians} which are
also introduced in the present paper in that the image of the
relevant Hamiltonian-perturbed pseudoholomorphic curves
are confined away from the preferred direction of $\del \Mliou$.
(See~\eqref{eqn.pi-to-C} and Remark \ref{rem:difficulty}
for the explanation on what the main
relevant issues and difficulties are to identify such a class.)

\subsection{Pseudoconvex pairs and end-profile functions}

To apply a maximum principle argument--- thereby preventing pseudoholomorphic curves from coming close to
the boundary $\del_\infty M \cup \del M$--- we need to use a barrier function
$$
\psi : \Int(M) \bigcap \nbhd(\del_\infty M \cup \del M) \to \RR
$$
that is plurisubharmonic with respect to given almost complex structure $J$.
The following is one of our main results of the present paper.

\begin{theorem}\label{thm:barrier-function-intro}
Let $M$ be a Liouville sector (possibly with corners) with Liouville form $\lambda$.
Then there exists a pair $(\psi,J)$ such that
\begin{enumerate}
\item The function $\psi : \Int(M) \bigcap \nbhd(\del_\infty M \cup \del M) \to \RR$ is an exhaustion function
with compact level sets that blows-up along $\del_\infty M \cup \del M$.
\item $J$ is a $d\lambda$-tame almost-complex structure defined on $\Int M$
such that $\psi$ is plurisubharmonic with respect to $J$, i.e., satisfies
$$
-d(d\psi \circ J) \geq 0
$$
as a $(1,1)$-current.
\item The function $\psi \circ u$ is amenable to the strong maximum principle for any
$J$-holomorphic map $u:\Sigma \to M$ satisfying the boundary condition $u(\del \Sigma) \subset L$ for
any Lagrangian submanifold $L$ that is $Z$-invariant at infinity.
\end{enumerate}
\end{theorem}
We call a pair $(\psi,J)$ satisfying the first two conditions a \emph{pseudoconvex pair}.
(Definition~\ref{defn:J-convex}).

\begin{remark}
Recall that the second condition is the well-known standard condition both in several complex variables and
in symplectic geometry in relation to the study of pseudoholomorphic curves.
The first condition is also a  natural condition that often appears in the $L^2$-theory of
several complex variables, especially on
\emph{pseudoconvex domains}. The function of the type $\psi$ is used as the weight function in the definition of
relevant weighted Sobolev spaces. (See \cite{hormander:SCV} for example.)
\end{remark}

The fundamental difficulty of finding such a pair $(\psi, J)$ in relation to the
wrapped Fukaya category on Liouville sectors from \cite{gps} lies in the requirement
of making the  pseudoconvex pair $(\psi, J)$ also amenable to the strong maximum principle
\begin{itemize}
\item both for the boundary value problem of $J$-holomorphic curves under
the boundary condition of the \emph{$Z$-invariant-at-infinity} Lagrangian submanifolds $L$, and
\item for the Hamiltonian-perturbed Cauchy-Riemann equations appearing in the
construction of Floer cohomology and its continuation maps.
\end{itemize}
Confinement results for the relevant pseudoholomorphic curves are crucial for the construction of
\emph{wrapped Fukaya category} obtained by the wrapping localization via the relevant wrapping Hamiltonians.
Note that the confinement result for the latter case cannot be proved by a direct monotonicity argument
which complicates the matter and forces \cite{gps} to employ rather complex notion of
dissipative Hamiltonians for the purpose.

For the purpose of making the pseudoconvex pair $(\psi,J)$ to satisfy the first requirement, it is enough
for $(\psi,J)$ to  the following stronger condition.

\begin{defn}[Liouville pseudoconvex pairs]
\label{defn:liouville-pseudoconvex-intro} We call a pseudoconvex pair $(\psi,J)$ \emph{Liouville pseudoconvex}
if it satisfies
\eqn\label{eq:liouville-pseudoconvex}
-d\psi \circ J = \lambda.
\eqnd
\end{defn}

We note from the definition of Liouville sectors \cite{gps} that geometric natures of
$\del_\infty M$ and $\del M$ are quite different: $\del_\infty M$ is of contact-type while $\del M$
is `flat' in that it is contained in a level set of `linear' function
$R: \nbhd(\del M) \to \RR$. More generally for the Liouville sectors
with corners, the boundary $\del M$ carries the structure of
\emph{pre-symplectic manifolds with corners} as any \emph{coisotropic submanifold of $ M$
with corners} does. This different geometric nature of $\del_\infty M$ and $\del M$
is one source of difficulty constructing
such a pseudoconvex pair and requires us to carefully analyze the geometry of neighborhoods of $DM$,
especially when the sector $\del M$ has corners.
We write
\eqn\label{eq:DM}
DM: = \del_\infty M \cup \del M
\eqnd
and call it the \emph{full boundary} of the ideal completion of $M$.

Construction of aforementioned pseudoconvex pair will be carried out by carefully designing the \emph{end-profile function},
that provide a smoothing---of both $\del M$ (as $M$ may be a manifold with corners) and $\del_\infty M \cap \del M$ (where the ideal contact boundary ``meets'' the boundary strata of $M$). Our Floer data $(J,H)$,
depend on this choice, so we outline the construction of this end-profile function and
this dependency briefly. Details are in Section~\ref{sec:barrier-functions}.

\begin{defn}[Smoothing profile]\label{defn:smoothing-profile-intro}
Given a Liouville sector $M$, a {\em smoothing profile} for $M$ is the data of the following choices:
\enum
\item A splitting data $(F_\delta,\{(R_{\delta,i},I_{\delta,i})\}_{i=1}^{k_\delta})$ at each corner $C_\delta$ of $M$
(See Proposition \ref{prop:gps} and Condition \ref{cond:smoothing-profile} below),
\item Convex functions $\varphi= \varphi_k: \RR^k \to \RR$ for every codimension $k$ stratum of $\del M$
which is a family defined on the universal model $\RR^{k+1}$ for a chart of manifolds with corners.
(See Definition~\ref{defn:symmetric-convex} for the precise definition), and
\item A choice of a contact-type hypersurface $S_0$ and its associated symplectization coordinate $s$ satisfying $Z[s] =1$.
\enumd
\end{defn}
Given a smoothing profile, we will first define a canonical end-profile function denoted by
$$
{\frak s}_\varphi
$$
which is obtained as follows.
With a choice of a symplectization radial function $s$,
we consider functions which, near a point of the corner $\del M \cap \del_\infty M$, take the form
	\eqn\label{eqn. local form of frak s}
	- \log \varphi \left(R_1,\ldots, R_k,e^{-s}\right).
	\eqnd
We remark that \emph{these functions become the radial function $s$ away from a neighborhood
of the given corner and so they canonically glue together.}
By gluing these locally defined functions together, which become automatically
pseudoconvex on $\nbhd(\del M \cup \del_\infty M)$, we obtain a single end-profile function
\eqn\label{eqn. frak s phi}
{\frak s}_\varphi: \Int M \cap \nbhd(\del_\infty M \cup \del M) \to \RR.
\eqnd
Its level sets provide smoothings not just near $\del M$, but also near the ideal corners $\del_\infty M \cap \del M$.

\begin{remark}[Sectorial almost complex structures versus $\lambda$-sectorial ones]
Given a smoothing profile and its associated end-profile function, we will
define a plentiful and contractible collection of almost complex structures, called the \emph{sectorial almost-complex structures}
in \cite{oh:gradient-sectorial}, where they are paired with a class of Lagrangian branes that are called \emph{gradient-sectorial}
Lagrangians. In the present paper, we still adopt the standard class of \emph{cylindrical (i.e., $Z$-invariant-at-infinity)
Lagrangians} as the branes of $\Fuk(M)$. However this choice of sectorial ones \emph{cannot} be paired with
the \emph{$Z$-invariant Lagrangians}  for the
application of the \emph{strong maximum principle}. It turns out constructing such a Liouville pseudoconvex pair is
a highly nontrivial task and forces us to consider a deformation of Liouville forms which changes
the given Liouville structure of $(M,\lambda)$ along the way as a step towards the final construction thereof.
\end{remark}

\subsection{$\lambda$-sectorial almost complex structures}

Let $(M,\lambda)$ be equipped with a smoothing profile and
let ${\mathfrak s}_\varphi$ be its end-profile function.

For the construction of wrapped Fukaya category $\Fuk(M)$ on a Liouville manifold or a sector,
 the choices of almost complex structures and Hamiltonians are made so that
some confinement theorems for the relevant (perturbed) pseudoholomorphic curves hold.
To achieve these confinement theorems via
the maximum and the strong maximum principle, we will  require $J$ to satisfy
\eqn\label{eq:-dsJ=lambda}
-ds\circ J = \lambda
\eqnd
at infinity. It is important to separate out this property out of the common definition of
$Z$-invariant tame almost complex structure in the literature as follows.

\begin{defn}[$J$ being contact type]\label{defn:contact-type-intro} Let $(M,\lambda)$ be a Liouville
manifold equipped with a symplectization radial function $s$
near infinity.  We call a tame almost complex structure $J$
\emph{$\lambda$-contact type}, if it satisfies
$$
-ds \circ J = \lambda.
$$
We denote by $\CJ_{\lambda}$ the set of all contact-type almost
complex structures on $(M,\lambda)$.
\end{defn}

We also require $J$ to satisfy
\eqn\label{eq:-dRJ=lambda-df}
-dR \circ J = \pi_F^*\lambda_F + \pi_\C^*\lambda_\C (= \lambda -df)
\eqnd
near boundary, away from the corner $\nbhd(\del_\infty M \cap \del M)$.
(See Proposition \ref{prop:gps} for the unexplained notations here.)
We mentioned that the expression of \eqref{eq:-dRJ=lambda-df} is made only \emph{near $\del M$}.

The main task then is to interpolate
the two requirements on the intersection $\nbhd(\del_\infty M \cap \del M)$ in the way that relevant maximum
and strong maximum principles are still applicable.
\emph{This construction of $J$ could not be done without going through a careful pointwise consideration
of $J$} to reveal what presents the obstruction to interpolating the two conditions so that $Z$-invariant
Lagrangian boundary condition becomes amenable to the strong maximum principle. In this regard, \emph{Section 8 is the central
section of the paper where basically Theorem \ref{thm:interpolation} and Proposition \ref{prop:linear-independence} are the main
 technical results to be established.}
In some way, most of the
preparatory geometric results on the Liouville sectors given
in Part I are geared to prove these two results. In addition,
a carefully constructed convex smoothing of the corner of $\R^n_{\geq 0}$
in Section \ref{sec:corner-smoothing-functions} is a key ingredeint for the proof of Proposition \ref{prop:linear-independence}.

\begin{remark} In hindsight, our construction of barrier function $\psi$ seems to be
in the sprit similar to the construction of $i$-convex hypersurfaces in
\cite[Chapter 4]{cieliebak-eliashberg}. In our case, we will need
to construct the pair $(\psi,J)$ simultaneously so that
\begin{enumerate}
\item $\psi$
becomes  $J$-pseudoconvex  for a large contractible
family of tame almost complex structures, not just one.
\item $(\psi,J)$ is a Liouville pseudoconvex  pair in the sense of
Definition \ref{defn:liouville-pseudoconvex-intro} so that the
strong maximum principle for the (asymptltically) $Z$-invariant
Lagrangian submanifolds not intersecting $\del M$, i.e., for any
object of the wrapped Fukaya category $\CW(M)$ of \cite{gps}.
\end{enumerate}
\end{remark}

We will achieve the construction of such $J$ in two steps.

\subsubsection{Definition of $\kappa$-sectorial almost complex structures}

In the first step, we interpolate the two requirements \eqref{eq:-dsJ=lambda}, \eqref{eq:-dRJ=lambda-df}
along the corner $\nbhd(\del_\infty M \cap \del M)$ by considering the interpolation
\eqn\label{eq:lambda-kappa}
\lambda_\kappa: = \lambda - d((1-\kappa) f)
\eqnd
for a suitable cut-off function $\kappa: \RR^2 \to [0,1]$ satisfying
that $\kappa = 0$ on a neighborhood of
$\nbhd(\del M) \cap F_0 \times \{|I| \leq N_0\}  \cap \times \{s \geq N_2\}$
and $\kappa = 1$ further away from $\nbhd(\del M)$ and on $(F \setminus F_0) \times \C_{{\text{\rm Re}} \geq 0}$.
(See Definition \ref{defn:kappa} for the
precise description of $\kappa$.) We emphasize the facts that this form is globally defined and that
the perturbation term $- d((1-\kappa f))$ is \emph{not} compactly supported. In particular,
$\lambda_\kappa$ is not Liouville equivalent to $\lambda$ in the usual sense.

The following is the first step towards the definition of sectorial (resp. $\lambda$-sectorial) almost complex structures. This class of almost complex structures is introduced only as a convenient vehicle for the proof of
existence of  $\lambda$-sectorial almost complex structures, the main class of our interest
in the present paper. Note that each splitting data
$$
U: = \nbhd(\del M) \cong F \times \C_{\text{\rm Re}\geq 0}
$$
provides an obvious foliation $\cF_F$ on $U$ whose leaves are given by
the submanifolds $F \times \{(x,y)\}$, $(x,y) \in \C_{\text{\rm Re}\geq 0}$.
 Throughout the paper, we will denote by
$$
\nbhd^Z(\cdot)
$$
a $Z$ neighborhood of any subset $(\cdot) \subset M$ that is invariant
at infinity.

\begin{defn}[$\kappa$-sectorial almost complex structures]\label{defn:sectorial-J-intro}
Let $\kappa:\RR \to [0,1]$ be a cut-off function as above.
We call a $d\lambda$-tame almost complex structure $J$ on $M$ \emph{$\kappa$-sectorial}
if $J$ satisfies the following:
\begin{enumerate}
\item {\textbf{[$\cF_F$ is $J$-complex]}} In a
neighborhood $\nbhd^Z(\del \Mliou)$ of $\del \Mliou$, we require
    \eqn\label{eq:J-versus-JF-intro1}
    J\left(T^*F \oplus 0_{\text{\rm span}\{d\mu_{i}, d\nu_{i}\}_{i=1}^k}\right)
    \subset T^*F \oplus 0_{\text{\rm span}\{d\mu_{i}, d\nu_{i}\}_{i=1}^k},
    \eqnd
    and $J$ restricts to an almost complex structure of contact-type on $F$.
\item {\textbf{[$d{\mathfrak s}_\varphi$ is $J$-dual to $\lambda_\kappa$]}}
On $\nbhd(\del \Mliou \cup \del_\infty M) \cap \Int M$,
\eqn\label{eq:ds-dual-J}
- d{\mathfrak s}_\varphi \circ J = \lambda_\kappa
\eqnd
for the deformed Liouville form $\lambda_\kappa$.
\end{enumerate}
\end{defn}
Obviously any $\kappa$-sectorial almost complex structure forms a pseudoconvex pair
with $\psi = \frak s_{\varphi}$ since the above
duality requirement implies
$$
-d(d{\frak s}_\varphi \circ J)  = d\lambda \geq 0.
$$
Therefore the main task then is to adjust the pair $({\frak s}_\varphi,J)$ with $\kappa$-sectorial $J$
to a new pair $(\psi,J')$ that \emph{also satisfies}
$$
-d\psi \circ J' = \lambda.
$$
\begin{remark}\label{rem:difficulty}
\begin{enumerate}
\item It is easy to check that the condition \eqref{eq:J-versus-JF-intro1} is equivalent to
the property of $\cF_F$ that its leaves are $J$-complex submanifolds of $U =\nbhd_{2\varepsilon_0}(\del M)$.
(See Section \ref{subsec:kappa-sectorial} for more discussion on this.)
\item We would like to alert readers that the hypersurfaces ${\mathfrak s}_\varphi^{-1}(r)$
may not be of contact-type for the originally given Liouville form $\lambda$
in that they may not be transversal to the Liouville vector field $Z$
near the smoothing corners. This will create some difficulty in
applying the strong maximum principle against the end-profile function
because \emph{the intersection $L \cap {\mathfrak s}_\varphi^{-1}(r)$ may not be Legendrian}.
\item
We compare this failure of being of contact-type with the fact that
every Liouville sector $M$ \emph{admits a convex completion} $\overline M$
(i.e., admits an exhaustion function $\psi$ whose level sets are compact hypersurfaces of contact-type),
\emph{after the originally given Liouville sector
is sufficiently enlarged in the horizontal direction}. (See \cite[Lemma 2.31]{gps}.)
It would be interesting to quantify how much room we need to enlarge to make the boundary of
an extended Liouville sector inside the convex completion transversal to the
Liouville vector field in the completion defined in \cite[Lemma 2.31]{gps}.
\end{enumerate}
\end{remark}

The above definition already makes the function ${\frak s}_\varphi$
amenable to the maximum principle for a $J$-holomorphic curves.
\emph{However since the right hand side of \eqref{eq:ds-dual-J} is perturbed from $\lambda$,
one cannot apply  the strong maximum principle for the class of usual
$Z$-invariant-at-infinity Lagrangian branes on a neighborhood $\del_\infty M$ especially near the corner
$\del_\infty M \cap \del M$.}

\begin{remark}
At this stage, there are two possible routes one could take: one is to change the
objects of $\Fuk(M)$ from the usual $Z$-invariant Lagrangians to some other types of
Lagrangians amenable to the strong maximum principle. This route is
taken in the present author's article \cite{oh:gradient-sectorial}, \cite{choi-oh:construction}.
\end{remark}

In the present paper,
we do not change and keep the standard $Z$-invariant Lagrangians as the objects of $\Fuk(M)$.
Because of that, we need one more preparation  by deforming $\lambda_\kappa$ back to $\lambda$
before giving the final definition of the class of almost complex structures that we are searching for.

\subsubsection{Infinity-moving Liouville isotopy}

By definition,
$\lambda_t = \lambda$ on $M \setminus \nbhd(\del M)$,
$\lambda_t = \pi_\C^*\lambda_\C + \pi_F^*\lambda_F (= \lambda - df)$
 on $(F \setminus F_0) \times \C_{\text{\rm Re} \geq 0}$ and their $\C$-components
coincide with $\pi_\C^*\lambda_\C$. (We note $d_\C f = 0$ on $\{|I| \geq N_0\}$ and $\supp f \subset F_0 \times \C_{\text{\rm Re} \geq 0}$.)
We also have $d\lambda_t = d\lambda = \omega$ for all deformations $\lambda_t$ with $0 \leq t\leq 1$.
We refer readers to Subsection \ref{subsec:gray-contactification} for the precise description of this deformation.

The following deformation result is an important ingredient needed towards our final
definition.

\begin{prop}[Proposition \ref{prop:deformation}]\label{prop:deformation-intro}
[Liouville diffeomorphism]
Let $\omega = d\lambda (= d\lambda_\kappa)$ with $\lambda, \, \lambda_\kappa$
be as above. Then there exists a
(\emph{not necessarily compactly supported}) Liouville diffeomorphism
$$
\phi_\kappa: (M, \lambda) \to (M,\lambda_\kappa)
$$
that satisfies the following on a neighborhood $\nbhd(\del_\infty M \cup \del M)$:
\begin{enumerate}
\item $\phi_\kappa^*\lambda_\kappa =  \lambda$.
\item $\phi_\kappa$ is the time-one map of a time-dependent, not necessarily compactly supported,
 vector field $X_t$.
\end{enumerate}
\end{prop}

Utilizing this proposition, we now arrive at the final definition of
the class of $\lambda$-sectorial almost complex structures that includes all
the pull-backs of $\kappa$-sectorial ones provided in
Definition \ref{defn:J-for-corners} by the diffeomorphisms.

\begin{defn}[$\lambda$-sectorial almost complex structures]\label{defn:sectorial-J-intro}
We equip $(M,\lambda)$ with a smoothing profile and a $\kappa$-deformation $\lambda_\kappa$ of $\lambda$.
Consider the \emph{$\lambda_\kappa$-wiggled end-profile function} defined by
$$
{\mathfrak s}_{\varphi,\kappa}: = {\mathfrak s}_\varphi \circ \phi_\kappa.
$$
We call an $\omega$-tame almost complex structure $J$ on $M$ \emph{$\lambda$-sectorial}
if $J$ satisfies
$$
- d{\mathfrak s}_{\varphi,\kappa} \circ J = \lambda
$$
on a neighborhood $\nbhd^Z(\del_\infty M \cup \del M)$.
We denote by
$$
\cJ_{\lambda}^{\text{\rm sect}} = \cJ_{\lambda}^{\text{\rm sect}}(M)
$$
the set of $\lambda$-sectorial almost complex structures.
\end{defn}
In other words, $J$ is \emph{Liouville pseudoconvex} with respect to
$\psi = {\mathfrak s}_{\varphi,\kappa}$. (See Definition \ref{defn:liouville-pseudoconvex-intro}.)

\begin{theorem}
Let $M$ be a Liouville sector and fix a smoothing profile for $M$ (Definition \ref{defn:smoothing-profile-intro}).
The space $\cJ_{\lambda}^{\text{\rm sect}}$ is a
nonempty contractible infinite dimensional smooth manifold.
\end{theorem}

\begin{remark}
Here are a few notable points in our definition of $\cJ_{\lambda}^{\text{\rm sect}}$.
We assume that $M$ is a Liouville sector without corners for the moment.
\begin{itemize}
\item In contrast to~\cite{gps}, we do not demand the holomorphicity of the projection map
	\eqn\label{eqn.pi-to-C}
\pi_\C = R + iI: \nbhd^Z(\del \Mliou) \to \C
	\eqnd
near $\del \Mliou$ except near the ceiling corner $(F \setminus F_0 \times \C) \cap \nbhd_{\varepsilon_0}(\del M)$.
However, we do impose the condition that  the foliation $\cF_F$ mentioned before
is a $J$-complex foliation.
\item While we demand the common contact-type condition $\lambda \circ J = ds$ away from $\del M$, near the ceiling corner $\del_\infty M \cap \del M$, we utilize the choice of smoothing profile to ``interpolate'' the contact-type ceiling
    $\del_\infty M$ and to the flat-type wall $\del M$ to construct
    a pseudoconvex pair $(\psi, J)$ with $\psi$, especially on $\nbhd(\del_\infty M \cap \del M)$.
    Here we say that $\del M$ is `of flat-type' in that the characteristic foliation of $\del M$ forms a trivial fibration
$$
\xymatrix{\del M \ar[d]^{\pi_{\del M}} \ar[r]^\Psi & F \times \RR \ar[d]^{\pi_F}\\
\cN_{\del M} \ar[r]^\psi & F
}
$$
(See \cite[Theorem 1.4, Theorem 1.7]{oh:intrinsic}.)
\end{itemize}
\end{remark}
Note that by construction the pair $({\mathfrak s}_{\varphi,\kappa},J)$ with $\lambda$-sectorial $J$ is
a Liouville pseudoconvex pair.

\subsection{Sectorial Hamiltonians and their simplicial family}

Another important ingredient in defining the wrapped Fukaya category and symplectic cohomology
for Liouville sectors is to identify a suitable class of wrapping Hamiltonian functions whose Floer trajectories can only pass
through $\del \Mliou$ ``in the correct direction''. (See \cite[Section 2.10.1 \& Lemma 4.21]{gps}. Also see \cite{oh-floer-continuity,kasturirangan-oh} for an early discussion of such directedness
of Floer trajectories.) 

Our effort of constructing pseudoconvex pairs $({\frak s}_{\varphi,\kappa},J)$,
a wiggled end-profile function ${\frak s}_{\varphi,\kappa}$ and its adapted
$\lambda$-sectorial almost complex structures $J$ pays off since it makes \emph{very simple} the definitions of Hamiltonians and
nonnegative Hamiltonian isotopies or
more generally $n$-simplices thereof
that are amenable to the maximum principle, when paired with
$\lambda$-sectorial $J$'s. Furthermore the pair $(H,J)$ is also
amenable to the strong maximum principle
for the $Z$-invariant Lagrangian boundary conditions. We would like to
compare the simplicity and naturality of our definitions of
sectorial Hamiltonians (Definition~\ref{defn:sectorial-H}) and
of their simplicial family
\cite[Section 4]{oh-tanaka:liouville-bundles} with the complexity of the definition of
the dissipative $n$-simplices of Hamiltonians used
in \cite[Definition 4.5]{gps}: They need
the dissipative family to achieve the confinement results via
the monotonicity arguments. In the sectorial framework,
such complications are subsumed
in our careful geometric preparation of the description of background geometry of Liouville sectors related to smoothing of
corners of the Liouville sectors.
These are all structural results on the underlying geometry
of Liouville sectors (with corners) whose statements
have little to do with pseudoholomorphic curves, although they are motivated by our aim of achieving
the $C^0$-estimates via the maximum principle. (See Part I of the present paper.)

The following is the key identity that holds for any $\lambda$-sectorial almost complex structure,
which motivates our simple definition of sectorial Hamiltonians below.

\begin{prop}[Proposition \ref{prop:energy-identity}] Assume $J$ is a $\lambda$-sectorial almost complex structure.
Consider the function $H$ of the type
$$
H = \rho({\mathfrak s}_{\varphi,\kappa})
$$
for a smooth function $\rho: \RR \to \RR$. Then for any solution
$u: \R \times S^1  \to M$ (or $u: \R \times [0,1] \to M$ for the case with
boundary) of the equation
\eqn\label{eq:CRJH-intro}
(du - X_H(u) \otimes dt)^{(0,1)}_J = 0,
\eqnd
we have
\eqn\label{eq:Laplacian-intro}
\Delta({\mathfrak s}_{\varphi,\kappa} \circ u) =
\frac{1}{2}|du - X_H(u) \otimes dt|^2_J - \rho'({\mathfrak s}_{\varphi,\kappa}\circ u)\frac{\del}{\del \tau}({\mathfrak s}_{\varphi,\kappa}\circ u).
\eqnd
\end{prop}

\begin{defn}[Sectorial Hamiltonians] \label{defn:sectorial-Hamiltonians-intro}
Let $(M,\lambda)$ be a  Liouville sector and fix a smoothing profile
(Definition \ref{defn:smoothing-profile-intro}) for $M$.
Let ${\mathfrak s}_{\varphi,\kappa}$ be the associated end-profile function.
We call a Hamiltonian $H: M \to \RR$ \emph{sectorial} (with respect to the smoothing profile) if
$$
H = \rho({\mathfrak s}_{\varphi,\kappa})
$$
on a neighborhood $\nbhd(\del_\infty M \cup \del M)$ for some smooth function $\rho: \RR \to \RR_+$ with $\rho' > 0$.
\end{defn}

We will give the precise description of the above mentioned neighborhoods
$$
\nbhd(\del_\infty M\cup \del M)
$$
and $\nbhd(\del_\infty M\cap \del M)$
in Section \ref{subsec:sectorial-Hamiltonians}.

\begin{remark} Our usage of these smoothing profiles and sectorial packages replace the smoothing
operations omnipresent in \cite{gps} and other literature to study the K\"unneth-type theorems by absorbing
the smoothing process into the definition of sectorial package.
(See \cite{oancea}, \cite[Subsection 2.5]{gps}, \cite[Subsection 9.3.3]{gps-2} for example.)
See \cite[Conjecture 3.40 \& Conjecture 4.39]{gps} for the relevant conjectures.
We will study these conjectures elsewhere.
For this purpose, we have already identified a correct class of branes that is \emph{monoidal},
which we call \emph{gradient-sectorial Lagrangians} in \cite{oh:gradient-sectorial}.
\end{remark}

\subsection{Confinement theorems for Floer's equations}

The main utility of sectorial Hamiltonians and almost complex structures
is to establish confinement results for (Hamiltonian-perturbed) pseudoholomorphic curves
using the maximum principle, and the strong maximum principle relative to the $Z$-invariant
Lagrangian submanifolds on $\nbhd(\del_\infty M \cup \del M)$.

We establish confinement results all by the maximum principle with respect to
the $\lambda$-sectorial $J$ (and to the cylindrical Lagrangians for the strong maximum principle) in:
\begin{itemize}
    \item Theorem \ref{thm:autonomous-confinement} (for symplectic cohomology)
    \item Theorem \ref{thm:nonautonomous-confinement}  (for continuation via nonnegative Hamiltonian isotopies)
    \item Theorem \ref{thm:unwrapped} (for $A_\infty$ structure maps) and
    \item Theorem \ref{thm:closed-open} (for closed-open maps).
\end{itemize}
Once these confinement results are established, construction of wrapped Floer cohomology
\label{sec:construction-wrapped-Floer-cohomology}
and of covariantly functorial inclusion function $\Fuk(X) \to \Fuk(X')$ for the
inclusion of Liouville sectors $X \hookrightarrow X'$ follows the same procedure as
given in \cite[Section 3.4]{gps}, with the caveat that $J$ is defined only on $M \setminus \del M$.
This singularity-near-the-boundary introduces minor inconveniences in the proof
which are easily overcome by being careful about choices along the way.
(We elaborate this point in Section \ref{sec:inclusion-functor}.)

\subsection{Relation to other works}

In~\cite[2.10.1,2.10.3]{gps}, the authors opine it would be an ``important technical advance'' to identify a class of almost complex structures and Hamiltonians on Liouville sectors that guarantee the maximum principle and render the projection to $T^*[0,\infty)^k$ holomorphic. (See~\eqref{eqn.pi-to-C} and Remark \ref{rem:difficulty} for some difficulties in pinpointing a class.) As mentioned above, we accomplish this task, with the caveat that $J$ is defined only on $M \setminus \del M$.

The present author came to need this technical advance in \cite{oh-tanaka:liouville-bundles}, where the
authors thereof construct an unwrapped Floer theory
for bundles of Liouville manifolds, and made the (then-incorrect) claim that this theory extended to bundles of Liouville sectors. As a referee pointed out, one did not know of a good class of almost-complex structures on sectors suitable for maximum-principle-type arguments (as opposed to monotonicity arguments) at the time when the original version
of \cite{oh-tanaka:liouville-bundles} first appeared. By constructing a $\lambda$-sectorial
almost complex structure, we provide such a class which also pairs with the sectorial Hamiltonian $H$ in the way
that all confinement theorems for perturbed $(J,H)$-curves can be easily proved for $Z$-invariant Lagrangian branes
by the maximum principle. (See Part III of the present paper.)

\begin{remark}\label{rem:maximum-principle}
We now mention two important anticipated advantages of our sectorial package $(J,H)$ over
the dissipative package $(J,H)$ used in \cite{gps}:
\begin{enumerate}
\item
The proof of $C^0$-estimates via the maximum principle
will be important when considering bundles of Liouville sectors.
It is somewhat awkward to state the
relevant $C^0$ estimates by monotonicity arguments on bundles: This is because the \emph{total space}
of a family will be neither symplectic nor almost complex so
no monotonicity argument can be applied on the nose.
\item As mentioned before, our sectorial package provides a \emph{very simple} definitions
of Hamiltonians and nonnegative Hamiltonian isotopies or more generally $n$-simplices thereof
that are amenable to the maximum principle when paired with
$\lambda$-sectorial $J$'s the pair $(H,J)$  of which is also amenable to the strong maximum principle
for the $Z$-invariant Lagrangian boundary. We anticipate this simpleness
of the sectorial package will
facilitate the future simplicial study of wrapped Fukaya categories of Liouville sectors. See  \cite{gps-2}, \cite{asplund:simplicial} for such studies and
\cite{oh-tanaka:actions,oh-tanaka:liouville-bundles} for a
usage of the sectorial package for a simplicial study of Liouville bundles.
\end{enumerate}
\end{remark}

Finally, we would like to mention that there is a much simpler way of defining the sectorial Floer data
in the \emph{asymptotically-cylindrical-at-infinity (ACI) analytical framework} of Bao \cite{bao1,bao2}
(see also \cite{oh-wang2}) in terms of Giroux's notion of \emph{ideal Liouville domains and completions}.
We call them the \emph{ACI-sectorial Floer data}.
In this ACI framework, all ACI Floer data is defined
through Giroux's notion of ideal Liouville domains (or sectors) and then just restrict them to
$W \setminus \del_\infty W$. This restriction process becomes natural in the ACI framework but
not in the standard cylindrical-at-infinity (CI) analytical framework.
Nonetheless, mainly because the analysis of all existing literature --- except
\cite{bao1,bao2}, \cite{oh-wang2} as far as we are aware---is based on the CI analytical framework,
the present paper is written conforming to the CI framework.
We refer readers to Appendix \ref{sec:ACI-framework} for more detailed discussion on this.
We hope to elaborate this point of view elsewhere.

\medskip

{\bf Acknowledgments.}
We thank the unknown referee of the paper \cite{oh-tanaka:liouville-bundles} for pointing
out the difficulty of applying the maximum principle for Liouville sectors. This prompted the author to search for the present notion of sectorial almost complex structures.
We thank Hiro Tanaka for his collaboration on the study of Liouville geometry and for making
numerous useful comments on the early drafts of the present paper.
We also thank the referee for her/his careful reading and pointing out
many careless mistakes and writing of the paper. It helps us much
improve the presentation of Part I which is the most technical
part of the paper. We are especially grateful to the referee for kindly providing an elegant
idea of the proof of Lemma 3.13 which is
presented in Appendix with some modification and addition.
The author is supported by the IBS-project IBS-R003-D1.

\bigskip

\noindent{\bf Notations and terminologies:}
\begin{enumerate}
\item $\alpha$: eccentricity of a splitting data $\nbhd(\del M) \cong F \times \C_{\text{\rm Re} \geq 0}^k$,
\item $\varphi, \, \varphi_k$: a convex (corner) smoothing function,
\item $s_{k,\varphi}$: smoothing function of a sectorial corner $C_\delta$ of codimension $k$,
\item $s_{k+1,\varphi}$: a smoothing function of the ceiling of a sectorial corner $C_\delta$,
which we call an end profile function at a corner,
\item ${\mathfrak s}_\varphi$: end-profile function of $\nbhd(\del_\infty M \cup \del M)$,
\item ${\mathfrak s}_{\varphi,\kappa}$: wiggled end-profile function of $\nbhd(\del_\infty M \cup \del M)$,
\item $\varepsilon_0$: the size of a $Z$-invariant neighborhood $\nbhd^Z(\del M)$ of $\del M$ in
\eqref{eq:epsilon0}, which is fixed once and for all,
\item $C_\delta$: a sectorial corner,
\item $\del_\infty M \cap \del M$: the ceiling corner,
\item $C_\delta \cap \del_\infty M$: the ceiling of the sectorial corner $C_\delta$.
\item $\lambda$-sectorial almost complex structures: the main class of our interest in the present paper,
\item Sectorial almost complex structures: the main class of interest in \cite{oh:gradient-sectorial}.
\end{enumerate}

\bigskip

\noindent{\bf Locations of the definitions of key constants}
\begin{itemize}
\item $\varepsilon_0 > 0$; Notation \ref{nota:epsilon0}.
\item $\epsilon_1 > 0$; Condition \ref{cond:s}, Theorem \ref{thm:interpolation}
 and \eqref{eq:epsilon1}. It concerns the function $R$ and $s$,
 but mainly on $R$.
\item $N_0 > 0$;  Condition \ref{cond:s}. It concerns the function $I$.
\item $N_1 > 0$;  Definition \ref{defn:kappa}. It concerns the function $\kappa_2 = \kappa_2(I)$.
\item $N_2 > 0$; Condition \ref{cond:s} and Lemma \ref{lem:ceiling-corner}.
  It concerns the symplectization radial function $s$.
\end{itemize}

\section{Recollections and preliminaries}

\subsection{List of conventions}

In the present paper, we follow
the conventions of the book \cite{oh:book1} as summarized in \cite[p.xxi]{oh:book1} except
that there is no mention of symplectization there.

\begin{convention}
\begin{enumerate}
\item Definition of Hamiltonian vector field: $X_H \rfloor \omega = dH$,
\item The standard symplectic form on the cotangent bundle
$$
\omega_0 = - d\theta, \quad \theta = pdq,
$$
and hence the Liouville form $\lambda = -\theta$ and Liouville vector field is the Euler vector field
$\vec E = p\frac{\del}{\del p}$.
\item Symplectization: For given contact-type hypersurface $S_0 \subset M$, we have an identification
$$
\nbhd(\del_\infty M) = S_0 \times [0,\infty).
$$
\end{enumerate}
\end{convention}

The convention of \cite{gps} is the one which we call Entov-Polterovich's convention in
\cite{oh:book1}. Here are the corresponding list
of \cite{gps}:
\begin{enumerate}
\item [{(a)}] Definition of Hamiltonian vector field: $X_H \rfloor \omega = - dH$
\item [{(b)}] The standard symplectic form on the cotangent bundle: $\omega_0 = d\theta$,
and hence the Liouville form $\lambda = \theta$ and Liouville vector field is the Euler vector field
$\vec E = p\frac{\del}{\del p}$: The two negatives cancel each other.
\item [{(c)}] Symplectization: For given contact-type hypersurface $S_0 \subset M$, \cite{gps} has an identification
$$
\nbhd(\del_\infty M) = [0,\infty) \times S_0.
$$
(See \cite[p.9]{gps}.)
\end{enumerate}
\medskip
In the discussion of convexity of $\del_\infty M \cap \del M$,
these differences of the definitions cancel the differences when the meaning of
the \emph{outward pointing direction of $X_I$ along $\del X$} is determined. In particular
the outward pointing direction is that of
$$
- \frac{\del}{\del R}
$$
in both conventions. This gives rise to the identification
$$
\nbhd(\del M) \cong \del M \times [0, \infty) \cong F \times \C_{\text{\rm Re} \geq 0} \cong F \times T^*\RR_{\geq 0}
$$
in both cases. (See \cite[Section 2.1]{oh:intrinsic} for detailed discussion on this orientation issue.)

\subsection{Splitting data for boundaries}

\begin{notation}
Given $0<\alpha\leq1$, we write
\eqn\label{eq:lambda-alpha}
\lambda_\C^\alpha := (1-\alpha)xdy - \alpha y dx.
\eqnd
We also write $\omega_\C = dx \wedge dy$ and $J_\C$ for the standard symplectic and complex structures on $\C$ (with complex coordinate $z = x+i y$).
\end{notation}

We recall that, by definition, a Liouville sector allows us to split a neighborhood of $\del M$ as a direct product. We refer to~\cite{gps} for the notion of an $\alpha$-defining function, and begin with the case when $M$ only has boundaries (and no corners):

\begin{prop}[Proposition 2.25 \cite{gps}]\label{prop:gps}
Let $\Mliou$ be a Liouville sector and fix a real number $0<\alpha \leq 1$. Every $\alpha$-defining function
$I: \nbhd^Z(\del \Mliou)\to \RR$ extends to a diffeomorphism (valid over a cylindrical neighborhood of the
respective boundaries) $\nbhd^Z(\del\Mliou) \cong F \times \C_{{\text{\rm Re}}\geq 0}$
in which $I = y\circ \pi_\C$ is the imaginary part of the $\C_{{\text{\rm Re}} \geq 0}$-projection and $F$ is a Liouville manifold.

This diffeomorphism also gives an identification
\eqn\label{eq:splitting}
(\nbhd^Z(\del\Mliou),\lambda_\Mliou|_{\nbhd^Z(\del\Mliou)})
\cong
(F \times \C_{{\text{\rm Re}}\geq 0}, \pi_F^*\lambda_F + \pi_\C^* \lambda_\C^\alpha + df)
\eqnd
where $\lambda_F$ is the Liouville form on $F$, and $f: F \times \C_{{\text{\rm Re}} \geq 0} \to \RR$ satisfies the following properties:

\newenvironment{splitting-f-props}{
	  \renewcommand*{\theenumi}{(f\arabic{enumi})}
	  \renewcommand*{\labelenumi}{(f\arabic{enumi})}
	  \enumerate
	}{
	  \endenumerate
}
\begin{splitting-f-props}
\item\label{item. f supported in F0 x C} $f$ is supported inside $F_0 \times \C$ for some
compact Liouville domain $F_0 \subset F$.
\item $f$ is independent of the $\C_{{\text{\rm Re}} \geq 0}$-coordinate for $|I|$ sufficiently large.
\end{splitting-f-props}
\end{prop}

\begin{notation}[$R$]
In the course of the proof of Proposition \ref{prop:gps}, \cite{gps} defines another function
$R: \nbhd^Z(\del \Mliou) \to \RR$. It satisfies
\eqn\label{eq:ZRZI}
	R|_{\del X} = 0,
	\qquad ZR = (1-\alpha)R,
	\qquad  \omega(X_R,X_I) = 1.
\eqnd
\end{notation}

Since the choice of the pair $(F,(R,I))$, and the constant $\alpha$, will appear frequently, we give them names.

\begin{defn}[Splitting data and eccentricity]\label{defn:eccentricity} Let $(\Mliou,\lambda)$ be a Liouville sector.
We call a pair
$(F,(R,I))$ a \emph{splitting data} of $(\Mliou,\lambda)$, and the associated $\alpha$ the
\emph{eccentricity} of the splitting.
\end{defn}

\begin{notation}
As already insinuated in Proposition~\ref{prop:gps}, we denote the projections by
$$
\pi_F: \nbhd^Z(\del \Mliou) \to F, \quad \pi_\C: \nbhd^Z(\del \Mliou) \to \C.
$$
\end{notation}

\begin{remark}
The map $\pi_\C$ can be written as
$$
\pi_\C = R+ iI : \nbhd^Z(\del \Mliou) \to \C_{{\text{\rm Re}} \geq 0}
$$
with $F = (R,I)^{-1}(0,0)$
\end{remark}

\begin{notation}
We let $Z_\C^\alpha$ be the Liouville vector field of $\lambda_\C^\alpha$
on $\C$, so that
$$
Z_\C^\alpha = (1-\alpha) x\frac{\del}{\del x} + \alpha y \frac{\del}{\del y}.
$$
Then, near $\del M$, we have the decomposition
\eqn\label{eq:Z}
Z = Z_F \oplus Z_\C^\alpha -X_f
\eqnd
where $Z_F$ is the Liouville vector field on $F$.
\end{notation}

\begin{notation}
We note that $\lambda_\C^\alpha \circ J_\C$ is an exact one-form. A choice of primitive is
\eqn\label{eq:hCCalpha}
h_\C^\alpha(x,y) = \frac12\left((1-\alpha)x^2 + \alpha y^2\right).
\eqnd
\end{notation}
\begin{remark}\label{rem:eccentricity}
The eccentricity of the splitting data measures how much the level set $\{h_\C^\alpha = 2\}$---which is an ellipse for $0 < \alpha < 1$---deviates from the circle. The associate level set degenerates to the $y$-axis
as $\alpha \to 0$ while it degenerates to the $x$-axis as $\alpha \to 1$ and it becomes a circle when
$\alpha = \frac12$. Indeed the precise \emph{geometric eccentricity} of the relevant ellipse is given by
$$
\sqrt{1- \left(\frac{1-\alpha}{\alpha}\right)^2}
\qquad
\left(\text{resp., }\,
\sqrt{1- \left(\frac{\alpha}{1-\alpha}\right)^2}\right)
$$
when $\alpha > \frac12$ (resp., $\alpha < \frac12$).
\end{remark}

\begin{notation}[$\lambda_\C$ (without the $\alpha$)]
We will often denote $\lambda_\C := \lambda_\C^\alpha$ and likewise for $h_\C$,
dropping the superscript $\alpha$ when explicit mention of it is unnecessary.
\end{notation}

\subsection{Splitting data for sectorial corners}
The above can be generalized to the case of with corners.

\begin{lemma}[Lemma 9.7 \& Lemma 9.8 \cite{gps-2}]\label{lem:splitting-corner}
Let $H_1, \ldots, H_n \subset M$ be a sectorial collection of hypersurfaces (\cite[Definition~9.2]{gps-2}), equipped with choices of
$I_i: \nbhd^ZH_i \to \RR$. Then the following holds:
\begin{enumerate}
\item  There exist unique functions $R_i: \nbhd^ZH_i \to \RR$ satisfying
$Z[R_i] = 0$ near infinity, $H_i = \{R_i = 0\}, \, \{R_i, I_j\} = \delta_{ij}$ and $\{R_i,R_j\} = 0$.
\item The map
$$
\pi_{T^*\RR^k}: = ((R_{i_1}, I_{i_1}), \ldots, (R_{i_k}, I_{i_k})): \nbhd^Z(H_{i_1} \cap \cdots \cap H_{i_k}) \to
T^*\RR^k
$$
is a symplectic fibration whose symplectic connection is flat, thus giving a symplectic product
decomposition
	\eqnn
	M \supset \nbhd^Z(H_{i_1} \cap \cdots \cap H_{i_k}) \cong F \times T^*\RR^k
	\eqnd
identifying $Z$-invariant neighborhoods of
$\nbhd^Z(H_{i_1} \cap \cdots \cap H_{i_k})$ and $F \times T^*\RR^k$.
\end{enumerate}
\end{lemma}

\begin{defn}\label{defn:sectorial-corner}
Following~\cite{gps-2}, if $M$ is a Liouville manifold-with-corners whose boundary forms a sectorial collection when viewed as a collection of hypersurfaces, we say $M$ is a Liouville sector with corners.
\end{defn}

\begin{defn}[Splitting data for corners]\label{defn:splitting-corner} Let $(M,\lambda)$ be a Liouville sector with corners and
consider a corner of $\del M$ with codimension $k-1$ with $k\geq 1$. Let $H_{i_1}, \ldots, H_{i_k}$
be the associated faces of the corner. We call the collection of data $(F,\{(R_{i_j}, I_{i_j})\})$
as in Lemma \ref{lem:splitting-corner} a \emph{splitting data} of the corner.
\end{defn}

\begin{remark}
The Liouville isomorphism type of $F$ of course depends on the connected component of the corner stratum, but we suppress this from the notation.
\end{remark}

\subsection{Neighborhoods of $\del_\infty M$ and of $\del M$}

We start with the following using the definition of Liouville sector.
\begin{notation}\label{nota:epsilon0} Let $M$ be a Liouville sector (without corners).
We fix a decomposition
\eqn\label{eq:epsilon0}
\Mliou = (\Mliou \setminus F \times \C_{0 \leq \textrm{\rm Re} < \varepsilon_0}) \bigcup
(F \times \C_{0 \leq \textrm{\rm Re} < 2 \varepsilon_0})
\eqnd
for some $\varepsilon_0$ which will be fixed once and for all.
\end{notation}

For given Liouville manifold $(M,\lambda)$, there is a canonical projection
$$
\nbhd(\del_\infty M) \to \del_\infty M
$$
of sending a point to the Liouville ray issued at the point. This map also induces a canonical
contact structure on $\del_\infty M$, but not the contact form which is defined only up to contactomorphism.
(See \cite{giroux} e.g., or Appendix \ref{sec:ACI-framework} for a brief discussion.)
\begin{choice}[$s,S_0$]
\label{choice. s S_0}
We choose a contact-type hypersurface $S_0 \stackrel{\iota}\subset M$ and its associated symplectization radial function
	\eqn\label{eqn. s}
	s: \nbhd(\del_\infty M) \to \RR_{\geq 0}
	\eqnd
is such that $S_0 = s^{-1}(0)$ defined on a neighborhood $\nbhd(\del_\infty M)$ and satisfying $Z[s] = 1$.
A choice of $S_0$ defines a Liouville embedding of
$$
(S_0 \times [0,\infty), d(e^s \iota^*\lambda)) \to (M,d\lambda)
$$
via the Liouville flow of $Z$, if we start with the contact-type hypersurface $S_0$ sufficiently at infinity.
In terms of~\eqref{eqn. s}, any real number $t\geq 0$ exhibits a contactomorphism $\del_\infty M \cong s^{-1}(t)$. This provides a local fibration
$$
s: \nbhd(\del_\infty M) \to [0,\infty)
$$
through the identification $\nbhd_{s \geq 0}(\del_\infty M) \cong  S_0 \times [0,\infty)$.
We call $s$ the (symplectization) radial function associated to the contact-type hypersurface,
and mention that $\{s \leq N\}$ is compact for all $N \geq 0$ by definition of the
symplectization end.

\end{choice}
\begin{remark}
By convexity of $\del_\infty \Mliou \cap \del \Mliou$ in $\del_\infty \Mliou$ and
since $Z$ is tangent to $\del \Mliou$ outside a compact set of $\Mliou$,
we have an embedding
$$
\left(S_0 \times [0, \infty), \del S_0 \times [0, \infty)\right) \to (\Mliou, \del \Mliou).
$$
This induces an identification
$$
\del \Mliou \cap \{s \geq 0\} \cong C_0 \times [0, \infty),
$$
where $C_0 = S_0 \cap \del M$. Note that Liouville flow also naturally induces a diffeomorphism
$$
C_0 \cong \del M \cap \del_\infty M
$$
by restricting the canonical projection $\nbhd(\del_\infty M) \to \del_\infty M$ to
$S_0 \subset \nbhd(\del_\infty M)$.
\end{remark}

\begin{notation}
 Unless otherwise said, we denote by
$$
\nbhd(\del_\infty \Mliou) = \nbhd_{s \geq 0}(\del_\infty \Mliou)
$$
the image of the above mentioned embedding.
We may assume that $s$ is defined on a slightly
bigger region than $\nbhd(\del_\infty \Mliou)$, meaning that we may from time to time extend the function $s$ from the ray $[0, \infty)$ to
$[-\delta, \infty)$.
\end{notation}
%

Mostly we will work with the Liouville sectors without sectorial corner for the actual constructions of
the seeked almost complex structures and indicate how to extend to the case of with
sectorial corners.

The following choice of the radial function will be used later in the proof of
Proposition \ref{prop:linear-independence}.

We fix the radial function $s_\C$ on $\C_{\text{\rm Re} \geq 0}$  to be
\be\label{eq:sC}
s_\C = \log h_\C = \log\left(\frac12 ((1-\alpha) R^2 + \alpha I^2)\right)
\ee
\begin{cond}\label{cond:s}
Let $N_0> 0$ and $F_0 \subset F$ be so large that
\be\label{eq:suppf}
\supp f \subset F_0 \times \C_{\text{\rm Re}\geq 0}
\ee
and $d_\C f = 0$ on $\{|I| \geq N_0\}$, which holds true by Proposition \ref{prop:gps}.
\begin{enumerate}
\item
We choose the $S_0$ so that the associated radial function $s$ has the following
properties:
\be\label{eq:choice-s}
s = \beta \pi_F^*s_F + (1-\beta) \pi^*s_\C
\ee
on
 $$
( F \setminus F_0)  \times \C_{\text{\rm Re}\geq 0}
 \bigcup
 F_0 \times \{(R,I) \mid |I| \geq N_0
 \}.
 $$
 \item We choose a sufficiently large $N_2 > 0$ and a small $\epsilon_1 > 0$
  such that
 \bea\label{eq:choice-N2}
 &{}& F_0 \times \{(R,I) \mid |I| \leq N_0, \, 0 \leq  R \leq 2 \varepsilon_0 \}
 \nonumber\\
 & \subset & M \setminus
( \nbhd_{s \geq N_2}(\del_\infty M) \cup \nbhd_{R \leq \epsilon_1}(\del M)).
 \eea
 \end{enumerate}
\end{cond}

\begin{lemma}
The choice of contact-type hypersurface required in
Condition \ref{cond:s}  can be made.
\end{lemma}
\begin{proof} The requirement (2) follows since the subset
$$
F_0 \times \{(R,I) \mid |I| \leq N_0, \, 0 \leq R \leq 2\varepsilon_0\}
$$
is compact.

For (1), we note that
$\lambda = \lambda_F + \lambda_\C$ and so
$$
Z = Z_F + Z_\C
$$
on $(F \setminus F_0) \times \C_{\text{\rm Re}\geq 0}
\subset M \setminus \supp f $. Therefore on $(F \setminus F_0) \times \C_{\text{\rm Re}\geq 0}$,
$$
Z(\beta \pi_F^*s_F + (1-\beta) \pi^*s_\C) = \beta Z_F[s_F]
+ (1-\beta) Z_\C[s_C] = 1.
$$

On the other hand, on the region  $F_0 \times \{(R,I) \mid |I| \geq N_0\}$,
we compute
\beastar
Z[\beta \pi_F^*s_F + (1-\beta) \pi^*s_\C]
& = & (Z_F \oplus Z_\C - X_f)[\beta \pi_F^*s_F + (1-\beta) \pi_\C^*s_\C]\\
& = & 1 - X_f[\beta \pi_F^*s_F ].
\eeastar
Since $X_f$ is bounded (on $\nbhd_{0 \leq R \leq 2\varepsilon_0}(\del M)$)
$$
|X_f(y,R,I)| \leq C_{f,N_0}
$$
for a constant $C = C_{f,N_0} > 0$. Therefore we can choose $\beta> 0$ so small that
\be\label{eq:choice-beta}
0 < \beta < \frac1{C_{f,N_0}}.
\ee
Then $Z$ is transversal to
$(\beta \pi_F^* s_F + (1-\beta) \pi_\C^* s_\C)^{-1}(0)$
on $F \times \C_{0 \leq \text{\rm Re}\leq 2\varepsilon_0} \subset M$. Hence we can
take a contact-type hypersurface
$$
S_0^{\text{\rm loc}}: =  (\beta \pi_F^* s_F + (1-\beta) \pi_\C^* s_\C)^{-1}(0)
$$
thereon. Now we would like to extend $S_0^{\text{\rm loc}}$ to the outside of
$\{0 \leq R \leq 2\varepsilon_0\}$ to define
a (global) contact-type hypersurface on $M$ which we denote
by $S_0$.

For this purpose, let $s^{\text{\rm ini}}$ be the radial function associated
to a contact type hypersurface $S_0^{\text{\rm ini}}$
 which exists by definition of Liouville sectors. We consider the
associated symplectization end
$$
S_0^{\text{\rm ini}} \times [0,\infty) \hookrightarrow M
$$
and write $(y, s^{\text{\rm ini}})$ for an element of the product and denote by
$$
\pi: S_0^{\text{\rm ini}} \times [0,\infty) \to S_0^{\text{\rm ini}}
$$
the associated projection. Then we have
$S^{\text{\rm ini}}_0 = \{s^{\text{\rm ini}} = 0\}$.
By the above established
transversality, we can write $S_0^{\text{\rm loc}}$ as the graph
$$
\Graph g = \{(y,g(y))\mid y \in \pi(S_0^{\text{\rm loc}})\}
$$
of some function
$g:  \pi(S_0^{\text{\rm loc}}) \subset S_0^{\text{\rm ini}} \to \R$.
Since $Z[g] > 0$, we can
extend the function to $M$ by setting
$$
s := (1-\chi) s^{\text{\rm ini}}
$$
for a suitable partition of unity $\{\chi(R),1-\chi(R)\}$ such that
 $\chi\equiv 0$ near $\del M$ and $\chi'(R) \leq 0$ with strict
 inequality away from $\{\chi(R) =1,\, 0\}$. Then we have
\beastar
Z[s] & = & -\chi'(R)Z[R]  s^{\text{\rm ini}} + (1-\chi(R)) Z[s^{\text{\rm ini}}]\\
& = & - \chi'(R)(1-\alpha) s^{\text{\rm ini}} + (1-\chi(R)) > 0
\eeastar
say, for $s^{\text{\rm ini}} > 1$. (The latter holds essentially because
the set of functions $t$ satisfying $Z[t] > 0$ is convex.)
This finishes the proof.
\end{proof}

Then with respect to this $S_0$,
we take a new symplectization embedding
$$
S_0 \times [0,\infty) \hookrightarrow M
$$
and define the corresponding radial function which we denote by $s$
in the rest of the paper.

\clearpage

\part{Geometry of Liouville sectors}

In this part, our goal is to unravel various geometry of Liouville sectors which will
enter our construction of relevant almost complex structures.

The first thereof is to construct a family of exhaustion functions which smoothen the corners
$\del M \cap \del_\infty M$. For this purpose, we use two stage-construction
(See Section \ref{sec:intrinsic-geometry} for more detailed explanation.):
\begin{itemize}
\item We utilize a universal family of corner-smoothing functions of the standard
corner structure of $\RR^k_+$ which is compatible for different $k$'s which is
constructed in \cite[Section 18.5]{fooo:book-virtual}. However we need to augment the \emph{convexity}
property to this family for our purpose of constructing almost complex structures that are
amenable to the maximum principle.
\item We then utilize another natural geometric construction of symplectic $\RR^k$-action
induced by the definition of \emph{sectorial collection} from \cite{gps-2} as exposed in
\cite{oh:intrinsic}
\end{itemize}
This will provide us with a canonical family of corner-smoothing function of
Liouville sectors with corners.

The second ingredient is the study of \emph{Liouville deformation} $\lambda_t$ of
Liouville forms $\lambda$. Such a deformation will be needed
to ensure that the strong maximum principle is applicable to $J$-holomorphic
curves with \emph{$Z$-invariant-at-infinity Lagrangian} boundary condition for \emph{$\lambda$-sectorial
almost complex structures}, whose construction is of our main interest in the present paper.

We start with the first step of constructing the family of convex functions $\varphi = \varphi_k$
on $\RR^k_+$ for varying $k$'s.

\section{Convex corner-smoothing functions of $\RR^k_+$}
\label{sec:corner-smoothing-functions}

We first construct a suitable collection of convex functions $\varphi: (0,\infty)^n \to \RR$ which accomplish this goal on the standard Euclidean corner $[0,\infty)^n$. That is, we will find $\varphi$ so that $\varphi^{-1}(N)$ is a smoothing of the Euclidean corner for large enough $N$, and so that $\varphi^{-1}(N)$ is suitably convex. We will also construct $\varphi$ exhibiting various symmetries; this will allow us to turn these local models of smoothings to global smoothings of sectorial corners.

These functions---which we label by $\varphi$, and by $\varphi_k$ when we want to make the codimension $k$ explicit---are part of the data needed to define the functions ${\mathfrak s}_\varphi$ from the introduction.


\subsection{Octants}

\begin{notation}[Symmetric hyperplanes]
Consider the hyperplane
\eqn\label{eq:symmetric-R2}
\RR^n_{\bf y} := \{(x_0,\ldots, x_n) \in \RR^{n+1} \mid x_0 + \ldots + x_n =0\}.
\eqnd
We denote an element of this hypersurface by ${\bf y}$.
The permutation group ${\mathsf S}_{n+1}$ on $n+1$ letters acts on $\RR^n_{\bf y}$ by
permuting the coordinate variables $(x_0, \cdots, x_n)$.
\end{notation}

\begin{notation}[Octants]
We denote the positive octant by
$$
{\mathbf Q}^{n+1}_+: =  [0,\infty)^{n+1} \subset \RR^{n+1}.
$$
\end{notation}

Noting that $\RR^n_{\bf y} \subset \RR^{n+1}$ is preserved by the ${\mathsf S}_{n+1}$-action,
it also carries the induced action thereof. We define an ${\mathsf S}_{n+1}$-action so that
the action on the factor $\RR_{\bf y}^n$ is the induced one and on the $\RR$-factor is the trivial one
on the product $\RR^n_{\bf y} \times \RR$.
We now define an explicit ${\mathsf S}_{n+1}$-equivariant orthogonal linear isomorphism
\eqn\label{eq:Pi}
\Pi: \RR^{n+1} \to \RR^n_{\bf y} \times \RR \cong \RR^{n+1}
\eqnd
to $\Pi = (\pi_{\vec 1};t)$ which is given by
\begin{eqnarray}\label{eq:xi}
{\bf y} & = & \pi_{\vec 1}(x) \\
t & = & \frac{x_0 + \cdots + x_n}{\sqrt{n}}
\end{eqnarray}
where the map
$$
\pi_{\vec 1}: \RR^{n+1} \to \RR^n_{\bf y}
$$
is the orthogonal projection along the vector $\vec 1$. The map $\pi_{\vec 1}$  restricts to a homeomorphism
$$
\pi_{\vec 1}|_{\del{\mathbf Q}^{n+1}_+}:\del{\mathbf Q}^{n+1}_+ \to \RR^n_{\bf y}; \quad \vec 1 = (1,\cdots,1).
$$
Let $y= (y_1, \ldots, y_n)$ be the standard coordinates of $\RR^n_{\bf y} \cong \RR^n$ and $t$ the one on $\RR$.
Then we can write
\eqnn
\Pi^{-1}(y_1, \ldots,y_n; t) = (x_0(y),\ldots,x_n(y)) + t \vec 1 = (x_0(y)+t,\ldots, x_n(y)+t).
\eqnd
The map $\Pi$ restricts to an ${\mathsf S}_{n+1}$-equivariant homeomorphism
	\eqn\label{eqn. Phi}
\Phi: {\bf Q}_+^{n+1} \to \RR^n_{\bf y} \times [0,\infty)
	\eqnd
given by
$
\Phi(x): = ({\bf y}(x); t(x))
$
where ${\bf y}(x) = \pi_{\vec 1}(x)$ and $t(x) \geq 0$ is the constant determined by
$$
x - \pi_{\vec 1}(x) = t(x) \vec 1.
$$
This map $\Phi$ is smooth away from the corners.

\begin{remark}\label{rem:fooo-smoothing}
We refer readers to \cite[Section 18.5]{fooo:book-virtual} for some
systematic discussion on the corner smoothing of $[0,\infty)^n$
relevant to the smoothing process we apply below. One thing that was not
addressed in \cite[Section 18.5]{fooo:book-virtual} is the convexity property
of the smoothing which is needed for our present work.
\end{remark}

\subsection{Symmetric convex corner smoothing of $[0,\infty)^2$}
For the sake of exposition, we begin with the two-dimensional case
for which we write every map mentioned above explicit. For this purpose, we will
go back and forth between two coordinate systems $(x_0,x_1)$ in the description of
the gluing process because the coordinate system $(y;t)$ is more convenient for the visualization and
the coordinate system $(x_0,x_1)$ is convenient for calculations.

First, we describe the coordinate system $(y;t)$ explicitly.
By considering the orthonormal basis change of $\RR^2$ from the standard one to $\{e_{\bf y}, e_t\}$
given by
\begin{eqnarray}\label{eq:basis-change}
e_{\bf y} & = & \frac1{\sqrt{2}} \left(\frac{\del}{\del x_1} - \frac{\del}{\del x_2}\right),
\nonumber\\
e_t & = & \frac1{\sqrt{2}} \left(\frac{\del}{\del x_1} + \frac{\del}{\del x_2}\right).
\end{eqnarray}
We equip $\RR^2 = \RR_{\bf y} \times \RR$ with the coordinates $(y,t)$ associated to this basis.
We denote it non-symmetrically as $(y;t)$ to emphasize $y$ as the principle coordinate of our interest.

In terms of the coordinates of $(x_0,x_1)$ of $\RR^2$ and $(y;t)$ of $\RR_{\bf y} \times \RR$ the map
$$
\Pi: \RR^2 \to \RR_{\bf y} \times \RR
$$
is just given by the corresponding coordinate change map
\begin{eqnarray}\label{eq:yt}
y & = & \pi_{\vec 1}(x_0,x_1) = \frac1{\sqrt{2}}(x_1-x_2), \nonumber\\
t & = & \frac1{\sqrt{2}}(x_1 + x_2).
\end{eqnarray}
Its inverse is given by
$$
\Pi^{-1}(y;t) = \left(\frac{t+y}{\sqrt{2}},\frac{t-y}{\sqrt{2}} \right) = (x_0(y;t),x_1(y;t)).
$$

\begin{figure}[h]
\centering
\includegraphics[scale=0.5]{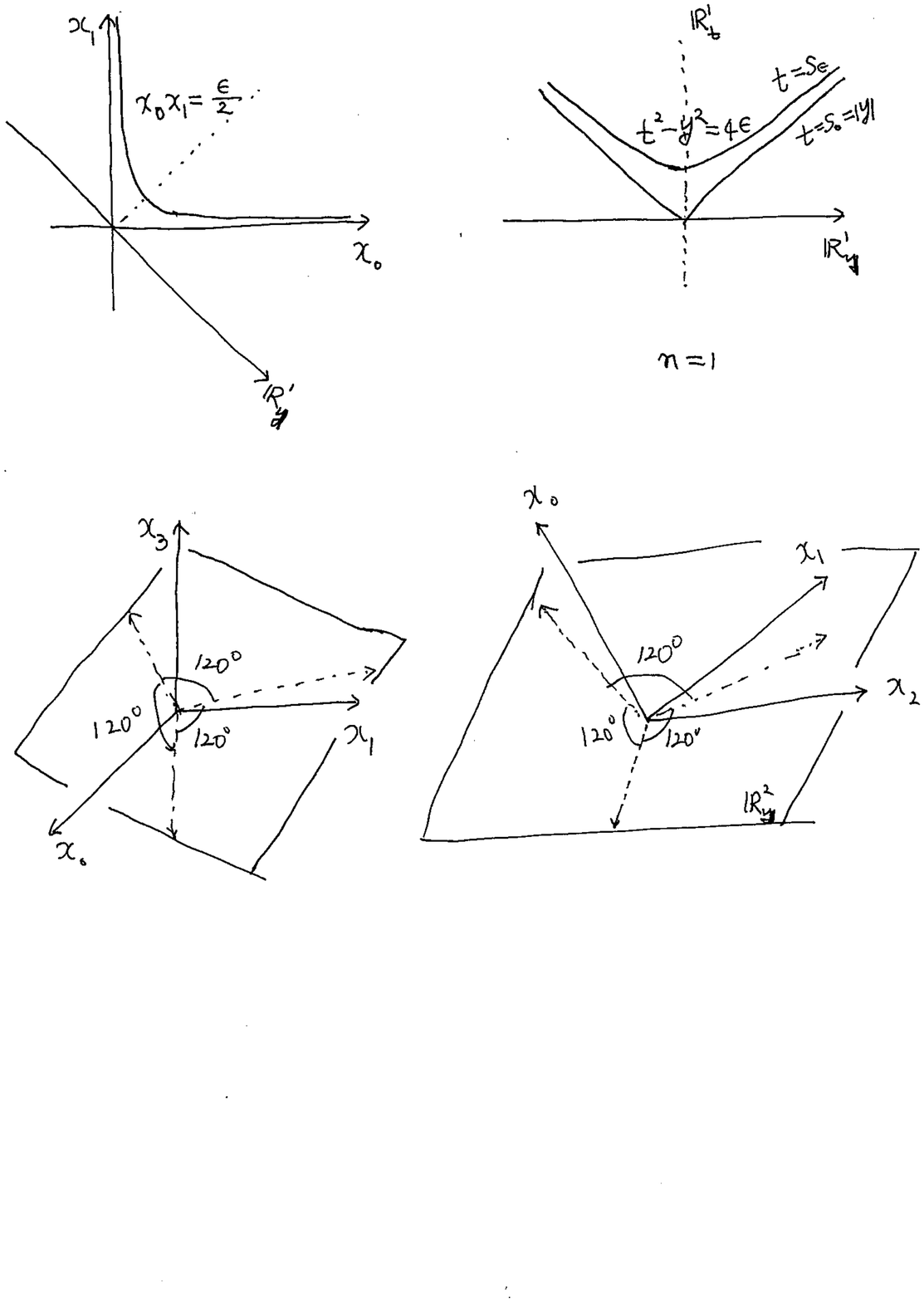}
\caption{Sections $s_\epsilon$ and $s_0$}
\label{fig:section-s0}
\end{figure}

Next we consider a one-parameter family of hyperbolas
\eqn\label{eq:S-epsilon}
S_{\epsilon} = \left \{(x_0,x_1) \in {\bf Q}^2_+ \, \big\vert \,  x_0 x_1 = \frac{\epsilon}{2}\right\} \subset \RR^2,
\qquad
\epsilon>0
\eqnd
as a smoothing of the corner of $\del {\bf Q}_+^2$.
This is a family of curves that has the following properties:
\begin{enumerate}
\item [{(i)}] they are convex in the direction towards $-\infty \cdot (1,1)$,
\item [{(ii)}] they are asymptotic to $\del {\mathbf Q}^2_+$ as $\epsilon \to 0$, and
\item  [{(iii)}] they are symmetric under reflection about the line $x_0 = x_1$ (i.e., about the $t$-axis).
\end{enumerate}
This motivates us to contemplate the orthogonal projection along the vector $(1,1)$
\eqnn
\pi_{(1,1)}: \RR^2 \to \RR^1_{{\bf y}}.
\eqnd
onto the line orthogonal to $(1,1)$. The restrictions
	\eqnn
	\pi_{(1,1)}|_{S_\epsilon} : S_\epsilon \xrightarrow{\cong} \RR^1_{{\bf y}}
	\xleftarrow{\cong} \del {\mathbf Q}^2_+ : \pi_{(1,1)}|_{\del {\mathbf Q}^2_+}
	\eqnd
exhibit a homeomorphism between $S_\epsilon$ and $\del {\mathbf Q}^2_+$ for all $\epsilon > 0$.
We regard
$$
\pi_{(1,1)}: \RR^2 \to \RR_{\bf y}
$$
as a trivial fiber bundle over $\RR_{\bf y}$ with a fiber $\RR$. (See Figure \ref{fig:section-s0}.) Then we express the
image of a section $s$ of this vector bundle in terms of the fiber coordinate $t$
$$
s({\bf y}): = t(x_0({\bf y}),x_1({\bf y})).
$$

Then in terms of the coordinate $y$ of $\RR_{\bf y}$, the following is obvious.

\begin{lemma} Let $s_\epsilon: \RR^1_{{\bf y}} \to \RR^2$ be the family of smooth sections
whose image is given by $S_\epsilon$ for $\epsilon \in (0,\infty)$, and
let $s_0$ be the piecewise linear limit $s_0$ whose image is given by
$\del {\mathbf Q}^2_+$. Then for $\epsilon > 0$
\eqn\label{eq:derivative-s-epsilon}
0< |s'_\epsilon(y)| < 1
\eqnd
for all $y$ and $s_0'(y) =\pm 1$ for $\pm y> 0$ respectively.
\end{lemma}
\begin{proof} Note that in $(y,t)$ coordinates of $\RR^2$, $S_\epsilon$ is expressed as
\eqn\label{eq:Sepsilon-inyt}
t^2 - y^2 = 4 \epsilon.
\eqnd
Therefore we have the formula
$$
s_0(y) = |y|, \quad s_\epsilon(y) = \sqrt{y^2 + 4 \epsilon}.
$$
Then straightforward calculation finishes the proof.
\end{proof}

\begin{remark}
We have a natural decomposition of $\RR^1_{\bf y}$ into $2$ rays
    $$
    \Delta_i: = \pi_{(1,1)}(\{x \in {\mathbf Q}^2_+ \mid x_i = 0\}) = \{(x_0,x_1) \mid x_0 + x_1 = 0, \, x_i> 0 \},
    \qquad
    i = 0,1,
    $$
which are naturally identified with the rays
$$
{\mathbf Q}^1_i: = \del {\mathbf Q}^2_+ \cap \{x_i = 0\} \subset \RR^2,
\qquad
i = 0,1,
$$
via the projection $\pi_{(1,1)}$.  Denote by $\tau : \RR^2 \to \RR^2$ the reflection
along the line $x_1 = x_0$ so that
$$
\tau(\Delta_i) = \Delta_{i+1} \quad \text{ for }\, i = 0, \, 1\mod 2.
$$
Note that $S_\epsilon$ is preserved as a set under the action by $\tau$ and
$$
\inf \left\{r = \sqrt{x_0^2 + x_1^2}\, \Big|\,  x \in S_{\epsilon}, \, \frac{\sqrt{\epsilon}}{4}
\leq x_0 \leq 2 \sqrt{\epsilon} \right\} = \sqrt{\epsilon}.
$$
\end{remark}

We now define an ${\mathsf S}_2$-equivariant smooth gluing, denoted by
$$
\widetilde S_\epsilon,
$$
by interpolating the curves
$$
S_\epsilon \cap D^2(2\sqrt{\epsilon}))
$$
and
$$
\del {\mathbf Q}_+^{2} \setminus D^2(2 T_0 \sqrt{\epsilon})
$$
along the region $D^2(2 T_0 \sqrt{\epsilon}) \setminus D^2(\sqrt{2\epsilon}))$.

We will do this expressing them as the graphs of functions in the $(x_0,x_1)$-coordinates as follows.
Recall that the hyperbolas $S_\epsilon$ are graphs of the function
$$
f_\epsilon(x_0) = \frac{\epsilon}{2x_0}
$$
for $\epsilon > 0$. We set $f_0$ to be the zero function.
First for the region $x_0 \geq \sqrt{\epsilon/2}$, we will smoothly interpolate them along the interval
$$
2 \sqrt{\epsilon} < x_0 < 2 T_0 \sqrt{\epsilon}
$$
for some choice of a small $\epsilon$ and a large $T_0> 0$, and then take the reflection
thereof along the line $x_1 = x_0$ to obtain
the other half of the interpolated curve over the region $0 \leq x_0 \leq \sqrt{\epsilon/2}$
near $\del {\bf Q}^2_+$.

\begin{choice}[$\epsilon$, $T_0$ and $\varepsilon_0$] \label{choice:Choice-T0}
The precise choices of $\epsilon$ and $T_0$
will be made later in Section \ref{sec:sectorial-J} so that---in terms of $\varepsilon_0$, another constant defined there---we will have
$$
0 < \varepsilon_0 \ll 2T_0 \sqrt{\epsilon} \ll \frac32 \varepsilon_0.
$$
We will take a sufficiently small $\epsilon$ and a sufficiently large constant $T_0 > 1$ relative to $\epsilon$.
\end{choice}

We want the graphs of $\widetilde f_\epsilon$ to be convex as follows. We postpone its proof
till Appendix.

\begin{prop}\label{prop:tildefepsilon}
 Let $0 < \epsilon < \frac12$ be given. Then there exists a sufficiently large $T_0 > 0$ independent of
  $\epsilon$ and a function $\widetilde f_\epsilon$ such that
\begin{enumerate}
\item $\widetilde f_\epsilon(x) = f_\epsilon(x)$ for  $\sqrt{\epsilon/2} \leq x \leq 2 \sqrt{\epsilon}$,
\item $-1 < \widetilde f_\epsilon'(x) \leq 0$ for all $x \geq \sqrt{\e/2}$,
\item $\widetilde f_\epsilon''(x) \geq 0$ for all $x$, and
\item $\widetilde f_\epsilon''(x) > 0$ for $x$ with $2 \sqrt{\epsilon}< x < 2T_0 \sqrt{\epsilon}$.
\item $f_\epsilon + x f_\epsilon'(x) \geq 0$ for all $x$.
\end{enumerate}
\end{prop}

\begin{figure}[h]
\centering
\includegraphics[scale=0.5]{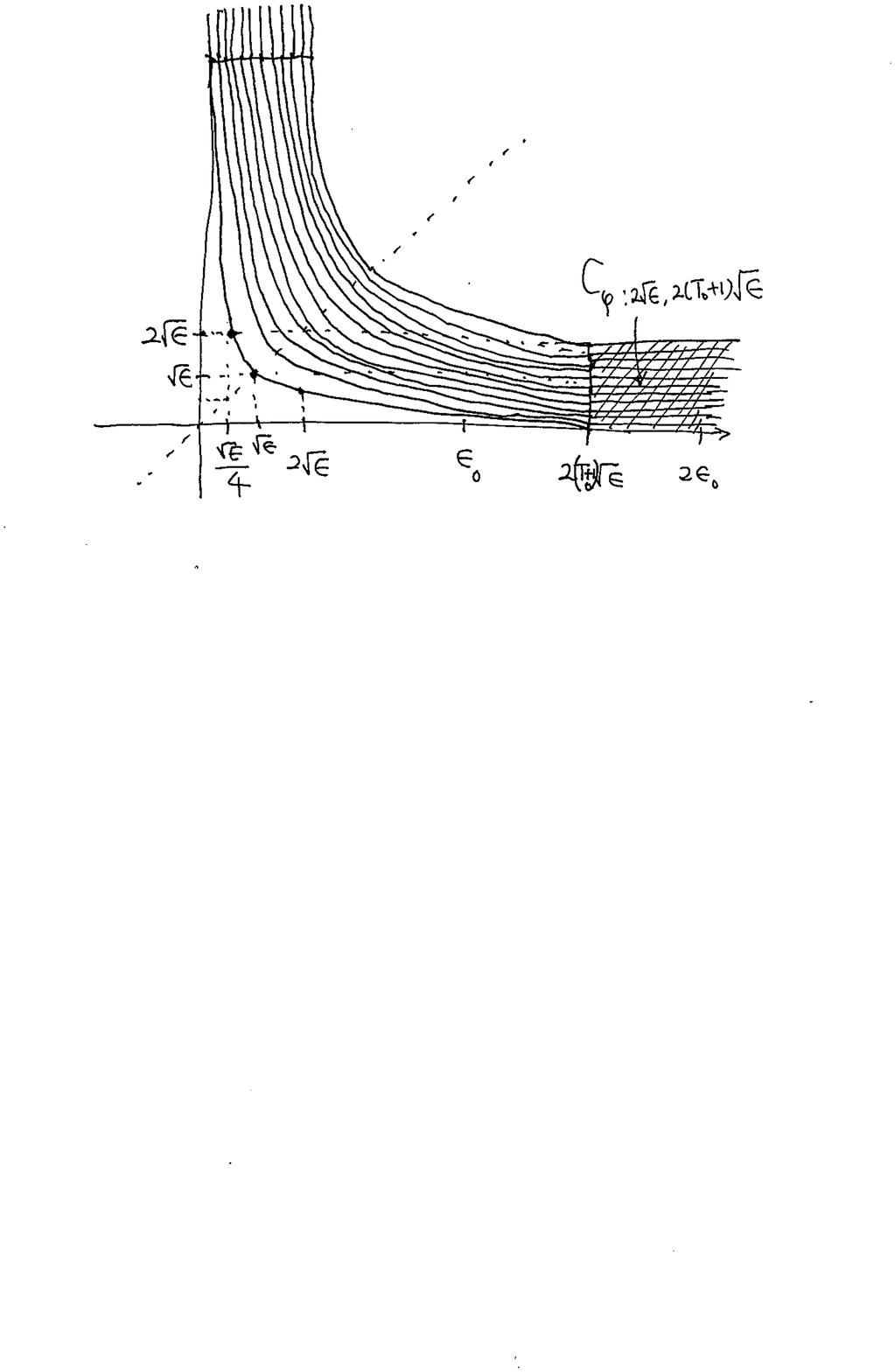}
\caption{Level sets of $\varphi_{2,\epsilon}$}
\label{fig:corner-smoothing}
\end{figure}

By construction, the curve
\eqn\label{eq:tildeS-epsilon}
\widetilde S_{\epsilon}
\eqnd
is a convex curve which is invariant under the reflection $\tau$ along the $t$-axis and satisfies
\eqn\label{eq:interpolation}
\widetilde S_{\epsilon}
=
\begin{cases}
\text{\rm Image } s_{\epsilon} & \quad \text{for } \, |y| \leq \sqrt{2 \epsilon} \}  \\
\text{\rm Image } s_0 & \quad \text{for } |y| \geq T_0 \sqrt{2 \epsilon} \}
\end{cases}
\eqnd
respectively. (See Figure \ref{fig:corner-smoothing}.)
In the $({\bf y};t)$ coordinate, we can express $\widetilde S_\epsilon$
as the graph of a section
\be\label{eq:tildesepsilon}
\widetilde s: \R_{\bf y} \to \R_t.
\ee
Then $\widetilde s_\epsilon$ induces a diffeomorphism between $\R_{\bf y}$
and $\widetilde S_\epsilon$ whose inverse is the restriction of $\pi_{(1,1)}$
to $\widetilde S_\epsilon \subset \R_{\bf y}$.

\subsection{Compatible corner smoothing of $[0,\infty)^n$}
\label{subsec:smoothcornermodel}

In this subsection, we provide a family of local models of
smoothing $[0,\infty)^n$ that is compatible with various $n$ and
with the $\mathsf S_n$-symmetry for the coordinate swapping.
We will also make the smoothing \emph{convex}---this will be
important to our study of pseudoconvexity later. (See Proposition \ref{prop:convex-interpolation}.)

We borrow the construction from \cite[Section 18.5]{fooo:book-virtual}.
The following is a slight variation of \cite[Condition 18.21]{fooo:book-virtual}
arising by stripping away some statements on the smooth structure on $[0,\infty)^k$, which
is not needed here.

\begin{figure}[h]
\centering
\includegraphics[scale=0.3]{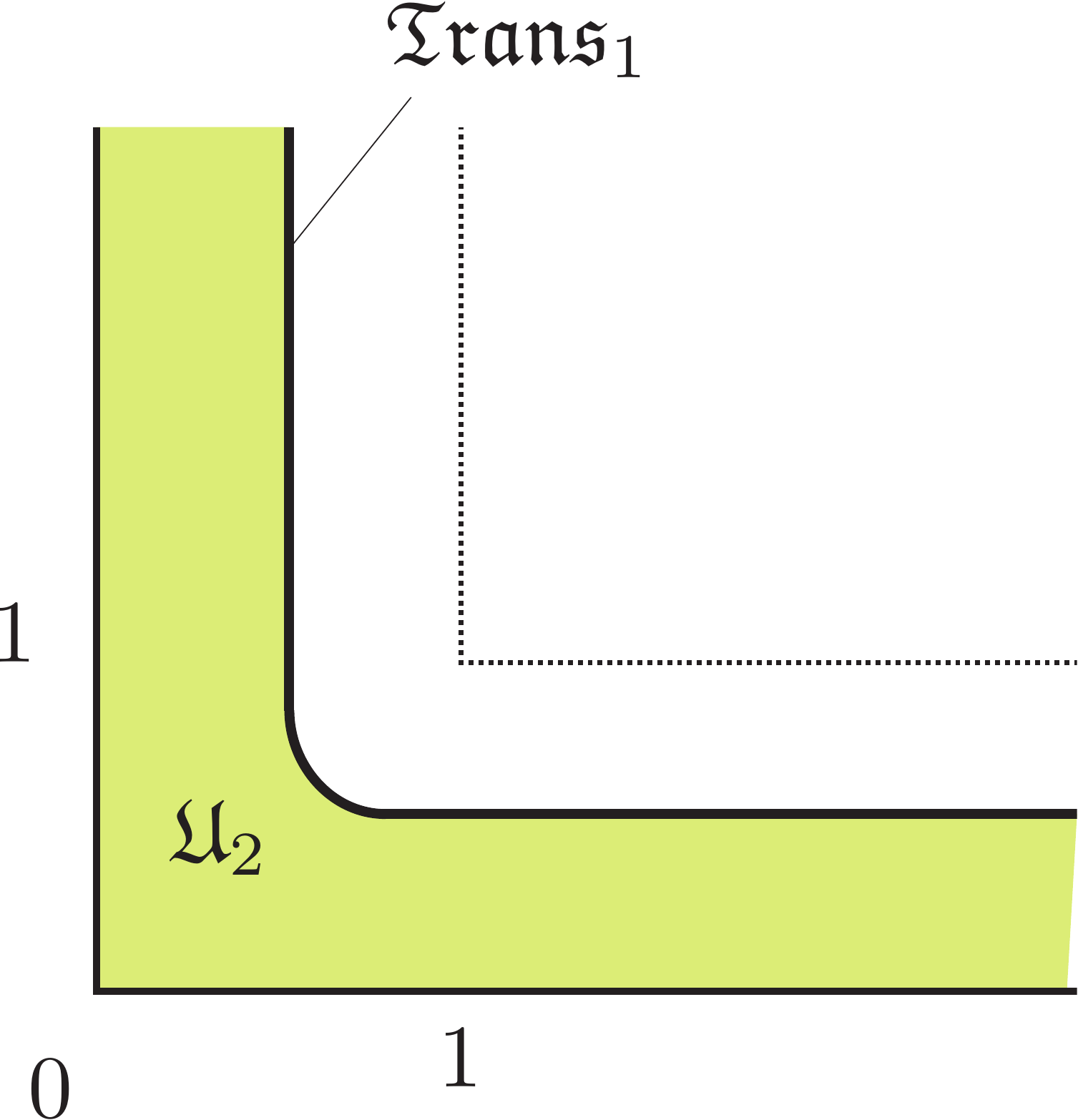}
\caption{$\frak{Trans}_{k-1}$ and $\frak U_k$ \cite[Figure 18.3]{fooo:book-virtual}}
\label{Figure17-3}
\end{figure}

\begin{cond}[Condition 18.21 \cite{fooo:book-virtual}]\label{cond:18.21}
For any $k \in \mathbb Z_{>0}$
we consider $\frak{Trans}_{k-1}$ and $\Psi_k$ with the following properties:
(See \cite[Figure 18.3]{fooo:book-virtual} which is duplicated below in \ref{Figure17-3}.)
\begin{enumerate}
\item
$\frak{Trans}_{k-1}$ is a smooth $(k-1)$-dimensional submanifold of $[0,\infty)^k$
and is contained in $(0,\infty)^k \setminus (1,\infty)^k$.
\item
$\frak{Trans}_{k-1}$ is invariant under the ${\rm Perm}(k)$ action on
$[0,\infty)^k$.
\item
$
\frak{Trans}_{k-1} \cap ([0,\infty)^{k-1} \times [1,\infty))
=
\frak{Trans}_{k-2} \times [1,\infty).
$
This is an equality
as subsets of $[0,\infty)^k
= [0,\infty)^{k-1} \times [0,\infty)$.
\item
$$
\Psi_k : [0,1] \times \frak{Trans}_{k-1}
\to [0,\infty)^k
$$
is a homeomorphism onto its image.
Let $\frak U_k$ be its image.
\item
The subset
$\frak U_k \subset [0,\infty)^k$ is a smooth $k$-dimensional submanifold with boundary
and corners and $\Psi_k$ is a homeomorphism. Moreover
$$
\partial \frak U_k
= \partial([0,\infty)^k) \cup \frak{Trans}_{k-1}
$$
and the restriction of $\Psi_k$ to $\{0\}\times \frak{Trans}_{k-1}$
is a homeomorphism onto $\partial([0,\infty)^k)$.
The restriction of $\Psi_k$ to $\{1\}\times \frak{Trans}_{k-1}$
is the identity map.
\item
$\Psi_k$ is equivariant under the ${\rm Perm}(k)$ action.
(The ${\rm Perm}(k)$
action on $\frak{Trans}_{k-1}$ is defined in Item (2) and
the action
on $[0,\infty)^k$ is by permutation of factors.)
\item
If $s \ge 1$, $t \in [0,1]$ and $(x_1,\dots,x_{k-1}) \in \frak{Trans}_{k-2}$,
then
$$
\Psi_k(t,(x_1,\dots,x_{k-1},s)) = (\Psi_{k-1}(t,(x_1,\dots,x_{k-1})),s).
$$
Here we use the identification in Item (3) to define the left hand side.
\end{enumerate}
\end{cond}

\begin{remark} In fact, the homeomorphism $\Psi_k$ is a diffeomorphism
with respect to some smooth structure on $[0,\infty)^k$, denoted by $\mathfrak{sm}_k$ in
\cite[Condition 18.21]{fooo:book-virtual}. We refer interested readers \cite[Section 18]{fooo:book-virtual} for
a detailed discussion.
\end{remark}

The following is a slight refinement of \cite[Lemma 18.22]{fooo:book-virtual} by the
addition of the convexity statement. We will make our exposition as close to that of
the proof of  \cite[Lemma 18.22]{fooo:book-virtual} as possible.

\begin{lemma}\label{lem:corner-smoothing}
For any $k \in \mathbb Z_{> 0}$ there exists a sequence of pairs
$(\frak{Trans}_{k-1},\Psi_k)$ satisfying Condition \ref{cond:18.21}.
Moreover, for each given $\delta > 0$, we may take them so that
$\frak U_k$ contains $[0,\infty)^k \setminus [1-\delta,\infty)^k$.
\end{lemma}
\begin{proof}
The proof is by induction.
If $k=1$, we define $\frak{Trans}_{0}$ and $\Psi_1: [0,1] \times \{ 1 - \delta/2\} \to [0,\infty)$ by
$$
\frak{Trans}_{0} = \{ 1 - \delta/2\}, \quad \Psi_1(t,1-\delta/2) = (1-\delta/2) t.
$$
in which case there is nothing to prove. When $k = 2$, we define $\frak{Trans}_1$ and
$$
\Psi_2: [0,1] \times \frak{Trans}_1 \to [0,\infty)^2
$$
as follows.

Recall the diffeomorphism $\widetilde s_\epsilon: \R_{\bf y} \to \R_t$
defined in \eqref{eq:tildesepsilon}. For a fixed $\delta > 0$,
using Proposition \ref{prop:tildefepsilon}, we can choose $\epsilon > 0$
and $T_0 > 0$ sufficiently large so that $2T_0 \sqrt{\epsilon} < \delta$.
Then we define
$\mathfrak{Trans}_1$ as the graph of the function
$$
\widetilde s_{\epsilon,\delta}({\bf y}): = \widetilde s_\epsilon({\bf y})
+ \sqrt{2}(1-\delta/2)
$$
in the $({\bf y};t)$ coordinates. Now, for $t \in [0,1]$ and
$z \in \mathfrak{Trans}_1$, we define the function $\Psi_1$ by
$$
\Psi_2(t,z): = tz + (1-t) s_0 \circ \pi_{(1,1)}(z)
$$
where we recall $ \mathfrak{Trans}_1$ depends on the choice of
$\epsilon$ and so does $\Psi_2$. (This being said,
we note that the second summand goes to zero as $\epsilon \to 0$
and hence $\Psi_2^\epsilon$ converges to the identity map
as $\epsilon \to 0$.)
This completes the construction for $k=2$.

Now for higher $k\geq 3$,
suppose we have $\frak{Trans}_{k'-1}$, $\Psi_{k'}$ for $k' < k$ as the induction hypothesis.
Conditions \ref{cond:18.21} (2) (3) determine
$\frak{Trans}_{k-1}$ outside $[0,1]^k$.
Conditions \ref{cond:18.21} (6) (7) determine $\Psi_k$ outside
$[0,1]^k$.
It is easy to see that we can extend them to $[0,1]^k$ and obtain $\frak{Trans}_{k-1}$, $\Psi_k$.
\end{proof}

\begin{remark}\label{rem:epsilon-dependence}
Under the above construction, the map $\Psi_k$ depends on the choice of some parameters such as
$\epsilon, \, T_0$ and $\rho$ which entered in the definition of $\widetilde s_\epsilon$ when $k = 2$.
We will denote $\Psi_k$ by $\Psi_k^\epsilon$ when we need to highlight this parameter dependence of
the definition of $\Psi_k$.
\end{remark}

Then we define a one parameter family of hypersurfaces
\eqn\label{eq:S-epsilon-n}
S_{\epsilon:k}^t: = \frak{Trans}_n - (t + 1/2) (1,1,\ldots,1), \quad t \geq 0.
\eqnd
We note that $S_{\epsilon:k}^0 = \widetilde S_{\epsilon}$ from \eqref{eq:tildeS-epsilon}
by construction, with the dimension $k$ specified.

Finally--recalling the change of coordinates in~\eqref{eqn. Phi}---we consider the change of coordinates
$\Phi^{-1} (= \Pi:) (x_0,\cdots, x_k) \mapsto ({\bf y}; t)$
and define a function $\varphi_k^\epsilon = \widetilde s_\epsilon \circ \Phi^{-1}|_{x_0 = 0}$.
which has the expression
\eqn\label{eq:varphi-epsilon}
\varphi^\epsilon_k(x_1,\cdots, x_k) = t \circ \widetilde \Psi_k^\epsilon\circ \Pi^{-1} (x_1,\ldots x_k)
\eqnd
We consider the function on
$$
{\bf Q}_\epsilon^k: = \{(x_1,\ldots, x_k) \mid \varphi^\epsilon_k(x_1,\cdots, x_k) \geq 0\},
$$
which is the corner-smoothing of ${\bf Q}_+^k$.  (See Figure \ref{Figure17-3}.) By construction,
\eqref{eq:derivative-s-epsilon} implies
\eqn\label{eq:phi-derivatives}
0 \leq \frac{\del \varphi_k^\epsilon}{\del x_i} \leq 1
\eqnd
for all $i = 1,\ldots,k$.

\begin{remark}\label{rem:varphi-epsilon-k} By construction,
the function $\varphi_k^\epsilon$ defined on
 ${\bf Q}_\epsilon^k$ restricts to the lower $\varphi_\ell^\epsilon$'s
with $\ell < k$ on the faces of ${\bf Q}_\epsilon^k$. For example,
\eqn
\varphi_k^\epsilon|_{\del_i {\bf Q}_\epsilon^k}
(x_1,\ldots,0, \ldots,  x_k) =
\varphi_{k-1}^\epsilon(x_1,\ldots, \widehat x_i, \ldots  x_k)
\eqnd
for $x$ with $ \min_i |x_i| \geq 2 \sqrt{\epsilon}, \, x_i = 0$,
and it has the value
$$
\varphi_k^\epsilon(x_1,\cdots, x_k) =  x_1 \cdots x_k - \frac{\epsilon}{2}
$$
for $x$ with $\max_i |x_i| \leq \frac{\sqrt{\epsilon}}{4}$
on ${\bf Q}_\epsilon^k$.
In particular, when $k = 2$, we have the following explicit expression
\eqn\label{eq:varphi-epsilon-2}
\varphi^\epsilon_2(x_1,x_2) = \begin{cases} x_1 \quad & \text{\rm for } \,
\quad x_2 \geq 2T_0 \sqrt{\epsilon}, \, x_1 \leq \frac{\sqrt{\epsilon}}{4} \\
 x_2 \quad & \text{\rm for } \, \quad x_1 \geq  2T_0 \sqrt{\epsilon} , \, x_2 \leq \frac{\sqrt{\epsilon}}{4}\\
x_1 x_2 - \frac{\epsilon}{2} \quad & \text{\rm for } \,
 \quad \frac{\sqrt{\epsilon}}{4} <  x_1, \, x_2 \leq  2\sqrt{\epsilon}
 \end{cases}
\eqnd
and on the level set $(\varphi^\epsilon_2)^{-1}(0)$, we have
$$
x_1 x_2 = \frac{\epsilon}{2}
$$
when $\frac{\sqrt{\epsilon}}{4} <  x_1, \, x_2 \leq  2\sqrt{\epsilon}$.
\end{remark}

The following property of $\varphi_2$ will play an important role later  in our  proof of Proposition \ref{prop:linear-independence}, which is the
key ingredient in our construction of a pseudoconvex pair $(\psi,J)$.
We will also make a precise choice of $\epsilon_1$ later
during the proof thereof.

\begin{lemma}\label{lem:x2dphi>x1dphi} Let $\varepsilon_0 > 0$ be the constant
fixed in Notation \ref{nota:epsilon0}.
Let $\epsilon_1 > 0 $ be sufficiently small and then choose $T_0$ so that
\eqn\label{eq:5/4epsilon}
\frac{5\varepsilon_0}{4} \leq 2T_0 \sqrt{\epsilon_1} \leq \frac{3\varepsilon_0}{2}.
\eqnd
Then we can choose the function $\varphi:= \varphi_2^{\epsilon_1}$ such that it satisfies \eqref{eq:varphi-epsilon-2}
and
\eqn\label{eq:x2dphi>x1dphi}
x_2 \frac{\del \varphi}{\del x_2}(x_1,x_2)\geq x_1 \frac{\del \varphi}{\del x_1}(x_1,x_2)
\eqnd
on the region
\eqn\label{eq:x1-region}
B_{\epsilon,\varepsilon_0}: = \left\{(x_1,x_2) \, \Big| \, 0 \leq x_2 \leq 2 \sqrt{\epsilon_1}, \,\,
0 \leq x_1 \leq \frac{\varepsilon_0}{2} \right\}.
\eqnd
By symmetry, the opposite inequality of \eqref{eq:x2dphi>x1dphi} holds on the
region $\tau(B_{\epsilon,\varepsilon_0})$.
\end{lemma}
\begin{proof} We will check the inequality on each level set
$$
\varphi(x_1,x_2) = c
$$
of the box $B_{\epsilon_1,\varepsilon_0}$ for $0 \leq c \leq 2\sqrt{\epsilon_1}$.
On each level set, the slope of the tangent line of the curve at $(x_1,x_2)$ is given by
\eqn\label{eq:slope}
\frac{d x_2}{d x_1}\Big|_{(x_1,x_2)} = - \frac{\del \varphi}{\del x_1}(x_1,x_2) \Big/ \frac{\del \varphi}{\del x_2}(x_1,x_2)
\eqnd
We start with the level $c = 0$. By the explicit formula
given \eqref{eq:varphi-epsilon-2} away from region
$$
2\sqrt{\epsilon} \leq x_1 \leq 2T_0 \sqrt{\epsilon}
$$
we see that the equality holds in \eqref{eq:x2dphi>x1dphi} for all $0 < \epsilon \leq \epsilon_1$.

Now we let $0 < \epsilon \leq \epsilon_1$.
On the region of our interest, we have the associated curve
$$
\widetilde S_\epsilon \cap \{x_1 \geq \sqrt{\epsilon/2}\}
= \Graph \widetilde f_\epsilon
$$
and hence
\be\label{eq:slope-quotient}
\frac{\del \varphi}{\del x_2}(x_1,x_2)
\Big/ \frac{\del \varphi}{\del x_1}(x_1,x_2)
= - 1\Big/ \widetilde f_\epsilon'(x_1)
\ee
on the region of our interest. In particular we have
$$
x_2 \frac{\del \varphi}{\del x_2}(x_1,x_2)
\Big/ x_1 \frac{\del \varphi}{\del x_1}(x_1,x_2)
= - \frac{\widetilde f_\epsilon(x_1)}{x_1 \widetilde f_\epsilon'(x_1)} \geq 1
$$
where the inequality follows from Proposition \ref{prop:tildefepsilon} (5).

For the positive level $c > 0$, the level set $\varphi = c$ is given by translating it by the direction $c\cdot (1,1)$.  Therefore we have
only to consider the region  with $x_1 \geq \sqrt{c}$. Then the level set $\varphi = c$ is also given by the graph of
the function, denoted by $\widetilde f_\epsilon^c$ whose domain
is shifted by $\sqrt{c}$ from that of $f_\e$ and so recalling the
reflection symmetry along $y = x$, we have only to consider the region
$$
 x_1 \geq \sqrt{c} + \sqrt{\epsilon/2}.
$$
The function $\widetilde f_\epsilon^c$ is given by
$$
\widetilde f_\epsilon^c(x_1) = \widetilde f_\epsilon(x_1 - \sqrt{c}) + \sqrt{c}
$$
and hence we have
$$
\widetilde f_\epsilon(x_1') = \widetilde f_\epsilon^c(x_1'+ \sqrt{c})
- \sqrt{c}
$$
for $x_1' = x_1 - \sqrt{c}$ with $x' \geq \sqrt{\e/2}$. Substituting this into
$$
\widetilde f_\epsilon(x_1')  + x_1' \widetilde f_\epsilon'(x_1') \geq 0
$$
we arrive at
$$
\widetilde f_\epsilon^c(x_1)  - \sqrt{c}
+ (x_1 - \sqrt{c})(\widetilde f_\epsilon^c)'(x_1)  \geq 0.
$$
This proves
$$
\widetilde f_\epsilon^c(x_1)
+ x_1 (\widetilde f_\epsilon^c)'(x_1)  \geq
\sqrt{c} (1 + (\widetilde f_\epsilon^c)'(x_1)) \geq 0
$$
where the last inequality holds since
$(\widetilde f_\epsilon^c)'(x_1) > -1$ by Proposition \ref{prop:tildefepsilon}
(2).
By the same reason as $c = 0$, this implies
the required inequality \eqref{eq:x2dphi>x1dphi} for the level curve $\varphi=c$.
By taking the union of the level curves $\varphi = c$ with $0 < c \leq 2\sqrt{\epsilon_1}$.
This finishes the proof of the lemma.
\end{proof}

Once we have constructed this compatible smoothing, we have the following proposition
whose proof is omitted which is intuitively
obvious because $\widetilde s_\epsilon$ is a smooth approximation of
$s_\epsilon$ as $\epsilon \to 0$. And the functions are obviously $\mathsf{S}_k$-symmetric.

\begin{prop}\label{prop:barrier-construct}
There exists a sufficiently small $\epsilon > 0$ and a sufficiently large
$T_0 > 0$  independent of  the choice of $\sqrt{\epsilon} > 0$
such that  they satisfy the following for all $0 < \sqrt{\epsilon} \leq \varepsilon_0$:
\begin{enumerate}
\item The restriction $\varphi^\epsilon|_{\RR^J}$ are $\mathsf S_{|J|}$-equivariant
for any subsets $J \subset \underline n = \{1,\ldots,n\}$. Here $\RR^J \subset \RR^n$ is the obvious copy of
$\RR^{|J|}$ where we put $\mathsf S_k = \text{\rm Perm}(k)$.
\item The Hessian $\text{\rm Hess}(\varphi^\epsilon)$ of $\varphi^\epsilon$ is positive semi-definite.
\item $\text{\rm Hess}(\varphi^\epsilon) = 0$ for $x_i \geq 2T_0 \sqrt{\epsilon}$ for some $i = 1,\cdots, n$.
\end{enumerate}
\end{prop}
We denote by $\varphi_k$ or $\varphi_k^\epsilon$ the function corresponding to
the corner of codimension $1 \leq k \leq n$. More generally we introduce the set of
such a symmetric convex functions $\varphi$.

\newenvironment{symmetric-convex}{
	  \renewcommand*{\theenumi}{(CV\arabic{enumi})}
	  \renewcommand*{\labelenumi}{(CV\arabic{enumi})}
	  \enumerate
	}{
	  \endenumerate
}

\begin{defn}\label{defn:symmetric-convex} A \emph{symmetric convex smoothing function} on $\RR^k$ is
a function $\varphi: \RR^k \to \RR$ satisfying the following:
\begin{symmetric-convex}
\item\label{item. CV symmetric}
The restriction $\varphi|_{\RR^J}$ is $\mathsf S_{|J|}$-invariant
for all subsets $J \subset \underline n = \{1,\ldots,n\}$. Here $\RR^J \subset \RR^n$ is the obvious copy of
$\RR^{|J|}$.
\item\label{item. positive definite}
 $\text{\rm Hess}(\varphi)$ is positive semi-definite everywhere.
\item\label{item. hessian}
$\text{\rm Hess}(\varphi|_{\RR^J})$ is compactly supported on $\RR^J$ for all subsets $J \subset \underline n$ for $|J| \geq 1$.
\end{symmetric-convex}
We denote the set thereof by $\mathfrak{Conv}^{\mathsf S_k}_{\mathfrak{sm}}(\RR^k)$. We also consider the
set of \emph{nonsymmetric} convex smoothing function consisting of those not necessarily satisfying \ref{item. CV symmetric},
which we denote by $\mathfrak{Conv}_{\mathfrak{sm}}(\RR^k)$
\end{defn}
Proposition \ref{prop:barrier-construct} shows that both sets are nonempty, which are also convex and so contractible.

For given $\varphi$, we define
\eqn\label{eq:linear}
U_{\varphi;i} = \{(x_1,\ldots, x_k) \mid \varphi(x_1, \ldots, x_k) = x_i\}
\eqnd
for each $i=1, \ldots, k$. (We alert readers that $U_{\varphi;i}$ is a codimension zero subset of
$\RR^k$ by the property of the function $\varphi$ that satisfies $\varphi = x_i$
near the coordinate planes of $x_j$ with $j \neq i$ and $x_i$ is away from $0$.)
For each $i$, we consider the rectangular half cylinder defined by
\eqn\label{eq:Cied}
C_{\varphi;\epsilon,\delta}: = \{(x_1,\ldots, x_k) \in \RR^k_+ \mid x_i \geq \delta, \,  0 \leq x_j \leq \epsilon \, \forall j \neq i \}
\eqnd
for each pair  $\epsilon, \, \delta > 0$ of positive constants.
It follows from Definition \ref{defn:symmetric-convex} \ref{item. hessian} that we can choose
$0 < \epsilon < \delta$ so that
$
C_{\varphi,i;\epsilon,\delta}
$
is contained in $ U_{\varphi;i}$. (See Figure \ref{fig:corner-smoothing}.)

The collection thereof, which we denote by $\cC_i(\varphi)$, is a partially ordered set by inclusion.
The following is obvious by the defining conditions on $\varphi$.

\begin{lemma}\label{lem:maximal-element}
The maximal element of any nonempty chain $\cD$ of \text{\rm POSET} $\cC_i(\varphi)$ is also
a rectangular half cylinder $C_{\varphi;\epsilon',\delta'}$ for some $0 < \epsilon' \leq \delta' < \infty$.
\end{lemma}

\begin{defn}[Width and Height of $\varphi$]\label{defn:height-width}
Let $\varphi \in \mathfrak{Conv}^{\mathsf S_k}_{\mathfrak{sm}}(\RR^k)$.
\begin{enumerate}
\item
For given $i$ and a chain $\cD_i$ of $\cC_i(\varphi)$,
we define the \emph{$i$-th height}, denoted by $\mathsf{ht}(\varphi,i;\cD_i)$, to be the $\delta'$, and
the \emph{$i$-th width}, denoted by $\mathsf{wd}(\varphi,i;\cD_i)$, to be the $\epsilon'$ appearing
in the above lemma. Denote by $\cD$ the symmetric collection consisting of $\cD_i$ for $i=1, \ldots, k$.
\item Choose the symmetric collection $\cD$ above so that $\cD_i$'s are  pairwise disjoint for $i = 1, \ldots, k$.
We define a \emph{height} of $\varphi$, denoted by $\mathsf{ht}(\varphi;\cD)$,
to be
$$
\mathsf{ht}(\varphi;\cD): = \mathsf{ht}(\varphi,i;\cD_i)
$$
for a (and so any) choice $i=1, \ldots, k$. We define the width $\mathsf{wd}(\varphi;\cD)$ by the same way.
\end{enumerate}
\end{defn}
Thanks to the symmetry of $\varphi$, $\mathsf{ht}(\varphi;\cD)$ is well-defined independent of $i$'s.
Since we will just need one choice of such a collection $\cD$ for our purpose, we will
omit the dependence on $\cC$ from the notations of
$$
\mathsf{ht}(\varphi;\cD) =:\mathsf{ht}(\varphi),
\quad \mathsf{wd}(\varphi;\cD) = : \mathsf{wd}(\varphi).
$$
By definition, we have
\eqn\label{eq:linearity-region}
\varphi(x_1, \ldots, x_k) = x_i \quad \text{when $x_i \geq \mathsf{ht}(\varphi)$,
 and $0 \leq x_j \leq \mathsf{wd}(\varphi)$ for $j \neq i$}
\eqnd
for any given $\varphi$.
The following is also obvious from definition.

\begin{cor}\label{cor:uniform-width}
\begin{enumerate}
\item For any compact collection of $\cC \subset \mathfrak{Conv}_{\mathfrak{sm}}^{\mathsf{S}_k}(\RR^k_+)$,
$$
0 < \inf_{\varphi \in \cC} \mathsf{wd}(\varphi) <
\sup_{\varphi \in \cC} \mathsf{ht}(\varphi) < \infty.
$$
\item For every $\eta > 0$ and $A > 0$, there exists $\varphi \in \mathfrak{Conv}_{\mathfrak{sm}}^{\mathsf{S}_k}(\RR^k_+)$ such that
 $\mathsf{wd}(\varphi) < \eta$ and
$\mathsf{ht}(\varphi)/\mathsf{wd}(\varphi) > A$.
\item $d\varphi(x_1, \ldots, x_k) = 0$ if and only if $(x_1, \ldots, x_k) = 0$.
\end{enumerate}
\end{cor}

\begin{example}[Width and height of the model]
By definition, we have
$$
\mathsf{wd}(\varphi_{k+1}^\epsilon) = \frac{\sqrt{\epsilon}}{4}, \quad
\mathsf{ht}(\varphi_{k+1}^\epsilon) = 2T_0 \sqrt{\epsilon}, \quad
\mathsf{ht}(\varphi)/\mathsf{wd}(\varphi) = 8 T_0
$$
for those $\varphi = \varphi_{k+1}^\epsilon$ explicitly constructed in Proposition \ref{prop:barrier-construct}.
\end{example}

\begin{remark}\label{rem:nonsymmetric-case} Even when $\varphi$ is not symmetric, we can define
the width and the height $\mathsf{wd}(\varphi)$ and $\mathsf{ht}(\varphi)$ as follows.
Lemma \ref{lem:maximal-element} still enables us to define
$\mathsf{wd}(\varphi;i,\cC_i)$ (resp. $\mathsf{ht}(\varphi;i,\cC_i)$). Then consider
a (nonsymmetric) collection $\cC = \{\cC_i\}$ which are
still mutually disjoint. Then we define
\beastar
\mathsf{wd}(\varphi;\cC) &: =&  \min_{i=1, \ldots, k} \mathsf{wd}(\varphi;i,\cC_i)\\
\mathsf{ht}(\varphi;\cC) &: =&  \max_{i=1, \ldots, k} \mathsf{ht}(\varphi;i,\cC_i).
\eeastar
Again we only need one such choice $\cC$, we will omit $\cC$ from their notations.
\end{remark}

\section{End-profile functions ${\frak s}_\varphi$ of Liouville sectors}
\label{sec:barrier-functions}

Let $\Mliou$ be a Liouville sector without corners. We will fix a smoothing profile (Condition \ref{cond:smoothing-profile}), so that we in particular have:
\begin{itemize}
\item a radial coordinate $s$ as in~\eqref{eqn. s} (giving rise to a contact-type hypersurface $S_0 \subset M$),
\item a convex smoothing function $\varphi$ for $\del \Mliou$, and
\item Splitting data $(F, \{(R_i, I_i)\}_{i=1}^k)$ for every sectorial corner $C_\delta$ of codimension $k$.

\end{itemize}
%
%

\subsection{Smoothing profiles}
\label{subsec:smoothing-profile}

Now we take $\varphi \in \mathfrak{Conv}_{\mathfrak sm}^{\mathsf{S}_2}(\RR_+^2)$ and consider
the function
\eqn\label{eq:s2-varphi}
s_{1+1,\varphi} : = - \log \varphi\left(R, e^{-s}\right)
\eqnd
on the union
$$
(\nbhd_{2\varepsilon_0}(\del M) \setminus \del M )\cup \nbhd_{s \geq N}(\del_\infty M).
$$
By the properties of $\varphi_{2,\epsilon}$, we have
$$
s_{1+1,\varphi} = \begin{cases}
-\log R (= \mu) \quad & \text{for }\, e^{-s} \geq \mathsf{ht}(\varphi), \, R \leq \mathsf{wd}(\varphi) \\
s \quad & \text{for }\,  R \geq \mathsf{ht}(\varphi), \, e^{-s}\leq \mathsf{wd}(\varphi)
\end{cases}
$$
We note that both $R$ and $s$ are defined on the interpolating region
and hence the function $s_{1+1,\varphi}$ is smooth.

Note that the level sets of this function
provide the asymptotic profile of the horizon $\del_\infty M \cup \del M$ under the coordinate
change $R \mapsto -\log R: = \mu$, which will play a fundamental role in our definition of
the sectorial almost complex structures.

\begin{defn}[End-profile function $s_{1+1,\varphi}$]
\label{defn. s2phi}
Let a contact-type hypersurface $S_0$ be given as above. For each sectorial corner $C_\delta$,
We fix a splitting data $(F,(R,I))$ and
a convex smoothing function $\varphi \in \mathfrak{Conv}_{\mathfrak sm}^{\mathsf{S}_2}(\RR^2_+)$
satisfying
$$
\varepsilon_0 \leq \mathsf{ht}(\varphi) < \frac{3}{2}\varepsilon_0
$$
so that their union defines a smoothing function by a partition of unity. We denote this glued function by
$$
{\mathfrak s}_\varphi
$$
and call an \emph{end profile function} of $M$.
\end{defn}
We consider the function ${\mathfrak s}_\varphi$ on a neighborhood $\del_\infty M \cup \del M$,
\eqn\label{eq:nbhd-s2epsilon-N}
\nbhd_{s \geq N}(\del_\infty M) \cup \nbhd_{2\varepsilon_0}(\del M)
\eqnd
for varying $N$ over $N \geq 0$.

\begin{lemma}\label{lem:exhaustion} Let $\varphi = \varphi_{2,\epsilon}
 \in \mathfrak{Conv}^{\mathsf{S}_2}_{\mathfrak{sm}}(\RR^2_+)$
with $\mathsf{wd}(\phi) \ll \varepsilon_0$ be as given before,
and consider the end-profile function ${\mathfrak s}_\varphi$ associated to $\varphi$.
Then each level set $({\mathfrak s}_\varphi)^{-1}(r) \cap \mu^{-1}((-\infty,K])$ is a
smooth compact hypersurface for each $r > 0$ and $K \geq -\log(2\varepsilon_0)$.
\end{lemma}
\begin{proof} At each sectorial corner $C_\delta$, we compute the differential of its associated
end-profile function $s_{1+1,\varphi}$.
We compute
and
\be\label{eq:LHS}
ds_{1+1,\varphi} = - \frac{d\varphi}{\varphi}
= \frac{1}{\varphi}\left(-\frac{\del \varphi}{\del x_1} dR
+ e^{-s}\frac{\del \varphi}{\del x_2} ds \right)
\ee
with $\frac{\del \varphi}{\del x_1}$ and $\frac{\del \varphi}{\del x_2}$ evaluated at
$(R,e^{-s})$. Since $\del M \pitchfork \del_\infty M$, the two one-forms
$dR, \, ds$ are linearly independent.
Since $\frac{\del \varphi}{\del x_2}(R,e^{-s})\neq 0$ in the region of
interest, this vector is also nowhere vanishing. In particular $ds_{1+1,\varphi}$ is
nowhere vanishing which proves smoothness of $s_{1+1,\varphi}^{-1}(r)$ by the regular value theorem.
The other compactness statement is obvious by definition.
\end{proof}

\begin{figure}[h]
\centering
\includegraphics[scale=0.5]{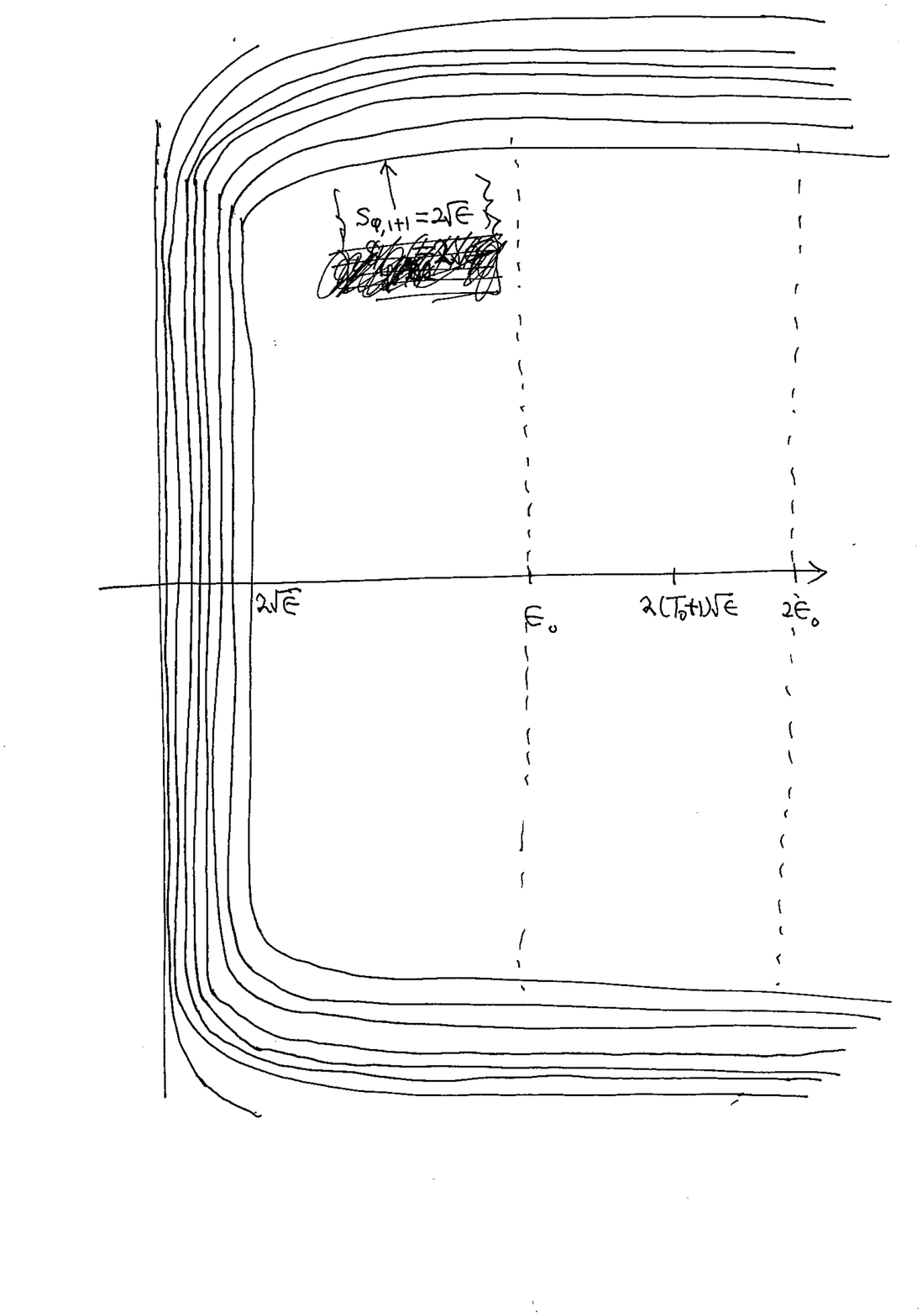}
\caption{Foliation by an exhaustion function}
\label{fig:foliation}
\end{figure}

\begin{remark}\label{rem:not-gradient-like}
\begin{enumerate}
\item As mentioned in the introduction,
 the hypersurfaces ${\mathfrak s}_\varphi^{-1}(r)$
may not be of contact-type but a neighborhood
$$
\nbhd (\del_\infty M \cup \del M)
$$
is exhausted by them by construction for sufficiently large $r$'s, which will be sufficient for our purpose.
We suspect that using the non-contact type hypersurfaces may be inevitable
to make the almost complex structure amenable to the application of maximum principle
in both directions of $\del_\infty M$ and $\del M$ of the Liouville sector $M$.
(See Remark \ref{rem:difficulty} (3).)
\item It is more natural to regard the above exhausting foliation as the
smoothing of the ceiling corner of the compact sectorial domain
$$
W = M \setminus (S_0 \times [0,\infty))
$$
or the ideal completion $\overline M$. (See Appendix \ref{sec:ACI-framework}.)
The coordinate change $s \mapsto e^{-s}$ or $s \mapsto e^{N-s}$ by translating the radial
coordinates reflects this point of view, which is also related to the coordinate change
$$
(R,I) \mapsto ( - \log R, I) = : (\mu,\nu)
$$
we will use in the next section. (See Figure \ref{fig:foliation}.)
\end{enumerate}
\end{remark}

\subsection{The case with corners}

We can also apply the above definition of $s_{1+1,\varphi}$ to $s_{k+1,\varphi}$ for the
case of Liouville sectors with corners which is now in order.

Recall from  Definition \ref{defn:sectorial-corner} that $\del M$ is a union of
a collection of cleanly intersecting hypersurfaces $H_1, \ldots, H_k \subset M$ (cylindrical at infinity)
near each sectorial corner of $\del M$ of codimension $k$ for some $k$.
We now focus on a corner with the codimension $1 < k < n$ with $n = \frac12 \dim M$
since remaining cases are easier.
Recall from Lemma \ref{lem:splitting-corner} that $\{R_i,I_i\}_{i=1}^k$ form a
canonical coordinates near the corner over the subset
$$
\bigcap_{i=1}^k \{I_i = 0\} \cong F \times [0,\infty)^k
$$
and $H_i = \{R_i = 0\}_{i=1}^k$. As before, we make the coordinate changes
$$
\mu_i = -\log R_i, \quad \nu_i = I_i
$$
as in \eqref{eq:coordinate change} and write $\mu = (\mu_1, \ldots, \mu_k)$ and $\nu = (\nu_1, \ldots, \nu_k)$.
Then we consider the
version of barrier functions that incorporates the mixture of the ceiling corner and sectorial
corners. By choosing a splitting data $(F, \{R_i,I_i\}_{i=1}^k\}$ with eccentricity $0 < \alpha < 1$, we have a splitting
$$
\nbhd^Z\left(\bigcap_{j =1}^k \nbhd_{\epsilon}(H_{i_k}))\right) \cong F \times T^*\RR_+^k
$$
of a neighborhood of the intersection $\bigcap_{j =1}^k \nbhd_{\epsilon}(H_{i_k})$.

As before we fix a contact-type hypersurface $S \subset M$ so that
and a Liouville embedding $S \times [0,\infty) \hookrightarrow M$ and let
 $s= s_{S_0}$ be the associated radial function.
Next, for each given $\varphi_k \in \mathfrak{Conv}_{\mathfrak{sm}}(\RR^k)$
we consider the function
\eqn\label{eq:s2-varphik+1}
s_{k+1,\varphi} =  - \log \varphi_{k+1}\left(R_1, \cdots, R_k,e^{-s}\right)
\eqnd
in the same way as in Subsection \ref{subsec:smoothing-profile}.

\emph{From now on, we will not further mention the possibility that there may be more
than one component of the sectorial corner, since there will be no essential differences in its exposition
from the case of one sectorial corner except notational complexity added.}

The following premises will be put for the definitions of all the sectorial things we
introduce in the rest of the paper.
\begin{cond}[Smoothing profile] \label{cond:smoothing-profile}
\begin{itemize}
\item We fix a contact-type hypersurface $S_0 \subset M$ and
the associated decomposition
$$
M = W \cup_{\del W} (S_0 \times [0,\infty)),
$$
and the associated radial function $s= s_{S_0}$ on $S_0 \times [0,\infty)$,
where $W$ is a compact Liouville domain with boundary  $\del W = S_0$.

At each sectorial corner $C_\delta$ of $\nbhd(\del_\infty M \cup \del M)$,
we fix the following data:
\item  a splitting data $(F_\delta, \{(R_{\delta,i},I_{\delta,i})\})$ with
$F_\delta= (R_\delta,I_\delta)^{-1}(0,0)$ with $R_\delta = (R_{\delta,1}, \ldots, R_{\delta,k})$ and
similarly for $I_\delta$.
\item convex smoothing function
$$
\varphi_\delta = \varphi_{k_\delta +1}  \in \mathfrak{Conv}_{\mathfrak{sm}}(\RR^{k_\delta+1}).
$$
and its associated end-profile function
$$
s_{k_\delta+1,\varphi}: \nbhd(C_\delta) \to \RR.
$$
\end{itemize}
We glue the above end-profile functions of the corners by a partition of unity and denote by
$$
{\mathfrak s}_\varphi, \quad \varphi := \{\varphi_\delta\}
$$
the resulting function, which we call an \emph{end-profile function} of the full boundary of $\overline M$,
$DM=\del_\infty M \cup \del M$.
\end{cond}
We have a natural diffeomorphism between $\del_\infty M$ and the level set
$$
s^{-1}(N) = : \del_{N}M \cong S_0
$$
for any $N \geq 0$, which is induced by the Liouville flow.

\begin{remark}
It is shown in \cite{oh:intrinsic} that a choice of section for
the leaf space fibration $\pi: \del M \to \cN_{\del M}$ canonically induces a splitting
$$
\nbhd(\del M) = F \times \C_{\text{\rm Re}^k \geq 0}, \quad \omega = \pi_F^*\omega_F + \sum_{i=1}^k dR_i \wedge dI_i
$$
via a natural $\RR^k$-equivariant symplectic diffeomorphism.
\end{remark}

\section{Intrinsic geometry of end-profile functions and splitting data}
\label{sec:intrinsic-geometry}

In this section, we explain intrinsic geometry of end-profile functions and
the splitting data laid out in \cite{oh:intrinsic} which
provide some intrinsic meaning and perspective with
the constructions of sectorial Floer package we perform in Part II
which much uses coordinate calculations.

We start with the splitting data for the Liouville sectors with corners.
In \cite{oh:intrinsic}, in the course of providing an intrinsic characterization of
Liouville sectors with corners, the present author considers \emph{clean coisotropic collection}
$$
\{H_1, \cdots, H_k\}
$$
of hypersurfaces, and the natural projections
$$
\pi_{H_i}: H_i \to \cN_{H_i}
$$
to the space of leaves of the characteristic foliation of $H_i$.
It is shown that at each sectorial corner
$$
C_\delta = H_1 \cap H_2 \cap \cdots \cap H_k
$$
associated to the collection $\{H_1, \cdots, H_k\}$, a choice of smooth sections
$$
\{\sigma_1, \cdots, \sigma_k\}
$$
with $\sigma_i: \cN_{H_i} \to H_i$ of the fibration $\pi_{H_i}: H_i \to \cN_{H_i}$,
induces a free $\RR^k$-action on $C_\delta$ which
simultaneously linearizes the characteristic flows of $H_i$ which leads to
the collections of functions $\{R_1, \cdots, R_k\}$ as the \emph{action variables}
in the integrable systems.

Then the $\RR^k$-action is canonically extended to a neighborhood $\nbhd(C_\delta) \subset M$
using Gotay's coisotropic embedding theorem \cite{gotay} of presymplectic manifolds $(H_i, \omega_{H_i})$
where $\omega_{H_i} = (d\lambda)|_{H_i}$.  Its moment map
$$
\phi_G:\nbhd(C_\delta) \to {\frak g}^* \cong \RR^k
$$
for the action of the group $G = \RR^k$ is precisely the map
$$
\phi_G = (R_1, \cdots, R_k): \nbhd(C_\delta) \to \RR^k
$$
where the functions $R_i$ are the Hamiltonian of the vector field $\widetilde {\bf e}_i$
associated to the basis element ${\bf e}_i \in {\frak g}$, i.e.,
$$
dR_i = \widetilde{\bf e_i} \rfloor \omega
$$
for $i=1, \ldots, k$, with the normalization
$
R_i|_{H_i} = 0
$
for each $i$. Then the condition that the Liouville vector field of M is outward
pointing and the equation $\{R_i, I_j\} = \delta_{ij}$ implies $R_i \geq 0$
for all $i$ on a neighborhood $\nbhd(C_\delta)$ and hence the codomain of $\phi_G$ is
indeed $\RR_+^k$. This not only provides a canonical splitting data $\{(R_i,I_i)\}$
in the description of neighborhood structure of $\del M$ given in \cite{gps-2}, but also
provides a collaring, or a coordinate chart, of the corner.

Then the functions
$$
s_{k+1, \varphi} = \varphi_{k+1}(R_1, \cdots, R_k, e^{-s})
$$
comprising the end profile function ${\frak s}_\varphi$, can be replaced by
$$
s_{k+1, \varphi} = \varphi_2(s_{k,\varphi}, e^{-s})
$$
instead  by breaking the symmetry from the direction of the ceiling $\del_\infty M$ but
keeping those of sectorial corners. Here the function $s_{k,\varphi}$ can be written as
$$
s_{k,\varphi} = -\log \varphi_k \circ \phi_G.
$$
Here the function $\varphi_k \circ \phi_G$ smoothly approximates the sectorial corner
of
$$
\nbhd(\del M) \cong F \times \C_{\text{\rm Re} \geq 0}^k
$$
in coordinates $(R_1, \cdots, R_k)$.
Note that both of two functions in the composition $\varphi_k \circ \phi_G$
carry natural geometric meaning: $\varphi_k$ is a canonically
constructed universal convex function which smoothly approximates the corner structure of $\RR^k_+$,
and $\phi_G$ is the moment map of the $\RR^m$ symplectic action of $\nbhd(\del M)$ which is
intrinsically defined in the coordinate free way once the collection of sections
$$
\{\sigma_1, \ldots, \sigma_k\}
$$
which substitutes the coordinate expression
$$
(R_1, \cdots, R_n).
$$
In this regard, the end-profile function $s_{k+1, \varphi}$ is intrinsic in that it depends only on the sectorial collection
$\{H_1, \cdots, H_k\}$ defining the sectorial corner $C_\delta$, and
\emph{a choice of these sections $\sigma_1, \cdots, \sigma_k$}.
Likewise, the (global) end-profile function ${\mathfrak s}_\varphi$ is intrinsic in this sense.

We would like to compare this with the \emph{$\lambda_\kappa$-wiggled end-profile function} ${\mathfrak s}_{\varphi;\kappa}$
which will be introduced in the next section. The latter depends on an additional datum of
\emph{infinity-moving Liouville deformation} of Liouville form $\lambda$ associated to the cut-off function
$\kappa$ to be defined in the next section.
This is the reason why we differentiate two types of end-profile functions ${\mathfrak s}_\varphi$ and
${\mathfrak s}_{\varphi,\kappa}$ for which we maintain different notations and name differently.

\section{Cutting-off tails of and deformation of Liouville forms}
\label{sec:cutting-off-form}

Recall that we would like $J$ to satisfy
$$
-ds\circ J = \lambda \quad \text{\rm on $\nbhd(\del_\infty M) \setminus \nbhd(\del_\infty M \cap \del M)$},
$$
and
$$
-dR \circ J = \pi_F^*\lambda_F + \pi_\C^*\lambda_\C \quad \text{\rm on $\nbhd(\del M) \setminus \nbhd(\del_\infty M \cap \del M)$}
$$
and to interpolate the two by considering a deformation of Liouville one-forms of the type
$$
\lambda_\kappa = \lambda - d((1-\kappa) f)
$$
where $\kappa$ is a function
as mentioned in the introduction.
\begin{remark} The reason why we write the deformation in this somewhat odd way is motivated by
another expression of $\lambda_\kappa$, which is
$$
\lambda_\kappa = \pi_F^*\lambda_F + \pi_\C^*\lambda_\C + d(\kappa \, f).
$$
\end{remark}

Because $\kappa$ may not vanish on the region
$$
F \times \{|I| \leq N_1\}
$$
it is not Liouville equivalent to $\lambda$ in the standard
sense by compactly supported Liouville diffeomorphisms.
(Recall that the subset $\{|I| \leq N_1\}$ is not compact.)

In Section \ref{sec:deformation-lambda}, we will deform $\lambda_\kappa$
back to $\lambda$ by applying a Moser-type deformation to an isotopy of one-forms
\eqn\label{eq:lambda-t}
t \mapsto \lambda_t:= \lambda - t d((1-\kappa) f).
\eqnd

 We start with providing the precise requirement for $\kappa$ we want.
 For this purpose, we summarize a few basic geometric properties
carried by the neighborhood of the ceiling corner
 $$
\nbhd( \del_\infty M \cap \del M)
$$
derived from Proposition \ref{prop:gps} in the following lemma.
Denote by $\pi_\C$ he projection to $T^*\C$ with respect to the decomposition
$$
T^*(F \times \C) = (T\C)^\perp \oplus (TF)^\perp \cong T^*F \oplus T^*\C
$$
on $\nbhd(\del M)$.

\begin{lemma}[On the ceiling corner]\label{lem:ceiling-corner}
Let $f$  and $F_0 \subset F$ be as in Proposition \ref{prop:gps}.\begin{enumerate}
\item The Liouville vector field
$Z$ is tangent to the level set $I^{-1}(0)$ in
$$
\{s \geq N_2\} \supset \nbhd( \del_\infty M \cap \del M).
$$
\item On $\{|I| \geq N_0\} \cap F \times \C_{\text{\rm Re}\geq 0}$,
$d_\C f = 0$ for a sufficiently large $N_0> 0$. In particular, we have
\be\label{eq:dCf}
\pi_\C(\lambda) = \pi_\C^*\lambda_\C = (1-\alpha)R\, dI - \alpha I\, dR
\ee
on
\be\label{eq:dCf-supp}
F  \times \{|I| \geq N_0\} \bigcup (F \setminus F_0) \times \C_{\text{\rm Re}\geq 0}.
\ee
\end{enumerate}
\end{lemma}
\begin{proof} Statement (1) follows from $Z[I] = \alpha I$ near infinity i.e.,
on $\nbhd(\del_\infty M)$ and hence on $\{s \geq N_2\}$ for
a sufficiently large $N_2 > 0$. Statement (2) is a restatement of combination of
Proposition \ref{prop:gps} (f1), (f2).
\end{proof}

It follows from (2) that  we do not need to make nontrivial
deformation on \eqref{eq:dCf-supp} even though $\pi_F(\lambda_\kappa)$ may be different from $\pi_F(\lambda)$ because
the proof of Proposition \ref{prop:linear-independence} does not
need it. We would like to mention that the lemma is the main
driving force for  our choice of the support of $1-\kappa$ that way.

This motivates us to introduce the following type of a cut-off function.

\begin{defn}[Cut-off function $\kappa$]\label{defn:kappa} Let $\epsilon > 0$ be given
and ${\mathfrak s}_\varphi$ be the
end-profile function of a smoothing profile of $\nbhd_{2\varepsilon_0}^Z(\del M) \cup \nbhd(\del_\infty M)$
that satisfies
$$
\frac{\sqrt{\epsilon}}{4} < \mathsf{wd}(\varphi) < \frac{\varepsilon_0}{8},
\quad \varepsilon_0 < \mathsf{ht}(\varphi) \leq \frac{3\varepsilon_0}{2}.
$$
Then we consider cut-off functions $\kappa$ of the type
\be\label{eq:kappa}
\kappa(y,R,I) = \kappa_0(y) \kappa_1(R) \kappa_2(I)
\ee
that satisfies the following:
\begin{enumerate}
\item  We choose another such Liouville domain $F_0'$
such that $F_0 \subset \overline F_0 \subset \text{\rm Int} F_0'$ and let
$\kappa_0$ satisfy
$\kappa_0(y) \equiv 0$ on $F_0$ and $\kappa_0 \equiv 0$ on $F \setminus F_0'$.
\item $\kappa_2$ satisfies
$$
\kappa_1(R) = \begin{cases} 0 \quad & \text{\rm for } 0 \leq R \leq \sqrt{2} \mathsf{wd}(\varphi) \\
1\quad & \text{\rm  for } \mathsf{ht}(\varphi) \leq R \leq \frac{3\varepsilon_0}{2}
\end{cases}
$$
for some $\epsilon_1$ with $ \sqrt{2} \mathsf{wd}(\varphi) < \epsilon_1$.
(Recall $\mathsf{ht}(\varphi) = 2T_0 \sqrt{\epsilon_1}$.)
\item $\kappa_3$ satisfies
$$
\kappa_2(I) =
\begin{cases} 0 \quad & \text{\rm for } |I| \leq N_0,\\
1 \quad & \text{\rm for } |I| \geq N_1
\end{cases}
$$
for $N_1$ so large that $e^{-N_1} < \sqrt{2}\mathsf{wd}(\varphi)$.
\end{enumerate}
We call $\kappa$  a \emph{deformation function} adapted to the smoothing profile.
\end{defn}

Using this family of cut-off functions \eqref{eq:kappa},
we introduce the associated deformations of
the given Liouville one form $\lambda$.

\begin{defn}[Deformed Liouville form $\lambda_\kappa$]
Let $\kappa$ be any \emph{deformation function} adapted to the smoothing profile.
We consider a deformation associated to $\kappa$ of the given Liouville one form $\lambda$
by
\bea \label{eq:lambda-kappa}
\lambda_\kappa & : = &  \lambda -d((1-\kappa) f) \\
& = & \pi_\C^*\lambda_\C + \pi_F^*\lambda_F + d(\kappa f)
\eea
and call it a \emph{$\kappa$-deformation} of $\lambda$.
\end{defn}

An upshot of this choice made in Condition \ref{cond:s}
for the cut-off function $\kappa$ is the following list of the properties
thereof:

\begin{remark}[Properties of $\kappa$] \label{rem:kappa}
\begin{enumerate}[{(a)}]
\item The choice of the supports of the functions
$1- \kappa$ and $\kappa$ are motivated by the
proof of Proposition \ref{prop:linear-independence} and Theorem \ref{thm:gray}.
\item For all choice of $\kappa$, we have
$d\lambda_\kappa = d\lambda = \omega$.\item
The support condition of $f$ implies $\lambda_\kappa = \pi_F^*\lambda_F + \pi_\C^*\lambda_\C$
on the region $F \times \C_{\text{\rm Re}\geq 0}$.
\item Condition (3) of Definition \ref{defn:kappa} implies
 $$
 \lambda_\kappa = \lambda
 $$
on the region $(F \setminus F_0) \times \C_{\text{\rm Re}\geq 0}$
\emph{near $\del_\infty M$ and away from $\del M$}.
\end{enumerate}
\end{remark}
%
%
\clearpage

\part{Sectorial Floer data and the maximum principle}

Let us recall that it is common to endow a Liouville manifold with an almost-complex structure $J$ \emph{of contact type}
for which we have the equation
\eqn\label{eq:bdy-condition}
- d(d(s\circ v) \circ j) = v^*d\lambda.
\eqnd
This enables one to prove $C^0$ estimate via the maximum principle as follows:
\begin{itemize}
\item First of all, the equation
implies subharmonicity $\Delta(s \circ v)\geq 0$ of the function $s\circ v$ so that the maximum principle applies.
\item  If a Lagrangian $L$ satisfies the property
that \emph{$L \cap s^{-1}(s_0)$ is Legendrian} for all $s_0 \in \RR$ large enough,
which always holds for cylindrical (i.e., $Z$-invariant) Lagrangians,
then the strong maximum principle also applies to $s\circ v$.
\end{itemize}
In general such an identity fails to hold for almost complex structures that are not of contact type---for example, the cylindrical-at-infinity almost complex structures used in \cite{gps}. In particular, one cannot easily apply
the maximum principle type arguments, if possible at all.
This is the reason why \cite{gps} uses monotonicity arguments in their $C^0$-estimates, and utilized the class of
\emph{dissipative Hamiltonians} to establish $C^0$ estimates in defining symplectic cohomology for sectors.

In this part, we will introduce the notion of $\lambda$-sectorial almost complex structures,
which exist in abundance and allow for maximum principle arguments. Then we introduce the
sectorial Floer package by identifying the classes of Hamiltonians and of Lagrangian submanifolds
which correctly pair with the $\lambda$-sectorial almost complex structures.

\section{Set-up for the maximum principle}
\label{sec:C0estimates}

We start with a well-known definition in several complex variables (and its almost-complex
version). (See \cite{cieliebak-eliashberg}, for example.)

\begin{defn} \label{defn:J-convex} Let $(M,J)$ be an almost complex manifold.
\begin{enumerate}
\item
We call a function $\psi: M \to \RR$ \emph{$J$-convex} (resp. weakly $J$-convex)
if the Levi-form
$$
-d(d\psi \circ J)
$$
is positive definite (resp. positive semi-definite.)
When $\psi$ is an exhaustion function of $\Int M$ in addition,
we call it a \emph{$J$-convex exhaustion function}.
\item
We also call such a pair $(\psi,J)$ a \emph{pseudoconvex pair}.
\end{enumerate}
\end{defn}
Such a function is more commonly called a plurisubharmonic (resp. weakly plurisubharmonic)
function in several complex variable theory.

We recall some general definitions of quantitative measure of $J$-pinching.

\begin{defn}\label{defn:pinching-bound} Let $M$ be given a Riemannian metric $g$ and denote by $S(TM)$ the unit
disc bundle of $TM$. Let $U$ be a neighborhood of $\del M$.
Consider an almost complex structure $J$ and a two-form $\eta$ on $U$.  We define
$$
[\eta]_{(J;U)}: = \inf\{ \eta(v,Jv) \mid v \in S_x(T\del M), \, x \in U\}
$$
and call it the \emph{$J$-pinching lower bound} of $\eta$ on $U \supset \del M$.
\end{defn}

\begin{remark} We always assume that the metric is of bounded geometry, i.e., that it has bounded curvature and
its injectivity radius $\text{\rm inj}(g) \geq \delta$ for some $\delta > 0$.
If $J$ is tame with respect to $\omega$ on $U$, it is easy to see that $[\omega]_{(J; U)}> 0$.
\end{remark}

\begin{remark}
A similar notion was considered  in \cite{cieliebak-eliashberg} to
provide a quantitative measurement of pseudo-convexity for the Levi form when the
Levi form is $J$-positive. In our case, we apply the notion to an arbitrary 2-form.
\end{remark}

Applying Definition \ref{defn:J-convex} to $- \psi = \log \varphi$, we write the standard
notion of pseudoconvex boundary in several complex variable theory as follows.
\begin{defn}[Pseudoconvex-type boundary]
\label{defn. pseudoconvex type}
Let $(M,J)$ be an almost complex manifold with boundary $\del M$.
We say an almost complex structure $J$ on $M$ is of \emph{$\del M$-pseudoconvex type}
if there exists a function $\varphi: M \to \RR$ such that
\begin{itemize}
\item $\del M = \varphi^{-1}(0)$ and $\varphi > 0$ on $\nbhd(\del M) \setminus \del M$.
\item
There is a constant $\varepsilon_0$ such that
\eqn\label{eq:varphi-convexJ}
[d(d(\log \varphi \circ J))]_{(J;U_{\varepsilon_0})} \geq 0
\eqnd
where
$$
U_{\varepsilon_0} = \{0 < \varphi \leq \varepsilon_0\} .
$$
(I.e., $- \log \varphi$ is weakly $J$-convex on $U_{\varepsilon_0}$.)
\end{itemize}
When \eqref{eq:varphi-convexJ} is replaced by the pinching condition
$$
[d(\log \varphi \circ J)]_{(J;U_{\varepsilon_0})} \geq C
$$
for some constant $C \in \RR$, we call such \emph{$\del M$ $J$-pinched below} and $C$ the
\emph{$J$-pinching lower bound}.
\end{defn}

\begin{remark}
When $J$ is given, the existence of $\varphi$ depends only
on the germ of the neighborhoods of $\del \Mliou$. In fact, if $C > 0$ (resp. $C = 0$), being pinched below by $C$
agrees with the standard notion of strong pseudoconvexity (resp. (weak) pseudoconvexity)
of the boundary (relative to $J$).
\end{remark}

\begin{remark}
\begin{enumerate}
\item
If $\del M$ is compact, the existence of the constant $\varepsilon_0> 0$ is automatic
once there exists a $\del M$-pseudoconvex type $J$, and is
unnecessary to state as one of the defining conditions above.
However, when $\del M$ is noncompact, the existence of a uniform size $\varepsilon_0$
is not automatic and requires some good behavior of $\del M$ at infinity.
The Liouville sectors have this good behavior at infinity by its $Z$-completeness.
Moreover any cylindrical complex structure in the sense of
\cite{gps} is of $\del \Mliou$-pseudoconvex type \emph{outside a neighborhood the corner}.
\item On the other hand, we will see that the sectorial almost complex structure is
of $\del \Mliou$-pseudoconvex type on \emph{a full punctured neighborhood of $\del M$}.
\end{enumerate}
\end{remark}

The $\del M$-pseudoconvexity (or more generally $\del M$-pinching) for an $\omega$-tame almost
complex structure on symplectic manifolds $(M,\omega)$ has the following consequence.

\begin{prop}\label{prop:logIv}
Let $(M,\omega)$ be a symplectic manifold, possibly with corners. Equip it with a an almost-complex structure $J$ for which the Riemannian metric $g = \omega(\cdot ,J \cdot))$ is complete and has bounded geometry, but is {\emph not necessarily tame to $\omega$}. Finally, we also assume that $J$ is of $\del M$-pseudoconvex type (Definition~\ref{defn. pseudoconvex type}).

Suppose $v: \Sigma \to \Mliou$ is a $J$-holomorphic map. Then for the classical Laplacian $\Delta$, and for any function $\varphi: M \to \RR$ as in Definition~\ref{defn. pseudoconvex type}, we have
$$
\Delta (\log \varphi\circ v) \omega_\Sigma \geq C \, v^*\omega
$$
for some constant $C \geq 0$.
In particular, $-\log \varphi \circ v$ is a subharmonic function.
\end{prop}

\begin{proof}
For an arbitrary function $\varphi$, we compute
$$
- d(d(-\log\varphi\circ v) \circ j) = d(d(\log \varphi)J dv) = v^*(d(d(\log \varphi)J)).
$$
For any unit tangent vector $\xi \in T_z\Sigma$, we evaluate
$$
v^*(d(d(\log \varphi)J))(\xi,j \xi)
= d(d(\log \varphi)J) \left(dv(\xi), dv(j\xi)\right)
= d(d(\log \varphi)J) \left(dv(\xi), J dv(\xi)\right).
$$
Therefore Definition~\ref{defn:pinching-bound} tells us
\beastar
v^*(d(d(\log \varphi)J))(\xi,j \xi) & = & d(d(\log \varphi)J) \left(dv(\xi), J dv(\xi)\right)\\
& \geq & [d(d(\log \varphi)J)]_{(J;U_\epsilon)} |dv(\xi)|^2_g.
\eeastar
We have the standard comparison inequality between two metrics of bounded geometry\footnote{Recall that a Riemannian metric is said to have bounded geometry if (i) its injectivity radius is bounded away from zero, and (ii) the norms of the covariant derivatives of the curvature tensor are uniformly bounded.}, namely $g$ and a choice of $\omega$-tame metric $g_{J'} = \omega(\cdot, J'\cdot)$:
	\eqnn
	|dv(\xi)|_g^2
	\geq
	C' |dv(\xi)|_{g_{J'}}^2
	\eqnd
where $C'$ depends only on $g$ and $g_{J'}$.  We also have a standard inequality
	\eqnn
	|dv(\xi)|_{g_{J'}}^2 \geq C''' v^*\omega(\xi,j\xi).
	\eqnd
(See \cite[Proposition 7.2.3]{oh:book1}
for example.)  By the assumption that $J$ is
of $\del M$-pseudoconvex type, we may as well take $\varphi$ to be a function satisfying the properties in Definition~\ref{defn. pseudoconvex type}. Thus we have:
	\eqnn
	[d((\log \varphi)\circ J)]_{(J;U_\epsilon)}]
	\geq C''
	\eqnd
for some constant $C''$. By setting $C = C'C''C'''$, we have finished the proof.
\end{proof}

\begin{example}
For the cylindrical $J$ in the sense of \cite[Section 2.10]{gps} we may take $\varphi = R$ and witness a vanishing
$$
d(d(\log \varphi) \circ J) = 0
$$
on a neighborhood  $\nbhd^Z(\del M)$.
\end{example}

The following general results are the reason why we consider the interpolation
by a convex function $\varphi$.
\begin{prop}\label{prop:convex-interpolation} Let $\varphi: \RR^k \to \RR_+$
be any (weakly) convex function with $\frac{\del \varphi}{\del x_i} > 0$.
Let $g_i: M \to \RR$ be a smooth functions for $i = 1, \ldots, k$.
Let $J$ be an almost complex structure $J$ tame to a symplectic form $\omega$ of $M$.
Then the following hold:
\begin{enumerate}
\item
The function
$$
G: = \varphi(g_1,\ldots, g_k):\RR \to \RR
$$
satisfies
\begin{eqnarray}\label{eq:ddGJ}
-d(dG\circ J) & = & \sum_j\sum_i \frac{\del^2 \varphi}{\del x_j \del x_i}(g_1,\ldots,g_k))\,
(dg_j\circ J) \wedge dg_i \nonumber\\
& {} &- \sum_i \frac{\del \varphi}{\del x_i}(g_1,\ldots,g_k))\,
d(dg_i \circ J).
\end{eqnarray}
\item Assume all $g_i$ are $J$-plurisubharmonic (resp. $J$-weakly plurisubharmonic)
Then if $v: \Sigma \to M$ is $J$-holomorphic, then $f: = G\circ v: \Sigma \to \RR$
is $J$-plurisubharmonic (resp. $J$-weakly plurisubharmonic)
\end{enumerate}
\end{prop}
\begin{proof} We compute $dG = \sum_{i=1}\frac{\del \varphi}{\del x_i}(g_1,\ldots,g_k))\,  dg_i$ and so
$$
dG \circ J = \sum_{i=1}\frac{\del \varphi}{\del x_i}(g_1,\ldots,g_k)\,  dg_i \circ J.
$$
Therefore by taking the differential of this and a rearrangement, we obtain
$$
-d(dG\circ J) = \sum_j\sum_i \frac{\del^2 \varphi}{\del x_j \del x_i}(g_1,\ldots,g_k)\,
dg_i \circ J \wedge dg_j - \sum_i \frac{\del \varphi}{\del x_i}(g_1,\ldots,g_k)\,
d(dg_i \circ J)
$$
which is \eqref{eq:ddGJ}.

For the second statement, let $v: (\Sigma,j) \to (M,J)$ be $J$-holomorphic. Then we compute
\beastar
-d(d(G\circ v)\circ j) & = & -d(dG\circ dv)\circ j) =  -d(dG\circ J) \circ dv)
= - v^*d(dG\circ J)\\
& = & \sum_j\sum_i \frac{\del^2 \varphi}{\del x_j \del x_i}(g_1(v),\ldots,g_k(v))\,
v^*(d(g_i \circ J) \wedge dg_j) \\
&{}& - \sum_i \frac{\del \varphi}{\del x_i}(g_1(v),\ldots,g_k(v))\,
v^*d(dg_i \circ J).
\eeastar
Then the first summand is semipositive $(1,1)$-form (with respect to $J$) since the Hessian matrix
$\frac{\del^2 \varphi}{\del x_j \del x_i}$ is positive semidefinite while the second
summand is so because $\frac{\del \varphi}{\del x_i} > 0$ by definition,
and each $g_i$ is assumed to be plurisubharmonic (resp. weakly plurisubharmonic)
 with respect to $J_i$, and so $- d(dg_i \circ J)$ is a positive (resp. semipositive) $(1,1)$-form.
\end{proof}
Once ${\mathfrak s}_\varphi$ is proven to be an exhaustion function on $\nbhd(\del_\infty M \cup \del M)$,
if we find an almost complex structure $J$ for which $({\mathfrak s}_\varphi,J)$ becomes
a pseudoconvex pair, this proposition will ensure uniform $C^0$-bound for
the $J$-holomorphic curves.
\begin{remark} When we apply this proposition to our
circumstance, one essential complication arises when we consider the interpolation
$$
\varphi_{k+1}(R_1,\cdots, R_k,e^{-s})
$$
It is easy to define almost complex structures in a disjoint component of the sectorial corner
$$
\nbhd^Z(\del M) \setminus \nbhd(\del_\infty M)
$$
and in the ceiling corner
$$
\nbhd(\del_\infty M)\setminus \nbhd^Z(\del M)
$$
respectively on  with respect to which each $R_i$ and $s$ are $J$-plurisubharmonic.
An essential task then is to glue them without destroying the pseudoconvexity of
the pair $({\mathfrak s}_\varphi,J)$. Because the hypersurface $\{s^{-1}(c)\}$
(nontrivially) intersects sectorial corners at infinity and also because the collection of hypersurfaces
$$
\{R_i^{-1}(c_i)\}_{i=1}^k \bigcup \{s^{-1}(c)\}
$$
no longer forms a \emph{coisotropic} collection, this gluing process is a nontrivial task. We accomplish it in the next section.
\end{remark}

We will prove the following existence theorem in the next two sections, whose precise form
will be given along the way.
\begin{theorem}\label{thm:critically-$J$-convex} Let $(M,d\lambda)$ be a Liouville sectors with corners.
Then there exists a pseudoconvex pair $(\psi, J)$ in a neighborhood
$$
\nbhd(\del_\infty M \cup \del M)
$$
such that the function $\psi\circ u$ is also amenable to the strong maximum principle for
any $J$-holomorphic curves with boundary condition $u(\del \Sigma) \subset L$ with
$Z$-invariant-at-infinity Lagrangian $L$.
\end{theorem}
In fact, we will see that  any pair of
\beastar
(\psi,J) & = & (\text{\rm wiggled end profile function},\\
&{}& \quad \text{\rm its associated $\lambda$-sectorial almost complex structure})
\eeastar
is one such pair.

\begin{remark} The pair
 $$
(\psi,J)= (\text{\rm end profile function},\text{\rm its associated sectorial almost complex structure})
$$
is a pseudoconvex pair but does not satisfy the additional condition on the amenability of strong maximum
principle with $Z$-invariant Lagrangian submanifolds. However it does
for a different class of Lagrangians called \emph{gradient-sectorial} Lagrangian submanifolds
with respect to the end-profile function ${\frak s}_\varphi$ (see \cite{oh:gradient-sectorial}).
This class of Lagrangians is more convenient than the $Z$-invariant ones for the purpose
of studying the K\"unneth-type maps.
\end{remark}

\section{Sectorial almost complex structures}
\label{sec:sectorial-J}

In this section and the next, we introduce  new classes of almost complex structures which weakens the properties of both
 contact-type and cylindrical almost complex structures.
 The first class, which we
 construct in the present section
 and call \emph{sectorial almost complex structures} in \cite{oh:gradient-sectorial},
 enables us to apply the maximum principle near  both $\del_\infty \Mliou$ and $\del M$ (and hence establishing Gromov compactness).
 For this construction, \emph{the Liouville deformations constructed in the last
 section are not needed but only the cut-off functions $\kappa$
 will matter in its construction.}

 The second class, which we call \emph{$\lambda$-sectorial
 almost complex structures} whose construction will be completed
 in the next section, enables us to also apply the strong maximum principle in addition when they are paired with the $Z$-invariant
 Lagrangian submanifolds.

\subsection{Remarks on the common choices of $J$}
\label{section. common J choices}
As usual, let $s : \nbhd(\del_\infty \Mliou) \to \RR$ be a cylindrical coordinate on $\Mliou$ associated to a
symplectization.
In the literature it is common to choose an almost-complex structure $J$ satisfying conditions (C1) and (C2) below:
\begin{enumerate}
\item[(C1)] $J$ is an extension of an
almost complex structure on the contact distribution $\xi \subset \ker ds$
on each level set of $s$, and
\item[(C2)] In some neighborhood of $\del_\infty M$,
$$
J(Z) = X_\thetacontact, \quad J(X_\thetacontact) = - Z
$$
where $X_\thetacontact$ is the Reeb vector field for the contact form $\alpha$ on a level set of $s$.
\end{enumerate}
We note that $J$ satisfies the above conditions if and only if it satisfies the following:
\enum
\item[(C1')] In some neighborhood of $\del_\infty \Mliou$, $J$ is invariant under the flow of the Liouville vector field $Z$, and
\item[(C2')] In some neighborhood of $\del_\infty M$, $ds \circ J = - \lambda$.
\enumd
(The equation
$ds \circ J = -\lambda$ itself, without condition (C1'), is weaker than the combination of (C1) and (C2).) The above conditions allow for the use of maximum principles in proving
$C^0$-estimates for $J$-holomorphic curves.

The following definition follows by taking the half  of the
defining conditions of the standard notion of
$Z$-invariant $\lambda$-adapted almost complex structures in the literature
\emph{by dropping the $Z$-invariance}.

\begin{defn}[$J$ being contact type] Let $(M,\lambda)$ be a Liouville
manifold equipped with a symplectization radial function $s$
near infinity.  We call a tame almost complex structure $J$
\emph{$\lambda$-contact type}, if it satisfies
$$
-ds \circ J = \lambda.
$$
We denote by $\CJ_{\lambda}$ the set of all contact-type almost
complex structures on $(M,\lambda)$.
\end{defn}
Recall that $(F,\lambda_F)$
itself is assumed to be a Liouville manifold with a compact Liouville subdomain
$F_0$ such that $\supp f \subset F_0$, and equipped with
the radial function $s_F$.

The authors of~\cite{gps} instead work with $J$ that satisfy the following:
\begin{enumerate}
\item[(i)] $J$ is cylindrical (i.e., invariant under $Z$ near infinity), and
\item[(ii)] The projection $\pi_\C$ is holomorphic on a neighborhood of $\del \Mliou$.
\end{enumerate}
This last requirement is not compatible with the requirement for $J$ to be of contact-type in a neighborhood
of the ideal boundary $\del_\infty \Mliou$. This is one reason why ibid.---after requiring (i) and (ii)---must employ monotonicity arguments (as opposed to maximum principle arguments) to study pseudoholomorphic curves on the Liouville sector.

In our construction of sectorial almost complex structures we will interpolate a `contact type-like' condition near $\del_\infty M$ with suitable convexity requirements near $\del M$. (See Theorem~\ref{thm:interpolation}.) We will use the end-profiles function ${\mathfrak s}_\varphi$ for this purpose.

\begin{remark}
It would have been enough to consider almost complex structures $J$ for which
the pair $({\mathfrak s}_\varphi,J)$ is a pseudoconvex pair, i.e., those $J$ for which
$$
-d(d{\mathfrak s}_\varphi \circ J) \geq 0
$$
if we would like to apply maximum principle arguments
only for unperturbed $J$-holomorphic curves
without boundary.  All complications in the process of defining sectorial almost complex structures
arise because we would like to apply the maximum principle also \emph{to
Hamiltonian-perturbed} pseudoholomorphic curves
(see the proofs of Proposition \ref{prop:energy-identity} and Theorem \ref{thm:nonautonomous-confinement})
as well as to apply the strong maximum principle
to the case {\em with boundary}. (See the proof of Theorem \ref{thm:unwrapped}.)

Compare the simplicity of our definition of
sectorial Hamiltonians (Definition~\ref{defn:sectorial-H}) with the complexity of the definition of
the dissipative Hamiltonians used in \cite{gps} which are needed to achieve the confinement results via
the monotonicity arguments. In our sectorial framework, such complications are subsumed
in our careful geometric preparation of
the description of background geometry of Liouville sectors related to smoothing of
corners of the Liouville sectors leading to the definition of sectorial almost complex structures.
Once these are achieved, our effort is compensated by
the simplicity and naturality of the definition of $\lambda$-sectorial Hamiltonians
(see Definition \ref{defn:sectorial-H}).
\end{remark}

\subsection{Definition of $\kappa$-sectorial almost complex structures}
\label{subsec:kappa-sectorial}

At the end of the day, we came up with the following definition of $\kappa$-sectorial
almost complex structures as a tool for constructing the final class of almost
complex structures in Subsection \ref{subsec:lambda-sectorial-J} that are amenable both to the maximum principle and to
\emph{the strong maximum principle for the associated pseudoholomorphic curves
with the boundary condition of the $Z$-invariant Lagrangians.}

Given a splitting data
$$
U: = \nbhd(\del M) \cong F \times \C_{\text{\rm Re} \geq 0}, \quad \{(R,I)\},
$$
we have two natural symplectic foliations whose associated involutive distributions are given by
$$
\cD_F: = \ker d\pi_\C, \quad \cD_{\C}: = \ker d\pi_F,
$$
respectively.
They induce a split short exact sequence
$$
0 \to \cD_F \to TU \to \cD_\C \to 0
$$
on $U$. Its dual sequence is given by
$$
0 \to \cD_\C^\omega \to T^*M|_U \to \cD_F^\omega \to 0
$$
via the identification $T^*F \cong \cD_\C^\omega$ and $T^*\C \cong  \cD_F^\omega$.
We denote by $\cF_F$ the associated foliation whose leaves are given by
$$
F \times \{(x,y)\}, \quad (x,y) \in \C_{\text{\rm Re} \geq 0}.
$$
Using this we give the definition of $\kappa$-sectorial almost complex structures
in terms of $T^*M$ not $TM$ which is more useful to study the pseudoconvexity of
the pair $(\psi,J)$ in general.

%
%

\newenvironment{corner-J}{
	  \renewcommand*{\theenumi}{(J$\kappa$\arabic{enumi})}
	  \renewcommand*{\labelenumi}{(J$\kappa$\arabic{enumi})}
	  \enumerate
	}{
	  \endenumerate
}

\begin{defn}\label{defn:J-for-corners-kappa} Let $(M,\lambda)$ be a Liouville sector with boundary and corners.
An $\omega$-tame almost complex structure $J$
on a Liouville sector is said to be \emph{$\kappa$-sectorial} (with respect to the given smoothing profile)
if $J$ satisfies the following:
\begin{corner-J}
\item {\textbf{[$\cF_F$ is $J$-complex]}}\label{item. piF is holomorphic}
In a neighborhood of $\nbhd^Z(\del \Mliou)$
of $\del \Mliou$, we require
    \eqn\label{eq:J-versus-JF}
    J\left(T^*F \oplus 0_{\text{\rm span}\{d\mu_{i}, d\nu_{i}\}_{i=1}^k}\right)
    \subset T^*F \oplus 0_{\text{\rm span}\{d\mu_{i}, d\nu_{i}\}_{i=1}^k},
    \eqnd
    and $J$ restricts to an almost complex structure of contact-type on $F$.
\item {\textbf{[$d{\mathfrak s}_\varphi$ is $J$-dual to $\lambda_\kappa$]}}\label{item. ds is dual to lambda_kappa}
    In a neighborhood $\nbhd^Z(\del \Mliou)\cup \del_\infty M)$ of
    $\del \Mliou \setminus \nbhd(\del_\infty M)$, we have
    $$
   -d{\mathfrak s}_\varphi \circ J = \pi_F^*\lambda_F + \pi_\C^*\lambda_\C + d(\kappa f)
    $$
    for the given cut-off function $\kappa$ given in Definition \ref{defn:kappa}.
\end{corner-J}
\end{defn}
Here we refer to Proposition \ref{prop:gps} for the definition of $f$.

Note that Condition \ref{item. ds is dual to lambda_kappa} implies that the pair
$({\mathfrak s}_\varphi, J)$ forms a pseudoconvex pair since it implies
$$
-d(d{\mathfrak s}_\varphi \circ J) = d\lambda \geq 0
$$
as a $(1,1)$-current. In fact, the following simple definition is used in \cite{oh:gradient-sectorial}
for a construction of wrapped Fukaya category whose objects are \emph{not} the $Z$-invariant Lagrangians
but those which are invariant-at-infinity \emph{under the gradient flow of ${\mathfrak s}_\varphi$}:

\begin{defn}[Sectorial almost complex structures \cite{oh:gradient-sectorial}]
\label{defn:J-for-corners}
Let $(M,\lambda)$ be a Liouville sector with boundary and corners.
An $\omega$-tame almost complex structure $J$
on a Liouville sector is said to be \emph{$\lambda$-sectorial} (with respect to the given smoothing profile)
if $J$ satisfies the following:
\begin{enumerate}
\item {\textbf{[$\cF_F$ is $J$-complex]}}
In a neighborhood of $\nbhd^Z(\del \Mliou)$
of $\del \Mliou$, we require
    \eqn\label{eq:J-versus-JF}
    J\left(T^*F \oplus 0_{\text{\rm span}\{d\mu_{i}, d\nu_{i}\}_{i=1}^k}\right)
    \subset T^*F \oplus 0_{\text{\rm span}\{d\mu_{i}, d\nu_{i}\}_{i=1}^k},
    \eqnd
    and $J$ restricts to an almost complex structure of contact-type on $F$.
\item {\textbf{[$({\mathfrak s}_\varphi, J)$ is a pseudoconvex pair]}}
    In a neighborhood $\nbhd^Z(\del \Mliou)\cup \del_\infty)$ of
    $\del \Mliou \setminus \nbhd(\del_\infty M)$, we have
    $$
   -d(d{\mathfrak s}_\varphi \circ J) \geq 0.
    $$
\end{enumerate}
\end{defn}
We will not use this definition at all in the present paper since we still employ the
$Z$-invariant Lagrangian submanifolds as objects of wrapped Fukaya category $\Fuk(M)$.

\subsection{Existence of $\kappa$-sectorial almost complex structures}

\subsubsection{A coordinate change}

Instead of using $(R,I)$ coordinates for the complex numbers, we will utilize coordinates $(-\log R, I)$ (which of course changes the form of both $\lambda^\alpha_\C$ and of $\omega_\C$). These coordinates have the intuition that approaching $R = 0$ is akin to approaching $\infty$---not in the Liouville flow direction, but in the direction toward $\del M$. We will later justify this intuition by showing that, for sectorial $J$, if any member of a continuous family of $J$-holomorphic curve has image constrained to the region $R >0$, then every member of that family is so constrained.  See Remark~\ref{remark. confinement away from del M}.

\begin{notation}[$\mu, \nu$]
Let
\eqn\label{eq:coordinate change}
\mu = -\log R, \quad \nu = I.
\eqnd
\end{notation}

In what follows, we will describe almost-complex structures $J$ in terms of the cotangent bundle
$$
T^*(\nbhd(\del \Mliou) \setminus \del \Mliou) \approx T^*F \times T^*\C_{{\text{\rm Re}}> 0} = T^*F \oplus \text{\rm span}_{\RR}\{-dR, dI\}
$$
instead of the tangent bundle $T(\nbhd(\del \Mliou) \setminus \del \Mliou) $.

\begin{remark}\label{remark. linear algebra applies}
Note that $d\mu$ and $d\nu$ are non-vanishing and linearly independent at every
point of
$\nbhd(\del \Mliou) \setminus \del \Mliou.$
So in this region, we
can uniquely express any co-vector $\gamma$---in particular,
$d\mu \circ J$ and $d\nu\circ J$---
in the form
\eqn\label{eq:dRJ-dIJ-first}
d\mu \circ J =  \eta_\mu + a \, d\mu + b \, d\nu,
\quad d\nu \circ J = \eta_\nu + c \, d\mu + d \, d\nu
\eqnd
for some constants $a, \, b, \, c, \, d$ and $\eta_\mu, \, \eta_\nu \in T^*F$.
The reader should later compare this to~\eqref{eqn. G version of eta}.
\end{remark}

\begin{remark}\label{rem:F-coordinates}
The factorization $F \times \C_{\{\text{\rm Re} \geq 0\}}$
is crucial for the definition of sectorial almost complex structures although
the equation \eqref{eq:dRJ-dIJ-first} itself is an equation in the tangent distribution level and
does not require integrability of the subbundle $TF \subset T\del X|_F$ into a submanifold $F$.
See Definition \ref{defn:J-for-corners} \ref{item. piF is holomorphic},
 and the beginning of the proof of Theorem \ref{thm:J-abundance}.
\end{remark}

 We consider the change of coordinates of $(0,\infty) \times \RR \to \RR^2$
$$
\sigma: (0,\infty) \times \RR \to \RR^2 ; \quad \sigma(x,y) = (-\log x,y) =: (u,v).
$$
\begin{notation}[$\widetilde \pi_\C$]
On $\nbhd(\del M) \setminus \del M$, we have two projections
$$
\pi_\C = (R,I): M \to (0,\infty) \times \RR, \quad
\widetilde \pi_\C = (\mu,\nu): M \to \RR^2
$$
where we define
\eqn\label{eq:tilde-versus-nil}
\widetilde \pi_\C = \sigma \circ \pi_\C.
\eqnd
\end{notation}

Recall the Liouville form  $\lambda_\C^\alpha$  \eqref{eq:lambda-alpha} and restrict it to
$\C_{{\text{\rm Re}}>0} \cong (0,\infty) \times \RR$. Denoting the standard coordinates of $\RR^2$ by
$(u,v)$, we denote the pushforward form of $\lambda_\C$ to $\RR^2$ by
    $$
    	\widetilde \lambda_\C^\alpha
		:= e^{-u}((1-\alpha) dv + \alpha v du).
    $$
Then we have
\eqn\label{eq:-tilde-pi-lambda}
    	\widetilde \pi_\C^* \widetilde \lambda_\C^\alpha
		= e^{-\mu}((1-\alpha) d\nu + \alpha \nu d\mu)
\eqnd
on $\nbhd(\del M) \setminus \del M$.

\begin{remark}
Note that $
		\pi_\C^* \lambda_\C^{\alpha}=\widetilde \pi_\C^* \widetilde \lambda_\C^\alpha$ on $\nbhd(\del M) \setminus \del M$. The tilde notation connotes both the restriction on the domain, and our preference to work with $(\mu, \nu)$ coordinates rather than $(R,I)$ coordinates.
\end{remark}

The following will be useful later  to guarantee the existence of certain almost-complex structures
when applying the consistency criterion laid out in Condition \ref{cond:consistency_for_J}.

\begin{prop} \label{prop. linear independence of forms}
On $\nbhd(\del M) \setminus M$,
\enum[(i)]
    \item For all $0 < \alpha < 1$,
    $\widetilde \pi_\C^* \widetilde \lambda_\C^\alpha$ is linearly independent of $d\mu$.
    (Though we will not utilize this:
    $\widetilde \pi_\C^* \widetilde \lambda_\C^\alpha$ is also linearly independent of $d\nu$ if we further assume $\nu \neq 0$.)
    \item  Let $0 < \alpha < 1$ and consider the function $h_\C^{\alpha}$ from \eqref{eq:hCCalpha}.
    Then
    	$d(\widetilde \pi_\C^* (\widetilde h_\C^{\alpha}) )$
    is linearly independent of $\widetilde \pi_\C^*\widetilde \lambda_\C^{\alpha}$ and $d\nu$.
    (Though we will not utilize this:
    $d(\widetilde \pi_\C^* (\widetilde h_\C^{\alpha})$ is also linearly independent of $d\mu$ if
    we further assume $\nu \neq 0$.)
\enumd
\end{prop}

\begin{proof}
(i) is obvious from \eqref{eq:-tilde-pi-lambda}.
For (ii), we compute
\beastar
        d(\widetilde \pi_\C^* (\widetilde h_\C^{\alpha})
        &= &\frac12 d\left((1-\alpha) e^{-2\mu} + \alpha \nu^2 \right) \\
        &= & - (1-\alpha) e^{-2\mu} d\mu + \alpha \nu d\nu.	
\eeastar
Comparing this with \eqref{eq:-tilde-pi-lambda}, the lemma follows.
\end{proof}

\subsubsection{Some linear algebra (a consistency check)}
\label{subsubsec:consistency}

We begin with a linear algebra computation.

\begin{notation}\label{notation. linear algebra}
Let $V$ be a finite-dimensional real vector space and equip $V$ with a splitting $V = V_F \oplus V_C$. We assume $V_C$ is two-dimensional, and we further fix a basis $\{u, v\}$ of $V_C$. If $G$ is an endomorphism of $V$, let us write
	\eqn\label{eqn. G version of eta}
	Gu = \eta_u + au + bv,
	\qquad
	Gv = \eta_v + cu + dv
	\eqnd
where $\eta_u$ and $\eta_v$ are the projections of $Gu$ and $Gv$, respectively, to $V_F$.
\end{notation}

The following shows that a block-upper-triangular almost-complex structure $G$ is uniquely extendable from $V_F$ by specifying one of $Gu$ or $Gv$.

\begin{prop}\label{prop. linear algebra}
The assignments
	\eqnn
	G \mapsto (G|_{V_F}, \eta_u, a, b),
	\qquad
	G \mapsto (G|_{V_F}, \eta_u, au+bv),
	\eqnd
define bijections between
	\enum
	\item The set of $\RR$-linear endomorphisms $G: V \to V$ satisfying $G^2= -1$ and $G(V_F \oplus 0) \subset V_F\oplus 0$.
	\item The set of quadruplets $(A, \eta, x, y)$ where $A$ is an endomorphism of $V_F$ satisfying $A^2 = -1$, $\eta$ is an element of $V_F$, $x$ is an arbitrary real number, and $y$ is a {\em non-zero} real number.
	\item The set of triples $(A,\eta, w)$ where $A$ and $\eta$ are as above, and $w$ is a vector of $V_C$ linearly independent from $u$.
	\enumd
Likewise, the assignments $	G \mapsto (G|_{V_F}, \eta_v, c, d)$ and $G \mapsto (G|_{V_F}, \eta_u, cu+dv)$ define bijections between the set of almost-complex structures on $V$ preserving $V_F$, and
	\begin{enumerate}
	\item[$(2)'$] The set of quadruplets $(A,\eta,x,y)$ where now only $x$ is required to be non-zero.
	\item[$(3)'$] The set of triples $(A,\eta, w)$ where $A$ and $\eta$ are as above, and $w$ is a vector of $V_C$ linearly independent from $v$.
	\end{enumerate}
\end{prop}

\begin{proof}
Assume $G$ is in the set (1) and write $G$ in matrix form as
	\be\label{eq:GJ}
	G =
	\left(
	\begin{array}{ccc}
	A & \eta_u & \eta_v \\
	0 & a & c \\
	0 & b & d
	\end{array}
	\right).
	\ee
We first claim that the data $(A,\eta_u,a,b)$  uniquely determine the third column. This proves that the assignment from (1) to (2) is an injection.

To prove the claim, note that $G^2 = -1$ implies that the lower-right 2-by-2 matrix also squares to -1. From this we deduce
	\begin{itemize}
	\item $a^2 + bc = -1$, and hence $b$ and $c$ must be non-zero (we have assumed that $V$ is a real vector space) and
	\item $a = -d$.
	\end{itemize}
(Indeed, these two requirements are equivalent to demanding that the 2-by-2 matrix squares to -1.) So the pair $(a,b)$ determines the pair $(c,d)$.

To finish the proof of the claim, we note that the two constraints
	\eqnn
	A\eta_u + a\eta_u + b\eta_v
	= 0
	= A \eta_v + c\eta_u + d\eta_v
	\eqnd
implies that a choice of $\eta_u$ uniquely determines $\eta_v$ and vice versa by the formulae
	\eqn\label{eqn. etas}
	\eta_v = -{\frac 1 b}  (A+a) \eta_u
	\qquad
	\eta_u = -{\frac 1 c} (A+d)\eta_v
	\eqnd
(keeping in mind that $bc \neq 0$).

To show that the map from (1) to (2) is a surjection, one needs only check that the two equations in~\eqref{eqn. etas} are consistent; this is straightforward. Noting that the condition ``$b \neq 0$'' is equivalent to the condition ``$w = au + bv$ is linearly independent from $u$,'' we see that the map from (1) to (3) is also a bijection.

That the maps from 1. to $(2)'$. and $(3)'$. are bijections follows analogously.
\end{proof}

\begin{prop}\label{prop:consistency-for-J}
Fix an almost-complex structure on $T^*F$ and a smooth 1-form $\tau$ whose $T^*\C$ component is linearly independent from $d\mu$. Then, on $\nbhd(\del \Mliou) \setminus \del \Mliou$, there exists a $T^*F$-preserving almost-complex structure $J$ for which $d\mu \circ J = \tau$.
Alternatively, if the $T^*\C$ component of $\tau$ is linearly independent of $d\nu$, on $\nbhd(\del \Mliou) \setminus \del \Mliou$ there exists a $T^*F$-preserving almost-complex structure $J$ for which $d\nu \circ J = \tau$.
\end{prop}

\begin{proof}
Apply Remark~\ref{remark. linear algebra applies} and Proposition~\ref{prop. linear algebra} by setting
	\eqnn
	V = T^*_{x}M,
	\qquad
	V_F = T^*_{\pi_F(x)}F,
	\qquad
	V_C = T^*_{\pi_\C(x)}\C,
	\qquad
	G = J_x .
	\eqnd
\end{proof}

\subsubsection{Strategy of the proof of existence}

\begin{notation}[$N,\epsilon,\varepsilon_0$]\label{nota:Nee0}
We consider mutually disjoint neighborhoods
\be\label{eq:near-infinity}
\nbhd_{2\varepsilon_0}(\del \Mliou)\setminus \nbhd(\del_\infty M) \subset
\del \Mliou \setminus \nbhd(\del_\infty M)
\ee
and
\eqn\label{eq:near-boundary}
\nbhd_{s\geq N}(\del_\infty M)\setminus \nbhd^Z_{\epsilon'}(\del \Mliou)
\eqnd
for some $N > 0$ and $0 < \epsilon' < 2 \varepsilon_0$
which will be determined at the end of the proof of
Proposition \ref{prop:linear-independence}.
\end{notation}

%

Then we define $J$ on these neighborhoods by putting the following conditions:
\begin{cond}\label{cond:consistency_for_J}
Consider $J$ on the above two neighborhoods \eqref{eq:near-infinity} and \eqref{eq:near-boundary}.
\begin{enumerate}
    \item [{(J1)}] In a neighborhood of $\nbhd^Z(\del \Mliou)$ of $\del \Mliou$, using the
    splitting~\eqref{eq:splitting}, we require
    \eqn\label{eq:J-versus-JF}
    J\left(T^*F \oplus 0_{\text{\rm span}\{d\mu_{i}, d\nu_{i}\}_{i=1}^k}\right)
    \subset T^*F \oplus 0_{\text{\rm span}\{d\mu_{i}, d\nu_{i}\}_{i=1}^k},
    \eqnd
    and $J$ restricts to an almost complex structure $J_F$ of contact-type on $F$.
\item [{(J2a)}]
    In a neighborhood $\nbhd_{2\varepsilon_0}(\del \Mliou)\setminus \nbhd(\del_\infty M)$ of
    $\del \Mliou \setminus \nbhd(\del_\infty M)$, we have
    $$
    - d\mu_i\circ J = \pi_{\C,i}^*\lambda_\C + \pi_F^*\lambda_F
    $$
for $i = 1, \cdots, k$.
\item [{(J2b)}] In a neighborhood
$\nbhd_{s\geq N}(\del_\infty M)\setminus \nbhd^Z_{\epsilon}(\del \Mliou)$ for some $N > 0$ and $0 < \epsilon < \varepsilon_0$,
we have
    $$
    - ds \circ J = \lambda
    $$
\item [{(J2c)}] In a neighborhood
$$
\nbhd(\del_\infty M) \cap \nbhd_{\varepsilon_0}^Z(\del \Mliou)
$$
    of the ceiling of each sectorial corner in $\del_\infty M \cap \del M$,  we have
    $$
    - ds_{k+1,\varphi} \circ J = \lambda_\kappa
    $$
for the given cut-off function $\kappa:[0,\infty) \to [0,1]$.
\end{enumerate}
\end{cond}

The following is the key result in which the special properties of the barrier
functions $\varphi$ crucially enters and enables us to interpolate the
almost complex structure of `contact type-like' near $\del_\infty M$ and the one
of `cylindrical type-like' on $F \times \C_{\text{\rm Re}\geq 0}$
used in \cite{gps} near $\del M$.

\begin{theorem}[Interpolation] \label{thm:interpolation}
Let $\kappa$ be the function  of the type \eqref{eq:kappa}
appearing in Condition (J2c) above and write
\eqn\label{eq:lambda-kappa}
\lambda_\kappa = \pi_\C^*\lambda_\C + \pi_F^*\lambda_F + d(\kappa f).
\eqnd
Then
\begin{enumerate}
\item We can choose $\varphi$ so that $\varepsilon_0 \leq \mathsf{ht}(\varphi) < \frac32 \varepsilon_0 < 2 \varepsilon_0$.
\item
There exists a sufficiently small $0 < \epsilon_1 \leq \frac{\varepsilon_0}{8}$ and a function $\kappa$ such that
Condition (J2c) is admissible, i.e., that there exists some $J$ that satisfies
$$
- ds_{k+1,\varphi} \circ J = \lambda_\kappa
$$
and is compatible with Conditions (J1) and (J2a)-(J2c) for all $\varphi$ satisfying
 $0 < \sqrt{2} \mathsf{wd}(\varphi) < \epsilon_1$.
\end{enumerate}
\end{theorem}

The proof of this theorem will occupy the next entire subsection.

\subsubsection{Proof of compatibility}

We start from Condition (J1)
by choosing a compatible $J_F \in \CJ_{\lambda_F}$ on $F$.
For the simplicity of exposition, we assume $k=1$.
For the case $k > 1$, not much is different and so left to the readers.

We decompose $ds$ into $T^*\C$ and $T^*F$ components in
$T^*M = T^*\C \oplus T^*F$
\beastar
ds & = & ds^\C \oplus ds^F, \quad ds^\C = ds|_{T\C \oplus \{0\}},
\quad ds^F: = ds|_{TF\oplus \{0\}}\\
\lambda_\kappa & = & \lambda_\kappa^\C \oplus \lambda_\kappa^F, \quad \lambda_\kappa^\C
= \lambda_\kappa|_{T\C \oplus \{0\}},
\quad \lambda_\kappa^F: = \lambda_\kappa|_{TF\oplus \{0\}}.
\eeastar
Then we recall from Remark \ref{rem:varphi-epsilon-k}
$$
\varphi(x_1, x_2) = \begin{cases} x_1 \quad &\text{for } \,  x_2 \geq \mathsf{ht}(\varphi),
\, x_1 \leq \mathsf{wd}(\varphi) \\
x_2 \quad &\text{for } \,  x_1 \geq \mathsf{ht}(\varphi), \, x_2 \leq \mathsf{wd}(\varphi)
\end{cases}
$$
and $\varphi(x_1, x_2)$ smoothly interpolates the two functions $x_1, \, x_2$.
Therefore in the rest of this section, we will make our discussion on the following
region
\bea\label{eq:outside-corner}
&{}& \{0 \leq  R \leq \mathsf{wd}(\varphi), \, 0< e^{-s} \leq \mathsf{ht}(\varphi)\}
\\
&\bigcup&  \{0 \leq e^{-s} \leq \mathsf{wd}(\varphi), \, 0 \leq R \leq \mathsf{ht}(\varphi)\}
\eea
focusing on the first summand region utilizing the
$\iota$-reflection symmetry of $\varphi$.

Then it remains to verify that the condition (J1)  can be continuously
deformed
\begin{enumerate}
\item to the rest of the interior so that the condition
becomes also compatible to Condition (J1),
\item and to a neighborhood of the ceiling corner
$$
\nbhd(\del_\infty M \cap \del M).
$$
\end{enumerate}

 We decompose
\be\label{eq:lambda-components}
\pi_\C^*\lambda_\C + \pi_F^*\lambda_F + d(\kappa f)
=  (\pi_\C^*\lambda_\C + d_\C (\kappa\, f))
+  (\pi_F^*\lambda_F + \, d_F(\kappa f)).
\ee
and
$$
ds_{1+1,\varphi} = \frac{1}{\varphi} \left(\frac{\del \varphi}{\del x_1}(R,e^{-s})dR
 - \frac{\del \varphi}{\del x_2}(R,e^{-s})e^{-s}\, ds \right)
$$
into
\be\label{eq:ds-components}
ds_{1+1,\varphi} = d_F s_{1+1,\varphi} + d_\C s_{1+1,\varphi} =: w_F + w_\C.
\ee

We first prove the following proposition.

\begin{prop}\label{prop:linear-independence} Let $0 < \alpha \leq 1$.
If we choose $N_0, \, N_1$ so that both $N_0$ and
$N_1 - N_0$ sufficiently large and that $F_0 \subset F$ is sufficient large,
then the two vectors
\be\label{eq:two-vectors}
\frac{\del \varphi}{\del x_1}(R,e^{-s})dR
- \frac{\del \varphi}{\del x_2}(R,e^{-s}) e^{-s} ds^\C,\quad
-\pi_\C^*\lambda_\C- d_\C (\kappa f)
\ee
are linearly independent on
a neighborhood $\nbhd(\del_\infty M \cup \del M)$.
\end{prop}
\begin{proof}
 We have already taken care of the region
$$
\nbhd(\del_\infty M \cup \del M) \setminus \nbhd(\del_\infty M \cap \del M)
$$
by the requirement imposed in Condition \ref{cond:consistency_for_J} (J2a) and (J2b). Therefore we
 will focus on the corner region $\nbhd(\del_\infty M \cap \del M)$
 concerned by (J2c) therein.

\medskip

{\bf{Step 1 (On $\nbhd(\del_\infty M \cap \del M) \bigcap (F\setminus F_0) \times \C_{\text{\rm Re}\geq 0}$):}}
\smallskip

We first consider the region $(F\setminus F_0) \times \C_{\text{\rm Re} \geq 0}$,
where we have $f = 0$ and so
$ \lambda_\kappa = \pi_F^*\lambda_F + \pi_\C^*\lambda_\C$.
Recall that we have made the choice
$$
 s = \beta \pi_F^*s_F  + (1- \beta) \pi_\C^*s_\C
$$
for some sufficiently small $\beta> 0$ (depending on $\alpha$) in \eqref{eq:choice-s}.
We also have $\lambda_\kappa = \lambda$.
This implies
$$
\lambda_\kappa^\C = \pi_\C^*\lambda_\C = (1-\alpha) R\, dI - \alpha I\, dR.
$$
We also recall from \eqref{eq:sC} that  $h_\C
= \frac12 ((1-\alpha)R^2 + \alpha I^2)$ and $s_\C = \log h_\C$.
Therefore we compute
\be\label{eq:dsC}
(ds)^\C = (1-\beta) ds_\C =(1- \beta) \frac{(1-\alpha) R}{h_\C} \, dR
+(1- \beta) \frac{\alpha I}{h_\C}\, dI.
\ee
On the other hand, we have
\bea \label{eq:LHS-2}
&{}& \frac{\del \varphi}{\del x_1}(R,e^{-s})dR  - \frac{\del \varphi}{\del x_2}(R,e^{-s})e^{-s} (ds)^\C
\nonumber\\
& = & \left(\frac{\del \varphi}{\del x_1}(R,e^{-s})
- \frac{\del \varphi}{\del x_2}(R,e^{-s})e^{-s} (1 - \beta)
\frac{\del s_\C}{\del R}e^{-s} \right) dR\nonumber\\
&{}& -  e^{-s} \frac{\del \varphi}{\del x_2}(R,e^{-s})(1 - \beta) \frac{\del s_\C}{\del I }\, dI\nonumber \\
& = & \left(\frac{\del \varphi}{\del x_1}(R,e^{-s})  -
\frac{\del \varphi}{\del x_2}(R,e^{-s})e^{-s} (1 - \beta) (1-\alpha)R/h_\C \right)dR\nonumber \\
&{}& \quad
- e^{-s} \frac{\del \varphi}{\del x_2}(R,e^{-s})\frac{ (1 - \beta) \alpha  I}{h_\C}\, dI.
\eea
  We also recall from \eqref{eq:phi-derivatives}
$$
0 < \frac{\del \varphi}{\del x_1}, \quad \frac{\del \varphi}{\del x_2} \leq 1.
$$
To check the linear independence of the two vectors \eqref{eq:two-vectors}, we evaluate the determinant
$$
\left|\begin{matrix} \alpha I  & \frac{\del \varphi}{\del x_1}
 - (1 - \beta)(1-\alpha)\frac{\del \varphi}{\del x_2} e^{-s}R/h_\C\\
 - (1-\alpha)R & -
(1 - \beta) \frac{\del \varphi}{\del x_2}e^{-s} \alpha I/h_\C
\end{matrix}
\right| =: D
$$
\emph{keeping the insertion of $(R,e^{-s})$ into $\frac{\del \varphi}{\del x_1}, \, \frac{\del \varphi}{\del x_2}$ in our mind.}  We compute
\bea\label{eq:D3}
D & =  &- (1 - \beta)\frac{ (\alpha I)^2}{h_\C}  \frac{\del \varphi}{\del x_2} e^{-s}
+ (1-\alpha)R  \frac{\del \varphi}{\del x_1}
- (1 - \beta) \frac{((1-\alpha)R)^2}{h_\C} \frac{\del \varphi}{\del x_2} e^{-s} \nonumber\\
& = &-(1 - \beta)\frac{\left( (\alpha I)^2  + ((1-\alpha)R)^2\right)}{h_\C}
 \frac{\del \varphi}{\del x_2} e^{-s}
+ (1-\alpha)R  \frac{\del \varphi}{\del x_1}.
\eea
Noting that
$$
\frac{\del \varphi}{\del x_2}e^{-s} = x_2 \frac{\del \varphi}{\del x_2}\Big|_{(x_1,x_2)= (R,e^{-s})},\quad
R \frac{\del \varphi}{\del x_1} = x_1 \frac{\del \varphi}{\del x_1}
\Big|_{(x_1,x_2) = (R,e^{-s})}
$$
we would like to apply Lemma \ref{lem:x2dphi>x1dphi} on the box region $B_{\epsilon_1,\varepsilon_0}$.
For this purpose, we consider the region
$$
\nbhd_{s \geq N_2}(\del_\infty M) \cap \{0 \leq R \leq \frac{3\varepsilon_0}{2}\}
$$
for  $N_2> 0$ sufficiently large
so that
\be\label{eq:epsilon1}
0 < e^{-s} < 2 \sqrt{\epsilon_1} \Longleftrightarrow s > N_2;
\quad N_2 = - \log 2 \sqrt{\epsilon_1}.
\ee
(We will also require $\epsilon_1$ to satisfy
$$
2T_0 \sqrt{\epsilon_1} \leq  \frac{3 \varepsilon_0}{2}.
$$
See Lemma \ref{lem:x2dphi>x1dphi}.)
Then we have
$$
\frac{\del \varphi}{\del x_2}e^{-s} \geq R \frac{\del \varphi}{\del x_1} \geq 0
$$
by Lemma \ref{lem:x2dphi>x1dphi} and hence obtain the inequality
\beastar
D  &\leq & \left( (1-\alpha) - (1-\beta) \frac{\left((\alpha I)^2
+ ((1-\alpha)R)^2\right)}{2h_\C}\right) \frac{\del \varphi}{\del x_2} e^{-s}\\
& = & ((1-\alpha) - (1-\beta))\frac{\del \varphi}{\del x_2} e^{-s}
= -(\alpha - \beta)\frac{\del \varphi}{\del x_2} e^{-s}
< 0
\eeastar
from \eqref{eq:D3}, if we choose $\beta < \alpha$.
Therefore provided
\be\label{eq:choice-beta2}
\beta < \min\left\{\alpha, \frac1{C_{f,N_0}}\right\}.
\ee
the two vectors \eqref{eq:two-vectors}
are linearly independent for all $0< \alpha \leq 1$ on the intersection
$$
(F \setminus F_0) \times \C_{\text{\rm Re}\geq 0} \cap \{s \geq N_2\}
$$
where we made the choice for $s$ as in \eqref{eq:choice-s}.

\medskip

{\bf Step 2 (On $\nbhd(\del_\infty M \cap \del M) \bigcap F_0 \times \C_{\text{\rm Re}\geq 0}$):}
\smallskip

Recalling from Condition \ref{cond:s} that we have chosen $F_0$ so that
\beastar
&{}& F_0 \times \{(R,I) \mid |I| \leq N_0, \, 0 \leq R \leq 2\varepsilon_0\} \\
& \subset&  M \setminus
(\nbhd_{s \leq N_2}(\del_\infty M)\cup \nbhd_{R\leq \epsilon_1}(\del M),
\eeastar
this region is out of the neighborhood of $\del_\infty M \cup \del M$
of our main interest for the study of linear independence.
Therefore we have only to consider the region
$$
F_0 \times \{|I| \geq N_0\} \cap \{s > N_2\}.
$$
Since $f$ does not depend on $R, \, I$ when $|I| \geq N_0$,
we still have
$$
\lambda^\C = \pi_\C^* \lambda_\C
$$
where we also made the same choice for $s$ as in \eqref{eq:choice-s}.
Then exactly the same proof as done in Step 2 applies to this case too.
\end{proof}

Now we go back to the proof of Theorem \ref{thm:interpolation}.

Using the matrix representation $G_J$ of $J$ given in \eqref{eq:GJ}, we write
$$
d_F s_{1+1,\varphi}  = d_F s_{1+1,\varphi} + d_\C s_{1+1,\varphi}; \quad
d_F s_{1+1,\varphi}=: w_F, \, d_\C s_{1+1,\varphi}=
\left(\begin{matrix} u_1\\u_2 \end{matrix} \right),
$$
and evaluate
\beastar
- ds_{1+1,\varphi} \circ J = \left(\begin{matrix} A & \eta_u & \eta_v \\
0 & a & b\\
0 & c & d \end{matrix}\right)
\left(\begin{matrix} w_F \\ u_1 \\ u_2 \end{matrix}\right)
= \left(\begin{matrix} A w_F + u_1\eta_u + u_2 \eta_v \\
\left(\begin{matrix} a & b \\ c & d \end{matrix}\right)
\left(\begin{matrix} u_1\\u_2 \end{matrix} \right)
\end{matrix}
\right)
\eeastar
Therefore the equation $-ds_{1+1,\varphi} = \lambda_\kappa \circ J$ for $J$ is reduced to
\be\label{eq:final-equation}
\begin{cases}
Aw_F + u_1 \eta_u + u_2 \eta_v = \lambda_\kappa^F\\
\left(\begin{matrix}a & b \\ c & d \end{matrix}\right)
\left(\begin{matrix} u_1\\u_2 \end{matrix} \right) = \lambda_\kappa^\C
\end{cases}
\ee
for the triple $(A,\eta_u,\eta_v)$.
\begin{lemma}
We can solve
the second equation for $(u_1,u_2)$ for the vector $(u_1,u_2) \in \R^2$
that is linearly independent of $\lambda_\kappa^\C$.
\end{lemma}
\begin{proof} Once the equation is solved, the second statement is obvious since
$$
\left(\begin{matrix}a & b \\ c & d \end{matrix}\right)^2 = - Id.
$$
By definition, the vector $(u_1,u_2)^T$ is
nothing but the coordinate expression of the  one-form
\bea\label{eq:vecu-in-coordinates}
u_1\, dR + u_2\, dI &: =& \frac{\del \varphi}{\del x_1}(R,e^{-s})dR
 - \frac{\del \varphi}{\del x_2}(R,e^{-s})e^{-s} ds^\C  \nonumber \\
 & = & \left( \frac{\del \varphi}{\del x_1}(R,e^{-s})
  - \frac{\del \varphi}{\del x_2}(R,e^{-s})e^{-s} \frac{\del s}{\del R}\right)\, dR
  \nonumber\\
&{}&  - \frac{\del \varphi}{\del x_2}(R,e^{-s})e^{-s} \frac{\del s}{\del I} \,dI
   \eea
in terms of the basis $\{dR,dI\}$. Then we compute
\be\label{eq:equation-u1u2}
\left(\begin{matrix}a & b \\ c & d \end{matrix}\right)
\left(\begin{matrix} u_1\\u_2 \end{matrix} \right) = \left(\begin{matrix}
a u_1 + b u_2 \\
c u_1 + d u_2 \end{matrix}\right) = \left(\begin{matrix}
a u_1 + b u_2 \\
c u_1 - a u_2 \end{matrix}\right).
\ee
Therefore it is enough to show that for any fixed nonzero vector $\left(\begin{matrix} u_1 \\ u_2 \end{matrix} \right)$ we have
\be\label{eq:ampleness}
\left\{\left(\begin{matrix}
a u_1 + b u_2 \\
c u_1 - a u_2 \end{matrix}\right) \, \Big|\, a^2 + bc = -1\right \} \supset \R^2
\Big \backslash \span_\R \left\{\left(\begin{matrix} u_1 \\ u_2 \end{matrix} \right)\right\}.
\ee
By a change of basis by an element of $GL_+(2,\R)$, it is enough to consider the case when
$(u_1,u_2) = (1,0)$. But in this case, the resulting vector
in \eqref{eq:equation-u1u2} becomes
$\left(\begin{matrix} a \\ c \end{matrix} \right)$.
By the only constraint $-bc = a^2 +1$ imposed on $a, \, c$, it follows that the above span
becomes
$$
\R^2 \setminus \span_\R \left\{\left(\begin{matrix} 1 \\ 0\end{matrix} \right)
\right\}.
$$
This finishes the proof.
\end{proof}
Now it remains to solve the first equation of \eqref{eq:final-equation}
\beastar
A w_F + u_1 \eta_u + u_2 \eta_v =  \pi_F^*\lambda_F
\eeastar
on $T^*F$ for the triple $(A, \eta_u, \eta_v)$ with $\eta_u, \, \eta_v \in T^*F$
and $A^* \in \CJ_{\omega_F}$.  (Here $A^*$ is the fiberwise complex structure
 on $TF$ induced from $A$ acted upon $T^*F$.)

 This implies $\lambda_\kappa^F = \pi_F^*\lambda_F$.
Since $F$ itself is a Liouville manifold, we can find a contact-type almost
complex structure $J_F$ that satisfies $-ds_F \circ J_F = \lambda_\C$.
Recalling
$$
w_F = - d_F s_{1+1,\varphi} = - \frac1{\varphi}\frac{\del \varphi}{\del x_2} e^{-s} d_F s
$$
and $\frac1{\varphi}\frac{\del \varphi}{\del x_2} \neq 0$ on $\nbhd_{R \leq \frac{3}{2}\varepsilon_0}(\del M)$, we can solve the equation by requiring $J$ to satisfy
$$
\begin{cases}
u_1 \eta_u + u_2\eta_v = d_F(\kappa f)\\
A w_F =  \pi_F^*\lambda_F.
\end{cases}
$$
Obviously we can solve the first equation for $(\eta_u,\eta_v)$
since we have $(u_1, u_2) \neq (0,0) $ and $d_F(\kappa f) = \kappa df \in T^*F \oplus \{(0,0)\}$.

The second equation can be solved by requiring $A$ to satisfy
$$
\begin{cases}
A(d_Fs) = - e^s \left(\frac1{\varphi}\frac{\del \varphi}{\del x_2}\right)^{-1} \lambda_F \\
A \lambda_F =  e^{-s} \frac1{\varphi}\left(\frac{\del \varphi}{\del x_2}\right) d_F s
\end{cases}
$$
on the subspace $\span_\R \{d_F s, \lambda_F\} \subset T^*F$ by decomposing
$$
T^*F = \span_\R \{d_F s, \lambda_F\} \oplus (\span_\R \{d_F s, \lambda_F\})^{\omega_F}.
$$
In other words, we require $A$ to satisfy
$$
\begin{cases}
A(d_Fs) =  (r \, J_F)(d_Fs)\\
A \lambda_F =  (\frac1r\, J_F)(\lambda_F)
\end{cases}
; \quad r : = - e^s \left(\frac1{\varphi}\frac{\del \varphi}{\del x_2}\right)^{-1}
$$
on $\span_\R \{d_F s, \lambda_F\} \subset T^*F$ which obviously defines a (sub) complex structure on
the subspace.
Such a choice is possible on the corner neighborhood
$$
\nbhd_{s \geq N_0}(\del_\infty M) \cap \nbhd_{R \leq \frac{3 \varepsilon_0}{2}}(\del M)
$$
which is the region of our current interest where the radial function is defined and
$\span_\R\{d_F s(y), \lambda_F(y)\}$ forms a 2 dimensional $\omega_F$-symplectic subspace of $T_y^*F$.

Finally combining this with Proposition \ref{prop:consistency-for-J} and Proposition \ref{prop:linear-independence},
we have finished the proof of Theorem \ref{thm:interpolation} by setting
$\epsilon_1$  given in Notation \ref{nota:Nee0} and
$N_0$ for the constant $N_0$ appearing in Condition \ref{cond:s},
and $\beta$ as in \eqref{eq:choice-beta2}.

\section{Construction of $\lambda$-sectorial almost
complex structures}
\label{sec:deformation-lambda}

In this section, we introduce the notion of $\lambda$-sectorial
almost complex structures.  We prove their existence by
utilizing a Liouville diffeomorphism $\phi_\kappa$ to be
constructed in the present section
to transform the $\kappa$-sectorial almost complex structures
into $\lambda$-sectorial ones by a simple coordinate change.

The main goal of the present section is to prove the following result.

\begin{prop}\label{prop:deformation} Let $\omega = d\lambda (= d\lambda_\kappa)$ with $\lambda, \, \lambda_\kappa$ and consider the family
$$
\lambda_t = \lambda - t\, d ((1-\kappa)f)).
$$
Then there exists an isotopy of diffeomorphisms
$
\phi_t: M\to M
$
such that $\phi_0 = id$, $\phi_1 = \phi_\kappa$ and
$$
\phi_t^*(\lambda_t) = \lambda
$$
for $0 \leq t \leq 1$,  and
$$
\supp  \phi_\kappa
 \subset \left\{0 \leq R \leq \frac{\varepsilon_0}{2} \right\}
 \bigcup \{F \setminus F_0\} \times \{|I| \leq N_2\}
$$
in $\nbhd(\del M) \cup \nbhd(\del_\infty M)$.
\end{prop}

To help readers grasp what we are going to do
and because we need the details of the proof, we recall the fine details
of the proof of Gray's theorem in Appendix \ref{sec:gray-theorem}.
Then our proof will be a variation of the stability theorem of
Gray-type  \cite{gray} applied to the $s$-dependent family of contactifications
\be\label{eq:alphas}
Q = M \times \R, \quad
\alpha_s = dt - \pi^*\lambda_s, \quad s \in [0,1]
\ee
on the `space-time' $M \times \R$ with $(y,t)$ as its coordinate system, which
arise from the homotopy of Liouville forms $\lambda_t$ given above.

\subsection{A variation of Gray's stability theorem on contactification}
\label{subsec:gray-contactification}

The proof of the following general result will occupy the whole of this subsection.
In our current situation, we have
\be\label{eq:kt}
k_t(x) = - t\, (1-\kappa) f.
\ee
\begin{theorem}\label{thm:gray}
Let $(M,\lambda)$ be a Liouville manifold and $\lambda_t$ be a
family of Liouville forms such that
$d\lambda_t = d \lambda$
with
 $$
 \lambda_t - \lambda = dk_t
 $$
 for some smooth functions $k_t$ satisfying the bound
 $$
\left\|\frac{\del k_t}{\del t}\right\|_{C^1} < C.
$$
Then there exists a, \emph{not necessarily compactly supported}, diffeomorphism
$\phi_t$ such that
$$
\phi_t^*(\lambda_t) = \lambda_0.
$$
Furthermore we have
$
\supp \phi \subset \supp \left(\frac{\del k_t}{\del t}\right)
  $
for all $t \in [0,1]$.
\end{theorem}

\begin{proof} We consider one-parameter family of contactifications on $Q = M \times \R$
with contact forms given by
$$
\alpha_s = dt - \pi^*(\lambda + dk_s)
$$
which are contact by the hypothesis $d\lambda_t = d\lambda$.
They define a family of contact structures on $Q$ given by
$$
\xi_s: = \ker \alpha_s \quad \text{\rm for } \, s \in [0,1].
$$
We write $\lambda_s = \lambda + dk_s$.

Considering the `space-time' $Q=M \times \R$,  we denote the new time parameter by $s$.
We note that the Reeb vector fields $R_{\alpha_s}$ of each $\alpha_s$ is given by
$$
R_{\alpha_s} = \frac{\del}{\del t} \quad \text{\rm for all }\, s \in [0,1].
$$
We lift the $s$-dependent function $k_s$ to the product
$$
\widetilde k_s(x,t) := \pi^*k_s(x)
$$
which we emphasize \emph{does not} depend on $t$-coordinate of
 the `space-time' $Q = M \times \R$.  This will be important when we
 go back to the study of the family $\lambda_t$ of Liouville one-forms
from our application of  Gray's stability theorem in the contactification.

 As in the general proof of the stability theorem (see Appendix \ref{sec:gray-theorem}), we will try to find a one-parameter
 family of contactomorphisms $\psi_s^*\alpha_s = e^{g_s}\alpha_0$.
 In the current context of our interest, we will try to find
 \emph{strict} contactomorphisms for which $g_s \equiv 0$.
Then we can choose a family of $s$-dependent vector fields
$$
\widetilde X_s \in \xi_s
$$
that we highlight satisfies
\be\label{eq:moser1}
d(\widetilde X_s \rfloor \alpha_s) + \widetilde X_s \rfloor d\alpha_s + \frac{\del \alpha_s}{\del s} = h_s \,  \alpha_s
\ee
for $h_s = \frac{\del g_s}{\del s} \circ \psi_s^{-1} \equiv 0$.
Then it follows that
\be\label{eq:hspsis-1}
\frac{\del\widetilde k_{s}}{\del t} = 0
\ee
and
$$
\alpha_s - \alpha_0 = \lambda - \lambda_s = - d\pi^*k_s.
$$
Therefore we have
$$
\frac{\del \alpha_s}{\del s} = - \pi^* d\dot k_s, \quad \dot k_s: = \frac{\del k_s}{\del s}.
$$
Therefore \eqref{eq:moser1} is equivalent to
\be\label{eq:moser-equation}
d(\widetilde X \rfloor \alpha_s) + \widetilde X \rfloor d\alpha_s
- \pi^* d\dot k_s = 0.
\ee

 Then the vector field $\widetilde X_s$ is uniquely  determined by the equation
\be\label{eq:tildeXs}
\widetilde X_s \in \xi_s, \quad \widetilde X_s \rfloor d\alpha_s = - \pi^* d_M
\left(\frac{\del \widetilde k_s}{\del s}\right).
\ee
By the hypothesis $\left\|\frac{\del k_t}{\del t}\right\|_{C^1} < C$, the vector field
$\widetilde X_s$ is globally Lipschitz.
This implies that the flow of $\widetilde X_s$ exists and satisfies
$$
\widetilde \psi_s^*\alpha_s = \alpha
$$
for all $s \in [0,1]$.

Now we write $\widetilde \psi_s(x,t) = (\psi_s(x,t),b_s(x,t))$ and its generating vector field
\be\label{eq:Xs}
\widetilde X_s(x,t) = X_s(x,t) \oplus a_s(x,t)\frac{\del}{\del t}.
\ee
The condition $\widetilde X_s \in \xi_s$ also implies $0 = \alpha_s(\widetilde X_s) = a_s - \lambda_s(X_s)$, i.e.,
\be \label{eq:das}
a_s  = \lambda_s(X_s).
\ee
Then Moser's deformation equation \eqref{eq:moser-equation} is equivalent to
\be
X_s \rfloor (-d\lambda) = d_M \dot k_s
\ee
where we utilize the identity $d\lambda_{s} = d\lambda$ for all $s$.
Hence we obtain
\be\label{eq:suppXs}
\supp X_s \subset \supp d_M \dot k_s
\ee
and $\supp \psi_s \subset \supp d_M \dot k_s$.

Now we rewrite the projection to $M$ of the equation $\widetilde \psi_s^*\alpha_s = \alpha$
into
$$
\widetilde \psi_s^*(dt - \pi^* \lambda_s) = dt - \pi^* \lambda
$$
which is equivalent to
$$
db_s - \psi_s^*\lambda_s = dt - \pi^*\lambda.
$$
From this, we derive
$$
\psi_s^*\lambda_s = \lambda + d_M b_s, \quad \frac{\del b_s}{\del t} \equiv 1.
$$
Noting the initial condition $(\psi_0(x,t), b_0(x,t)) = (x,t)$,
we in particular proved $b_s(x,t) = t$ for all $x$.
Then by setting $s = 1$, we define
$$
\phi_t(x): = \psi_1(x,t),
$$
which then satisfies
$$
\phi_t^*\lambda_t = \lambda.
$$
Obviously each $\phi_t: M \to M$ is invertible for each $t \in [0,1]$
since the diffeomorphism $(x,t) \mapsto \psi_1(x,t)$ maps each $t$-slice to itself.
This finishes the proof.
\end{proof}

Now we apply Theorem \ref{thm:gray} to the function
$$
k_t = - t (1-\kappa) f
$$
for which we have $\dot k_t = (\kappa -1) f$.

Finally we check the completeness of the flow $\phi_t$. For this purpose, we
go back to the defining equation \eqref{eq:tildeXs} of the vector field $\widetilde X_s$
$$
\widetilde X_s \rfloor d\alpha_s = - \pi^*d_M\widetilde k_s
$$
where we have
$$
\widetilde k_x(t,x) = k_s(x) = -s (1-\kappa(x)) f(x).
$$
Then we have
$$
d_M\widetilde k_s = s d((1-\kappa)f).
$$
We recall the shape of the cut-off function
$$
\kappa(y,R,I) = \kappa_0(y) \kappa_1(R) \kappa_2(I)
$$
where  $\kappa_1 \equiv 0$ near
$\del M = \{R =0\}$ and $\kappa_2 \equiv 1$ when $|I| > N_1$.
From this property of the cut-off functions $\kappa$,
and by the property of $f$ laid out in Proposition \ref{prop:gps},
the function $(1-\kappa)f$ is $C^2$-bounded which shows that the function
$k_t$ satisfies the hypothesis of Theorem \ref{thm:gray}.
This in particular implies that the vector fields $X_t$ are uniformly Lipschitz.
Therefore the flow exists for all time until it hits the boundary $\del M$
of our interest.  This now completes the proof of
Proposition \ref{prop:deformation}.

\begin{remark}\label{rem:noncompact-support}
We compute $\widetilde X_s $ from the result $b_s(t,x) = t$ and hence
$$
\widetilde X_s(t,x) = (X_s(x,t), 0)
$$
and $X_s(x,t) = s X_{d((1-\kappa) f)}(x)$. Recall that $\kappa = \kappa_0(y) \kappa_1(R)\kappa_2(I)$
and  $f$ is supported on $F_0 \times \C_{\text{\rm Re}\geq 0}$
and
$$
d(1-\kappa) = - d\kappa = -\left((d\kappa_0) \kappa_1 \kappa_2
+ \kappa_0 (d\kappa_1) \kappa_2 + \kappa_0 \kappa_1 (d\kappa_2)
\right)
$$
is supported in $F_0 \times \{(R,I) \mid |I| \geq N_0, \, 0 \leq
R \leq \mathsf{ht}(\varphi) \}$. Therefore
$d((1-\kappa)f) = -f\, d \kappa + (1-\kappa) df$ is supported  in the union
\beastar
&{}& F_0 \times \{(R,I) \mid |I| \geq N_0, \, 0 \leq R \leq \mathsf{ht}(\varphi) \} \\
& \bigcup &  F_0 \times \{(R,I) \mid  R \geq \mathsf{ht}(\varphi), \, |I| \leq N_1\}.
\eeastar
(See Definition \ref{defn:kappa}.)
The former set may \emph{not} be compact and so the above constructed Liouville isotopy
may not be compactly supported.
\end{remark}

\subsection{Definition of $\lambda$-sectorial almost complex structures}
\label{subsec:lambda-sectorial-J}

Suppose that such a diffeomorphism $\phi_\kappa$ constructed in Proposition \ref{prop:deformation}
is given. The following is an immediate corollary of this proposition.

\begin{cor} Define a function $\mathfrak{s}_{\varphi,\kappa}$ to be the composition
\be\label{eq:fraksphikappa}
\mathfrak{s}_{\varphi,\kappa} = \mathfrak{s}_{\varphi} \circ \phi_\kappa.
\ee
Then $\mathfrak{s}_{\varphi,\kappa}$ is an exhaustion function of the neighborhood
$$
\nbhd_{2\varepsilon_0}^Z(\del M) \cap \nbhd(\del M).
$$
We call it a \emph{$\lambda_\kappa$-wiggled end-profile function}.
\end{cor}

By applying $\phi_\kappa$ to the $J$-duality equation \ref{item. ds is dual to lambda_kappa},
we get
$$
- d({\mathfrak s}_{\varphi,\kappa} \circ \phi_\kappa) \circ \phi_\kappa^* J = \lambda.
$$
Motivated by this equation, we take a new end-profile function
\be\label{eq:psi}
\psi: = {\mathfrak s}_{\varphi,\kappa}.
\ee
Then we arrive at the final definition of sectorial almost complex structures.
%

\begin{defn}[$\lambda$-sectorial almost complex structures]\label{defn:J-for-corners}
Let $(M,\lambda)$ be a Liouville sector with boundary and corners.
We fix
\begin{itemize}
\item a smoothing profile, i.e., a splitting data and the end-profile function ${\frak s}_\varphi$,
\item the deformation function $\kappa$ adapted to the given smoothing profile as before.
\end{itemize}
An $\omega$-tame almost complex structure $J$
on a Liouville sector is said to be \emph{$\lambda$-sectorial} (with respect to the given smoothing profile and the choice of $\kappa$)
if $J$ satisfies
    $$
   -d\psi \circ J = \lambda, \quad \psi = {\frak s}_{\varphi,\kappa}
    $$
on a neighborhood $\nbhd^Z(\del \Mliou)\cup \del_\infty M)$.
We denote by $\cJ_{(\Mliou,\lambda)}^{\text{\rm sect}}$ the set of $\lambda$-sectorial almost complex structures.
\end{defn}
In the terminology of Definition \ref{defn:liouville-pseudoconvex-intro}, $J$ is Liouville-pseudoconvex
with respect to ${\frak s}_{\varphi,\kappa}$.

The following abundance and contractibility theorem is obvious.
Recall we have fixed a smoothing profile and in particular equip an
end profile function ${\mathfrak s}_{\varphi,\kappa}$ with $(M,\lambda)$
on which the definition of $\lambda$-sectorial almost complex structures depend.
\begin{theorem}\label{thm:J-abundance}
Let $(\Mliou,\lambda)$ be a Liouville sector.
The set of sectorial almost complex structures is
a nonempty and contractible infinite dimensional manifold.
\end{theorem}
\begin{proof}
Existence of sectorial $J$ from that of $\kappa$-sectorial almost complex structure since whenever
$J$ is $\kappa$-sectorial $\phi_\kappa^*J$ is $\lambda$-sectorial.
Certainly the choice of $\omega_F$-tame $J_F$ on $F$  is a contractible choice. Then the choice to be made for
the conditions given in Definition \ref{defn:J-for-corners}
 on the relevant partial neighborhoods of
$$
\nbhd(\del_\infty \Mliou) \cup \nbhd^Z(\del \Mliou)
$$
is unique on the germs of neighborhoods of $\del_\infty M$ and $\del M$. Finally
extending the one defined on the neighborhood to everywhere to $\Mliou$
is also contractible by Gromov's lemma  on the $\omega$-tame almost complex structures.
\end{proof}

\section{Sectorial Floer packages}
\label{sec:sectorial-package}

\subsection{Sectorial Hamiltonians and nonnegative isotopies}
\label{subsec:sectorial-Hamiltonians}

Our effort of constructing pseudo-convex pairs $({\frak s}_{\varphi,\kappa},J)$,
a wiggled end-profile function ${\frak s}_{\varphi,\kappa}$ and its adapted
sectorial almost complex structures $J$ pays off which makes simple the definitions of
Hamiltonians and negative Hamiltonian isotopies that are amenable to the maximum principle.
\begin{remark}
\cite[2.0.1]{gps} regards that it would be ``an important technical advance''
to enlarge the class of Hamiltonians defined near
$\del \Mliou$ for which similar confinement results for holomorphic curves can be proven.
\end{remark}

Let $(M,\lambda)$ be a Liouville sector equipped with a smoothing profile
and a deformation $\lambda_\kappa$ be given. We fix the associated wiggled end-profile function ${\mathfrak s}_{\varphi,\kappa}$.

\begin{defn}[Sectorial Hamiltonians] \label{defn:sectorial-H}
Fix a Liouville sector $M$ and a smoothing profile (Condition \ref{cond:smoothing-profile}) for $M$.
Let ${\mathfrak s}_{\varphi,\kappa}$ be the associated end-profile function.
We call a Hamiltonian $H: M \to \RR$ a \emph{sectorial} with respect to the smoothing profile if
$$
H = \rho({\mathfrak s}_{\varphi,\kappa})
$$
on a neighborhood $\nbhd(\del_\infty M \cup \del M)$ for some function $\rho: \RR \to \RR$.
\end{defn}

\begin{example}
If $H = c \, {\mathfrak s}_{\varphi,\kappa}$ for some constant $c$
on a neighborhood $\nbhd(\del_\infty M \cup \del M)$, then $H$ is sectorial.
\end{example}

It follows from Lemma \ref{lem:brane-and-endprofile-function} that
we can take the neighborhood $\nbhd(\del_\infty M \cup \del M)$ to be
$$
\nbhd_{s \geq N}(\del_\infty M) \cup \nbhd_{\epsilon}(\del M).
$$
for a sufficiently large $N$ and small $\epsilon > 0$.

\begin{prop}
The space of sectorial (with respect to a fixed choice of $\bf s_{\varphi,\kappa}$) Hamiltonians is nonempty and
is contractible.
\end{prop}
\begin{proof}
This is obvious. For example, the space of $H$ is convex.
\end{proof}

On noncompact manifolds like Liouville sectors,
the Floer continuation map (and hence Hamiltonian invariance) is available only in a preferred, `nonnegative' direction.

\begin{notation}
Let $H = H(s,t,x)$ be a one-parameter family of (time-dependent)
Hamiltonians $H^s(t,x): = H(s,t,x)$ and consider the $(s,t)$-family of
Hamiltonian diffeomorphisms $(s,t) \mapsto \phi_{H^s}^t$.
\end{notation}

The confinement theorem, Theorem \ref{thm:nonautonomous-confinement},
leads us to the following simple notion of \emph{sectorial nonnegative isotopies}
which is nothing but the natural definition of nonnegative isotopies in the
context of sectorial Hamiltonians.

\begin{defn}[Sectorial nonnegative isotopy]
\label{defn:nonnegative-isotopy} We say a one-parameter family of Hamiltonian isotopies $\psi^s = \phi_{H^s}^1$
is a \emph{sectorial nonnegative isotopy} if it satisfies
\begin{enumerate}
\item Each $H^s$ is sectorial. In particular,
on a neighborhood $\nbhd(\del_\infty M \cup \del M)$, $H^s = \rho^s({\mathfrak s}_{\varphi,\kappa})$  for some $s$-family $\rho^s: \RR \to \RR$.
\item The family $\rho^s$ satisfies
\eqn\label{eq:rhos'-nonnegative}
\frac{\del \rho^s}{\del s} \geq 0.
\eqnd
\end{enumerate}
\end{defn}

\begin{remark}
We would like to emphasize that we do not need to
put any other conditions on the behavior of Hamiltonians near $\del_\infty M \cup \del M$
for the confinement results. Compare this with the kind of Hamiltonians such as dissipative Hamiltonians
used in \cite{gps}.
\end{remark}

We note that the condition \eqref{eq:rhos'-nonnegative} is equivalent to the standard definition of
nonnegative Hamiltonians
$$
\frac{\del H^s}{\del s} \geq 0
$$
applied to sectorial Hamiltonians $H^s$.

\begin{example} Let $H_1, \, H_2$ be sectorial Hamiltonians. Then the linear interpolation
 $$
 s \mapsto H^s = (1-s)H_1 + sH_2
 $$
induces a sectorial nonnegative isotopy if $H_2 \geq H_1$ on
a neighborhood $\nbhd(\del_\infty M \cup \del M)$ of $\del_\infty M \cup \del M$.
\end{example}

We relate this definition in the more common notion of nonnegative isotopies in the literature.

\begin{example} \begin{enumerate}
\item In the literature on Floer theory for Liouville manifolds, the asymptotic behavior of
the (autonomous) Hamiltonians is required to be fixed; for example, by requiring $H$ to be a constant multiple of
the symplectization radial coordinate $r$. For any isotopy induced by an $s$-independent family
of sectorial Hamiltonians $H^s$ satisfying $H^s \equiv f(r)$, the isotopy induced by such a family
is trivially a nonnegative isotopy.
\item In \cite[Section 4.8]{gps}, some homotopy coherent diagrams are used to define the homotopy limits of
symplectic cohomology. In particular the authors of \cite{gps} used an explicit one-simplex $\{H^s\}$ of
Hamiltonians with certain requirements such as Equation (4.69) and others in the same section of \cite{gps}.
When we consider such a homotopy diagram in our sectorial framework, the only thing we need to care about
is a suitable continuation of background geometry of sectorial hypersurfaces and boundaries in the way that
is amenable to our sectorial framework. For this purpose, one may need to also consider
the higher simplex version of the nonnegativity requirement in the bundle setting.
(See \cite{oh-tanaka:liouville-bundles,oh-tanaka:actions} for such a practice.)
\end{enumerate}
\end{example}

\subsection{Sectorial Lagrangians}
\label{subsec:sectorial-Lagrangians}

We now identify the class of Lagrangians that define objects for a wrapped Fukaya category of Liouville sectors with corners.
We omit discussion of other brane data such as spin structures since they are not different
from the usual description in the literature.

\begin{defn}[Sectorial Lagrangians]\label{defn:sectorial-Lagrangians}
Let $(M,\lambda)$ be a Liouville sector with corners. Denote $\omega = d\lambda$.
We say an exact, properly embedded Lagrangian submanifold $L$ of $(M,\omega)$ is \emph{sectorial} if
\begin{itemize}
\item $L \subset M \setminus \del M$.
\item There exists a neighborhood of
$\del_\infty M$ whose intersection with $L$ is $Z$-invariant.
\end{itemize}
\end{defn}

\begin{remark}
The definition of sectorial Lagrangians
depends only on the Liouville structure $(M,\lambda)$
while the sectorial pair $(J,H)$
depends on additional data like $\psi= {\frak s}_{\varphi,\kappa}$
given in \eqref{eq:fraksphikappa}.
\end{remark}

\begin{lemma}\label{lem:brane-and-endprofile-function}
Assume that $\nbhd(\del_\infty M)$ is a $Z$-invariant neighborhood.
Fix ${\frak s}_{\varphi,\kappa}$ (and in particular, $\varphi$ and $\kappa$).
Then there exist a sufficiently
small $\epsilon > 0$ and a sufficiently large $N_0 > 0$ so that
$$
{\mathfrak s}_{\varphi,\kappa}^{-1}([N_0, \infty)) \cap \nbhd_\epsilon(\del M) \subset
\nbhd(\del_\infty M)
$$
\end{lemma}

\begin{remark}
This exhaustion property of end-profile functions of ${\mathfrak s}_{\varphi,\kappa}$ makes
the sectorial Lagrangian, which is $Z$-invariant on a (unspecified) neighborhood $\nbhd(\del_\infty M)$,
amenable to an application of the strong maximum principle for the function ${\mathfrak s}_{\varphi,\kappa}$
on the region ${\mathfrak s}_{\varphi,\kappa}^{-1}([N_0, \infty))$ for a sufficiently large $N_0 > 0$.
This is what enables us to make the definition of sectorial Lagrangian
independent of the choice of end-profile functions by unchanging the standard definition of
$Z$-invariant Lagrangians in the literature in our sectorial framework.
See \cite{oh:gradient-sectorial} for a different
class of Lagrangians that is amenable to the pseudoconvex pair $({\frak s}_\varphi,J)$ for
sectorial almost complex structures, not $\lambda$-sectorial ones as in the present paper.
\end{remark}

\clearpage
\part{Confinement theorems}
\label{part:confinement}

In this part, with our sectorial Floer package, we establish the fundamental confinement result
for various Floer-type equations: There are two parts of them, one the \emph{vertical}
uniform $C^0$ estimates and the other the \emph{horizontal} $C^0$ estimates. The former prevents
the relevant pseudoholomorphic curves from escaping to infinity and the latter from
coming too close to the boundary $\del M$ from the interior of $M$. Because of the different
natures of the geometries of $\del_\infty M$ and $\del M$ as mentioned before, the existing literature
(e.g., \cite{gps}) handle the two differently even for the case without background Hamiltonian $H$ turned on.
When the Hamiltonian is turned on, they also restrict the class of Hamiltonians to those with
somewhat complicated notion of dissipative Hamiltonians \cite{groman} to achieve the relevant confinement
theorems for the Hamiltonian-perturbed Floer trajectory equations which is one of the crucial
element in the functorial studies of wrapped Fukaya categories and of symplectic cohomologies on
Liouville sectors.

On the other hand, utilizing our sectorial framework,
we will prove all these confinement results by applying the maximum principle
to the $\kappa$-wiggled end-profile function ${\mathfrak s}_{\varphi,\kappa}$, and
the strong maximum principle for the boundary value problem of $Z$-invariant-at-infinity Lagrangians. We establish the
confinement results for the following types of Floer equations on the Liouville sectors with corners:
\begin{itemize}
\item the structure maps of the associated $A_\infty$ category,
\item the Hamiltonian-perturbed Floer equation,
\item the Floer continuity equation under the nonnegative isotopies of Hamiltonians,
\item the closed-open (the open-closed) maps from the symplectic Floer cohomology to
the wrapped Fukaya category.
\end{itemize}

\section{The curves we care about}
We lay out the relevant Floer equations for which we establish maximum principles.

\subsection{Structure maps of unwrapped Fukaya category}
\label{section. unwrapped disks}

In the study of (unwrapped) Fukaya category, we consider a disc $D^2$ with a finite number of boundary marked points $z_i \in \del D^2$ equipped with strip-like coordinates $(\tau,t)$ (or on the sphere $S^2$) with a finite number of marked points.
We denote by $\overline{z_iz_{i+1}}$ the arc-segment between $z_i$ and $z_{i+1}$, and $\tau = \infty_i$
the infinity in the strip-like coordinates at $z_i$.

Let $L_0,  \ldots, L_k$ be a $(k+1)$-tuple of Lagrangian submanifolds which are
$Z$-invariant near infinity. We
denote
$$
\Sigma = D^2 \setminus \{z_0, \ldots, z_k\}
$$
and equip $\Sigma$ with strip-like coordinates $(\tau,t)$
with $\pm \tau \in [0,\infty)$ and $t \in [0,1]$ near each $z_i$.

Then for a given collection of intersection points $p_i \in L_i\cap L_{i+1}$ for $i = 0, \ldots, k$,
we wish to study maps $u: \Sigma \to \Mliou$ satisfying the Cauchy-Riemann equation
\eqn\label{eq:unwrapped-structure-maps}
\begin{cases}
\overline \partial_J u = 0\\
u(\overline{z_iz_{i+1}}) \subset L_i \quad & i = 0, \ldots k\\
u(\infty_i,t) = p_i, \quad & i = 0, \ldots k.
\end{cases}
\eqnd

\subsection{Symplectic cohomology and its continuation map}
\label{subsec:continuation}
\label{section. SH}

Let $H$ be an (autonomous) sectorial Hamiltonian and $\ell_0, \, \ell_1$ be a pair of
periodic orbits of $H$.
We consider a map $u: \RR \times S^1 \to M$ satisfying
\eqn\label{eq:CRJH}
J(du - X_H(u) \otimes dt) = (du - X_H(u) \otimes dt) \circ j
\eqnd
which is equivalent to saying that the $TM$-valued one-form $du-X_H(u) \otimes dt$ is $(j,J)$-holomorphic;
i.e., that $u$ satisfies
\eqn\label{eq:du-XH01}
(du - X_H(u) \otimes dt)^{(0,1)}_J = 0.
\eqnd
By evaluating this equation against $\frac{\del}{\del \tau}$, we get the standard Floer equation
\eqn\label{eq:duJHu}
\dudtau + J\left(\dudt - X_H(u) \right) = 0
\eqnd
back. Then we put the asymptotic condition
\eqn\label{eq:asymptotic condition}
u(- \infty,t) = \ell_0(t), \quad u(\infty,t) = \ell_1(t).
\eqnd

For the continuation map, one commonly considers one-parameter family of sectorial Hamiltonians $\{H^s\}_{s \in [0,1]}$
and an elongation function $\chi: \RR \to [0,1]$ with $\chi(\tau) \equiv 0$ for $\tau \leq 1$ and $\chi(\tau) \equiv 0$
for $\tau \geq 1$ and consider the non-autonomous analog to \eqref{eq:CRJH},
\eqn\label{eq:nonoauto-duJHu}
\dudtau + J\left(\dudt - X_{H^\chi}(u) \right) = 0
\eqnd
where $H^\chi: \RR \times M \to \RR$ is the $\chi$-elongated one-parameter family of $\{H^s\}$ defined by
$$
H^\chi(\tau,x) := H^{\chi(\tau)}(x).
$$

\subsection{Closed-open maps and open-closed maps}
\label{subsec:closed-open-map}

These are the maps which connect the symplectic cohomology and the Hochschild (co)homology of
the Fukaya category of the underlying symplectic manifold in general.

The following is a variation of the similar moduli spaces that have appeared in
\cite[Definition 18.21]{fooo-memoir} (for $\ell = 1$) in the closed context (See also \cite{albers}.)

\begin{defn}
Consider a collection $\vec H = \{H_j\}_{j=1, \ldots, \ell}$ and let $\text{\rm Per}(H_j)$ be the set of periodic orbits.
We fix for all $j$ $\ell_j \in \text{\rm Per}(H_j)$.
We consider the moduli space of all pairs
$$
\left(u, (D^2; z_1^{+},\dots,z_{\ell}^+;z_0,\dots,z_k)\right)
$$
satisfying the following:
\begin{enumerate}
\item $z_1^{+},\dots,z_{\ell}^+$ are points in $\Int D^2$
which are mutually distinct.
\item
$z_0,\dots,z_k$ are  points on the boundary $\del D^2$
which are ordered counterclockwise on $S^1 = \del D^2$
with respect to the boundary orientation coming from
$D^2$.

\item Write $\Sigma = D^2 \setminus (\{z_j^+\} \cup \{z_i\})$ and equip $\Sigma$ with
cylindrical coordinates
$$
\epsilon_j^+: (- \infty, 0] \times S^1 \to \Sigma
$$
near $z_j^+$  and strip-like coordinates
$$
\epsilon_i: (-\infty, 0] \times [0,1]
$$
for $i = 1, \ldots, k$ and
$$
\varepsilon_0: [0, \infty) \times [0,1] \to \Sigma
$$
near $z_0$.
\item The map
$u : \Sigma \to M$ is a smooth map
such that $u(\overline{z_iz_{i+1}}) \subset L_i$.
\item
The map $u$ satisfies the equation
\eqn\label{eq:HJCR18sec}
(du - X_K(u) \otimes \gamma)^{(0,1)}_J = 0
\eqnd
for some subclosed one-form $\gamma$ on $\Sigma$, where $K = K(z,x):\Sigma \times M \to \RR$ is the domain-dependent
Hamiltonian such that
\eqn\label{eq:H-at-ends}
H(\epsilon_j(\tau,t),x) = H_j(t,x)
\eqnd
on each strip-like end $\epsilon_j: [0,\infty) \times [0,1] \to \dot \Sigma$ of $\dot \Sigma$.
\item We put the asymptotic condition
\eqn\label{eq:asymptotic condition}
u(z_i) = p_i \in L_i \cap L_{i+1}, \quad i = 0, \ldots k
\eqnd
for $z_i \in \del D^2$.
\end{enumerate}
\end{defn}

The open-closed map is defined similarly by considering the above maps
in the reverse direction of $\tau$.

\subsection{Quantitative measurements for the $C^0$-estimates}

In the study of each of the above Floer moduli spaces, one needs to
establish $C^0$-estimates  as the first step towards the study of analytic estimates
and compactification of the moduli spaces.

We now provide the quantitative measurements that enter in our
$C^0$-estimates uniformly for all the above moduli spaces.

\begin{defn}\label{defn:peak} Write $\vec L = \{L_0,\ldots, L_k\}$ and $\vec H_\epsilon = (H_0,\ldots, H_k)$.
\begin{enumerate}
\item Define
$$
\text{\rm Per}(\vec H_\epsilon):= \bigcup_{i} \text{\rm Per}(H_i)
$$
and
\eqn\label{eq:peak}
\mathfrak{ht}_{\text{\rm Per}}(\vec L,\vec H_\epsilon): = \max_{i} \left(\max\{s(\ell(t))\mid \ell \in \mathfrak{X}(L_i,L_{i+1};H_i)\}\right).
\eqnd
We call $\mathfrak{ht}_{\text{\rm Per}}(\vec L,\vec H)$ the \emph{$\vec H$-orbit height} of $\vec L$.
\item Let a symplectization radial function $s$ be given. Define the \emph{$Z$-invariance support} of $L$, denoted by $\supp^{Z}(L)$,
$$
\supp^{Z}(L) : =
\overline{M \setminus \bigcup_{N >0} \{x \in L \cap \nbhd_{\{s \geq N\}}(\del_\infty M) \mid  Z(x)\, \text{ \rm is tangent to } L\}}
$$
and its \emph{$Z$-invariance level} of $L$ by
$$
\text{\rm level}^Z(L) : = \max_{x \in \supp^{Z}(L)} s(x).
$$
Then we call
$$
\text{\rm level}^Z(\vec L): = \max_{i=0}^k \text{\rm level}^Z(L_i).
$$
the \emph{$Z$-invariance level} of the collection $\vec L$.
\item For each Lagrangian $L$ with $L \cap \del M = \emptyset$, we denote the distance
 between $L$ and $\del M$ by
$$
d(L;\del M) = \inf \{d(x,y) \mid x \in L, \, y \in \del M\}
$$
and $d(\vec L;\del M) = \min \{d(L;\del M)\mid L \in \vec L\}$ (with respect to
any given metric $g$ of bounded geometry, e.g., an $\omega$-tame metric $g$).
\item The quantities of the width $\mathsf{wd}(\varphi)$ and its height $\mathsf{ht}(\varphi)$
of the function $\varphi \in \mathfrak{Conv}_{\mathfrak{sm}}(\RR^k_+)$, which is fixed for each
sectorial corner,  will enter the measurements via the associated end-profile function $s_{k+1,\varphi}$.
\end{enumerate}
\end{defn}

\begin{remark}
By definition, if $L$ (resp. $\vec L$) is $Z$-invariant at infinity, $\text{\rm level}^Z(L) < \infty$
(resp. $\text{\rm level}^Z(\vec L) < \infty$).
\end{remark}

\section{Confinement for the structure maps of unwrapped Fukaya categories}
\label{sec:structure-maps}

We work in the framework of Section~\ref{section. unwrapped disks}.
\begin{defn}[Intersection height]\label{defn:peak} Let a symplectization radial function $s$ be given
Let $\vec L = \{L_0,\ldots, L_k\}$ and define
Define
\eqn\label{eq:peak}
\mathfrak{ht}(\vec L): = \max_{i} \left(\max\{s(p_i)\mid p \in L_i \cap L_{i+1} \}\right).
\eqnd
We call $\mathfrak{ht}(\vec L)$ the \emph{intersection height} of $\vec L$.
\end{defn}

We have
$$
\mathfrak{ht}(\vec L) \leq \text{\rm level}^Z(\vec L)
$$
when $L_i$'s are pairwise disjoint above the $Z$-invariance level.

With this definition, we establish the following $C^0$-estimate from the $J$-convexity of
${\mathfrak s}_{\varphi,\kappa}$ and the  maximum principle and the strong maximum principle.

\begin{remark}
For the usual $C^0$-estimates via the maximum principle in this bordered case (see, e.g., \cite{EHS})  we recall the importance of the role of \emph{$J$ being of contact-type},
(or more precisely $J$ satisfying the equation $-ds \circ J = \lambda$) and the $Z$-invariance of the Lagrangians at infinity.
\end{remark}

\begin{theorem}\label{thm:unwrapped}
Let $J$ be a sectorial almost complex structure of the Liouville sectors with corners $M$. Let $u$ be a solution to~\eqref{eq:unwrapped-structure-maps} as in~\ref{section. unwrapped disks}.
Then the following hold:
\begin{enumerate}
\item {(Horizontal $C^0$ estimates)} There exists some $\epsilon = \epsilon(\vec L)> 0$ independent of $u$'s such that
$$
\Image u \subset \Mliou \setminus \nbhd_\epsilon(\del \Mliou).
$$
\item {(Vertical $C^0$ estimates)} There exists $N_0= N_0(\vec L) > 0$ independent of $u$ such that
$$
\Image u \subset \Mliou \setminus \nbhd_{s \leq N_0}(\del_\infty \Mliou).
$$
\end{enumerate}
\end{theorem}
\begin{proof}
First since all Lagrangian submanifolds $L_i$'s are $Z$-linear at infinity, all $L_i$'s are contained in the interior
$\{R > \epsilon_1\}$ for $\epsilon_1 = d(\vec L,\del M)$.

Then by definition of sectorial almost complex structure $J$, $J$ is associated to a splitting data and
the wiggled end-profile function ${\mathfrak s}_{\varphi,\kappa}$.
Since a neighborhood of $\del_\infty M \cup \del M$ is exhausted by the family of hypersurfaces
$$
({\mathfrak s}_{\varphi,\kappa})^{-1}(r)
$$
for $r \geq 0$ (see Lemma \ref{lem:exhaustion}), it is enough to prove
\eqn\label{eq:confimennt}
\Image u \subset ({\mathfrak s}_{\varphi,\kappa})^{-1}((-\infty,r])
\eqnd
for some $r > 0$, which will establish both vertical and horizontal bounds
simultaneously: Note that the inequality
$$
{\mathfrak s}_{\varphi,\kappa} \leq r_0
$$
implies both $s \leq N + r_0$ and
$$
\max_{i=1, \ldots, n} |\log R_i| \leq r_0.
$$
We may take the vertical upper bound to be
$$
N_0(\vec L): =  \max\{\mathfrak{ht}(\vec L), \, \text{\rm level}^Z(L), N + r_0\}
$$
and the horizontal lower bound to be
$$
\epsilon(\vec L): = \min_{i=1, \ldots n} R_i \geq e^{-r_0}.
$$

We first recall that $du$ is $J$-holomorphic and
$
- d(d {\mathfrak s}_{\varphi,\kappa}  \circ J) = d\lambda
$
from Definition \ref{defn:J-for-corners}.
By definition of $[d\lambda_M]_{(J;U_{\epsilon_2})}$ (Definition~\ref{defn:pinching-bound})
and since $d\left({\mathfrak s}_{\varphi,\kappa} \circ u\right) \circ j = d{\mathfrak s}_{\varphi,\kappa}  \circ J \circ du$, we have
\eqn\label{eq:dd1-R>}
-d\left(d\left({\mathfrak s}_{\varphi,\kappa} \circ u\right) \circ j\right)\geq \frac{C}2  |du|_J^2\, \omega_\Sigma
\eqnd
with
$$
C = [-d{\mathfrak s}_{\varphi,\kappa} \circ  J]_{J;U_{\epsilon_2}} = [d\lambda]_{(J;U_{\epsilon_2})} > 0.
$$
In particular, the function ${\mathfrak s}_{\varphi,\kappa} \circ u$ is a subharmonic function and cannot carry an interior
maximum on $\RR \times [0,1]$ by the maximum principle.

Next we will show by the strong maximum principle that $u$ cannot have a boundary maximum
in a neighborhood of
$
\del_\infty M \cup \del M
$
either. This will then enable us to obtain a
$C^0$ bound $r_0 > 0$ such that
$$
\Image u \subset \{{\mathfrak s}_{\varphi,\kappa} \leq r_0\}
$$
for any finite energy solution $u$ with fixed asymptotics~\eqref{eq:asymptotic condition}.

Now suppose to the contrary that ${\mathfrak s}_{\varphi,\kappa} \circ u$ has a boundary local maximum
point $z' \in \del D^2\setminus \{z_0,\ldots, z_k\}$. By the strong maximum principle, we must have
\eqn\label{eq:lambda(dudtheta)}
0 < \frac{\del}{\del \nu}({\mathfrak s}_{\varphi,\kappa} \circ u)
= d{\mathfrak s}_{\varphi,\kappa}\left(\frac{\del u}{\del \nu}\right)
\eqnd
for the outward unit normal $\frac{\del}{\del \nu}$ of $\del \Sigma$,
unless ${\mathfrak s}_{\varphi,\kappa}$ is a constant function in which case there is nothing to prove.
Let $(r,\theta)$ be an isothermal coordinate of a neighborhood of $z_0 \in \del \Sigma$ in $(\Sigma,j)$
adapted to $\del \Sigma$, i.e., such that $\frac{\del}{\del \theta}$ is tangent to $\del \Sigma$ and
$|dz|^2 = (dr)^2 + (d\theta)^2$ for the complex coordinate $z = r+ i\theta$ and
\eqn\label{eq:normal-derivative}
\frac{\del}{\del \nu} = \frac{\del}{\del r}
\eqnd
along the boundary of $\Sigma$.
Since $u$ is $J$-holomorphic, we also have
$$
\frac{\del u}{\del r} + J \frac{\del u}{\del \theta} = 0.
$$
Therefore combining this with the duality requirement
$$
-d{\mathfrak s}_{\varphi,\kappa} \circ J = \lambda,
$$
we derive
$$
d{\mathfrak s}_{\varphi,\kappa}\left(\frac{\del u}{\del \nu}\right)
= d{\mathfrak s}_{\varphi,\kappa} \left(-J \frac{\del u}{\del \theta}\right)
= \lambda \left(\frac{\del u}{\del \theta}\right).
$$
By the $Z$-invariance of $L$ and the boundary condition $u(\del \Sigma) \subset L$,
both $Z(u(z_0))$ and $\frac{\del u}{\del \theta}(z_0)$ are contained in
$T_{u(z_0)}L$, which is a $d\lambda$-Lagrangian subspace. Therefore we have
$$
0 = d\lambda\left(Z,\frac{\del u}{\del \theta}\right) = \lambda\left(\frac{\del u}{\del \theta}\right)
$$
where the second equality follows from the definition of Liouville vector field $Z$ of $\lambda$.

This is a contradiction to \eqref{eq:lambda(dudtheta)} and hence
the function ${\frak s}_\varphi \circ u$ cannot have a boundary maximum either.
This finishes the proof.
\end{proof}


\begin{remark}
The same confinement result holds for
the Floer's continuation equation for the moving Lagrangian boundary condition of
nonnegative isotopies, which we refer readers to the most updated versions of
\cite{oh-tanaka:liouville-bundles} (for the unwrapped case),
and \cite{oh-tanaka:actions} (for the wrapped case) for the details.
\end{remark}

\section{Confinement theorems for symplectic cohomology}

We work in the setting of Section~\ref{section. SH}.
We consider two cases separately, one for the structure map for the symplectic cohomology (i.e., for the autonomous Hamiltonian)
and the other for the case of Floer continuation map (i.e., for the nonautonomous Hamiltonian).
%

In this section, we consider a (time-dependent) sectorial Hamiltonian
$H = H(t,x)$.


\begin{theorem}\label{thm:autonomous-confinement}
Let $(M,\lambda)$ be a Liouville sector.
Let $(J,H)$ be sectorial. Suppose $u$
satisfies the equation \eqref{eq:du-XH01} and $\Image u \cap \del M = \emptyset$.
Then the following hold:
\begin{enumerate}
\item {(Horizontal $C^0$ estimates)} There exists some $\epsilon = \epsilon(\vec L)> 0$ independent of $u$'s such that
$$
\Image u \subset \Mliou \setminus \nbhd_\epsilon(\del \Mliou).
$$
\item {(Vertical $C^0$ estimates)} There exists $N_0= N_0(\vec L, H) > 0$
 independent of $u$'s such that
$$
\Image u \subset \Mliou \setminus \nbhd_{s \leq N_0}(\del_\infty \Mliou).
$$
\end{enumerate}
\end{theorem}

The rest of this subsection will be occupied by the proof of this theorem.
Again it is enough to prove
\eqn\label{eq:confimennt}
\Image u \subset ({\mathfrak s}_{\varphi,\kappa})^{-1}((-\infty,r])
\eqnd
for some $r > 0$.

\begin{prop}\label{prop:energy-identity} Let $J$ be a $\lambda$-sectorial almost complex structure,
and $H = H(t,x)$ be any sectorial Hamiltonian.  Then for any solution $u$ of \eqref{eq:CRJH},
we have
\begin{eqnarray}\label{eq:Laplacian-u}
\Delta({\mathfrak s}_{\varphi,\kappa}\circ u) = \left|\frac{\del u}{\del \tau} - X_H(u)\right|^2_J -
\rho'({\mathfrak s}_{\varphi,\kappa} \circ u) \frac{\del}{\del \tau}({\mathfrak s}_{\varphi,\kappa} \circ u).
\end{eqnarray}
on a neighborhood $\nbhd(\del_\infty M \cup \del M)$.
\end{prop}
\begin{proof} We restrict ourselves to a neighborhood contained in
$$
\nbhd_{I\geq N}(\del_\infty M) \cap  \nbhd_{2\varepsilon_0}(\del M)
$$
of the ceiling of a sectorial corner $C_\delta$. Then
the function $f: F \times \C_{\text{\rm Re}\geq 0} \to \RR$ satisfies the properties
given in Proposition \ref{prop:gps}. We also focus on the neighborhood of a fixed sectorial corner
where the end-profile function ${\mathfrak s}_{\varphi,\kappa}$ is defined.
%
%

By taking the differential of $-u^*(d{\frak s}_{\varphi,\kappa} \circ J) = u^*\lambda$ for the
 $\lambda$-sectorial $J$, we obtain
$$
- u^*d(d{\mathfrak s}_{\varphi,\kappa} \circ J)) = u^*d\lambda.
$$
Using this and the equation \eqref{eq:duJHu}, we now evaluate $u^*d\lambda$ against the
tuple $(\frac{\del}{\del \tau}, \frac{\del}{\del t})$ and derive
\beastar
u^*d\lambda & = & d\lambda\left(\frac{\del u}{\del \tau}, \frac{\del u}{\del t}\right) \\
& =& d\lambda\left(\frac{\del u}{\del \tau}, \frac{\del u}{\del t} - X_H(u)\right)
+ d\lambda\left(\frac{\del u}{\del \tau}, X_H(u)\right)\\
& = & \left|\frac{\del u}{\del \tau} - X_H(u)\right|^2_J -
\rho'({\mathfrak s}_{\varphi,\kappa} \circ u) \frac{\del}{\del \tau}({\mathfrak s}_{\varphi,\kappa} \circ u).
\eeastar
For the last equality, we have also used $H = \rho({\mathfrak s}_{\varphi,\kappa} \circ u)$ and calculations
\beastar
d\lambda\left(\frac{\del u}{\del \tau}, X_H(u)\right) & = & - d_MH(u)\left(\frac{\del u}{\del \tau}\right)\\
& = & - \frac{\del}{\del \tau}(H(t,u)) + \left(\frac{\del H_t}{\del \tau}\right)(u) \\
& = &  - \frac{\del}{\del \tau}(H(t,u)) = -  \frac{\del}{\del \tau}\left(\rho({\mathfrak s}_{\varphi,\kappa} \circ u)\right)\\
& = & - \rho'({\mathfrak s}_{\varphi,\kappa} \circ u) \frac{\del}{\del \tau}({\mathfrak s}_{\varphi,\kappa} \circ u).
\eeastar
This finishes the proof.
\end{proof}

\begin{proof}[Proof of Theorem~\ref{thm:autonomous-confinement}.]
An immediate corollary of Proposition~\ref{prop:energy-identity} and \eqref{eq:CRJH} is that we can apply the
maximum (as well as the strong maximum) principle to conclude that
${\mathfrak s}_{\varphi,\kappa}\circ u$ cannot carry any local maximum point.
(See Appendix \ref{sec:maximum-principle} for some more details for how the strong maximum principle applies.)

In particular,  we have a uniform $C^0$ bound $r_0 > 0$ such that
$$
\Image u \subset \{{\mathfrak s}_{\varphi,\kappa} \leq r_0\}
$$
for any finite energy solution $u$.
\end{proof}

\begin{remark}\label{remark. confinement away from del M}
Theorem \ref{thm:autonomous-confinement} implies that the moduli space of $J$-holomorphic curves $v: \Sigma \to M$ is a disjoint union of those curves that touch $\del M$, and those curves that are bounded away from a neighborhood of $\del M$. In particular, if one member in a continuous family of $J$-holomorphic curves is bounded away from $\del M$, so are all other members of that family. In the Lagrangian Floer theory of Liouville sectors, we will always consider $J$-holomorphic curves whose images are contained in $\Int(M) = M \setminus \del M$.
Theorem \ref{thm:autonomous-confinement} ensures Gromov compactness for this class of curves.
\end{remark}

\section{Confinement for the continuation under nonnegative Liouville isotopies}
\label{sec:continuation}

In this section, we consider a one-parameter family of sectorial Hamiltonians $H^s = H^s(t,x)$
and consider Floer continuation maps for the symplectic cohomology.
The Floer continuation equation is given by
\eqn\label{eq:continuation-equation}
\dudtau + J\left(\dudt - X_{H_i^{\chi(\tau)}}(u)\right) = 0
\eqnd
for an elongation function $\chi:\RR \to S^1$ with $\chi(\tau) = 0$ such that
\begin{eqnarray}\label{eq:chi}
\chi(\tau) & = &
\begin{cases} 1, \quad & \text{\rm at }\, \tau = -\infty\\
0 \quad & \text{\rm at }\, \tau = \infty
\end{cases} \nonumber\\
\chi' & \leq & 0
\end{eqnarray}

\begin{theorem}\label{thm:nonautonomous-confinement}
Let $J$  and $H^s = H^s(t,x)$ be a sectorial nonnegative isotopy and consider the
Floer continuation equation \eqref{eq:continuation-equation}
for a map $u: \RR \times S^1 \to \Mliou$
with finite geometric energy and $\Image u \cap \del M = \emptyset$.
Then for each given periodic orbits $z_\pm$ of Hamiltonians $H_\pm = H_\pm(t,x)$ respectively
with $\Image z_\pm \cap \del M = \emptyset$, there exists some $r_0 > 0$ such that
any solution for $u$ with finite geometric energy satisfies
\eqn\label{eq:confirment-for-nonautonomous-H}
\Image u \subset ({\mathfrak s}_{\varphi,\kappa})^{-1}((-\infty,r_0]).
\eqnd
\end{theorem}
\begin{proof} As before, we compute $\Delta({\mathfrak s}_{\varphi,\kappa} \circ u) \, d\tau \wedge dt$
\eqn\label{eq:Laplacian-hu}
\Delta({\mathfrak s}_{\varphi,\kappa} \circ u) \, d\tau \wedge dt
= -d(d({\mathfrak s}_{\varphi,\kappa} \circ u)\circ j) = u^*d\lambda.
\eqnd
Then we have
$$
u^*d\lambda = \left(\left|\dudtau - X_{H^\chi}(u)\right|^2_J -
\frac{\del}{\del \tau}(\rho^\chi({\mathfrak s}_{\varphi,\kappa}\circ u))\right) \, d\tau\wedge dt.
$$
Combining the two, we have derived
$$
\Delta({\mathfrak s}_{\varphi,\kappa} \circ u) = \left|\dudtau - X_{H^\chi}(u)\right|^2_J -
\frac{\del}{\del \tau}(\rho^\chi({\mathfrak s}_{\varphi,\kappa}\circ u))
$$
similarly for \eqref{eq:Laplacian-u}, except $\rho$ replaced by $\rho^{\chi}$ which now
varies depending on the variable $\tau$.

We rewrite
$$
\frac{\del}{\del \tau}(\rho^\chi({\mathfrak s}_{\varphi,\kappa}\circ u)) =
(\rho^{\chi(\tau)})'({\mathfrak s}_{\varphi,\kappa}\circ u)\frac{\del}{\del \tau}({\mathfrak s}_{\varphi,\kappa}\circ u)+
\chi'(\tau) (\rho^{\chi(\tau)})'({\mathfrak s}_{\varphi,\kappa}\circ u).
$$
Since $\chi'\leq 0$ and $(\rho^s)'> 0$, we have derived
$$
- \frac{\del}{\del \tau}(\rho^\chi({\mathfrak s}_{\varphi,\kappa}\circ u)) \geq
- (\rho^{\chi(\tau)})'({\mathfrak s}_{\varphi,\kappa}\circ u)\frac{\del}{\del \tau}({\mathfrak s}_{\varphi,\kappa}\circ u).
$$
In particular, we have derived
\eqn\label{eq:Laplacian(su)}
\Delta({\mathfrak s}_{\varphi,\kappa} \circ u) \geq
- (\rho^{\chi(\tau)})'({\mathfrak s}_{\varphi,\kappa}\circ u)\frac{\del}{\del \tau}({\mathfrak s}_{\varphi,\kappa}\circ u).
\eqnd
Therefore we can apply the (interior) maximum principle. We postpone till
Appendix \ref{sec:maximum-principle} to explain how this differential inequality together with
the $Z$-invariance of the given Lagrangian boundary condition enables us to apply the strong maximum
principle too.

This finishes the proof.
\end{proof}

The above theorem can be rephrased as that the Floer continuation equation is amenable to
the maximum principle for the sectorial nonnegative homotopies $\{H^s\}$.

\begin{remark} It is worthwhile to highlight the pair of conditions in Definition \ref{defn:nonnegative-isotopy}
enables us to establish all the necessary $C^0$-estimates in our sectorial framework
by the maximum principle, using only the standard geometric calculations that have been
used in the Hamiltonian geometry and in the analysis of pseudoholomorphic curves.
All of the proofs of the $C^0$-estimates in the literature involve some type of analytic comparison estimates
between the growth rate of the Hamiltonian $H$, the symplectization radial coordinates $s$ and
other derivative estimates of the almost complex structures $J$.
(See \cite[Definitiion 4.5]{gps} and \cite[Appendix]{ganatra}, for example, and the relevant references
therein for the proofs of the $C^0$-estimates.)
On the other hand, under the aforementioned pair of the conditions,
our $C^0$-estimates in this sectorial framework are closely tied to the Hamiltonian geometry and are derived
from some tensorial calculations and the sign considerations \emph{without doing any estimates} for
the application of the maximum principle to ${\frak s}_{\varphi,\kappa}$.
\end{remark}

\section{Confinement for closed-open and open-closed maps}

As described in Subsection \ref{subsec:closed-open-map}, the relevant
Floer equation for the closed-open and open-closed maps is given as follows:
We will focus on the case when there is only one interior puncture.
In this case, we can identify $D^2 \setminus z_0^+$ with the semi-cylinder
$$
\RR \times (-\infty, 0].
$$
Under this identification, let $K :(-\infty,0] \times M \to \RR$ be
the given $\tau$-dependent sectorial Hamiltonians $\{K^\tau\}$.
\eqn\label{eq:H-at-ends}
K(\tau,x) = \begin{cases} H(x) \quad & \text{for $\tau > R$ with $R$ sufficiently large}\\
0 \quad & \text{for $0 \leq \tau \leq \delta$ for some $0 < \delta < R$}.
\end{cases}
\eqnd

We consider the equation
\eqn\label{eq:duXKJ}
\frac{\del u}{\del \tau} + J \left(\frac{\del u}{\del t} - X_K(u)\right) = 0.
\eqnd
The following is the confinement result for the solutions for \eqref{eq:duXKJ}.

\begin{theorem}\label{thm:closed-open}
Let $K = K(\tau,t)$ be a $\tau$-dependent (autonomous) sectorial Hamiltonians
$K_\tau = K(\tau,x)$ such that
$$
K(-\infty, x) = \chi(\tau) H(x)
$$
with an elongation function $\chi: (-\infty,0] \to [0,1]$ such that
\beastar
\chi(\tau) & = &
\begin{cases}1 \quad & \text{ near $\tau = -\infty$}\\
0 \quad  &\text{ near $\tau = 0$}
\end{cases}\\
\chi'(\tau) & \leq & 0.
\eeastar
Then there exists some $r_0 > 0$ such that
any solution for $u$ with finite geometric energy satisfies
\eqn\label{eq:confirment-for-closed-open}
\Image u \subset ({\mathfrak s}_{\varphi,\kappa})^{-1}((-\infty,r_0]).
\eqnd
\end{theorem}
\begin{proof} By the same calculation as in the proof of Theorem \ref{thm:nonautonomous-confinement},
we derive the following differential inequality
$$
\Delta ({\mathfrak s}_{\varphi,\kappa}\circ u) \geq - \chi'(\tau) \frac{\del}{\del \tau} ({\mathfrak s}_{\varphi,\kappa} \circ u).
$$
Then we can apply the maximum and
the strong maximum principle respectively again as in Appendix \ref{sec:maximum-principle} and finish the proof.
\end{proof}

\begin{remark}
The case of more than one interior punctures do not carry the global coordinate $(\tau,t)$
as in the above case of one interior puncture and requires the usage of general $TM$-valued one form
$P(u)$ and involves the general equation
$$
(du- P(u))^{(0,1)}_J = 0
$$
in the setting of Hamiltonian fibration as in \cite{oh-tanaka:liouville-bundles}. We will elaborate this
elsewhere.
\end{remark}

\section{Construction of wrapped Floer cohomology}
\label{sec:construction-wrapped-Floer-cohomolgy}

In the present section, we duplicate the definition of wrapped Floer cohomology
$HW^*(L, K)$ provided in \cite[Section 3.4]{gps} employing our sectorial package
$(H,J)$ in place of the pair $(H,J)$ with dissipative Hamiltonian
$H$ and cylindrical almost complex structure $J$ used by Ganatra-Pardon-Shende in
\cite[Section 2.10 \& 2.11]{gps}. \emph{To be consistent with the notations from
\cite{gps}, we denote by $X$ a symplectic manifold and by $Y$ a contact manifold
in this section and the next.}

We start with borrowing the following discussion from \cite[Section 2.9]{gps} on
the cut-off Reeb dynamics on general contact manifold $(Y,\xi)$ with convex boundary.

\begin{defn} An element of the Lie algebra of the automorphism group $\Aut(Y,\del Y)$
is called a \emph{cut-off contact vector field} , i.e., those vanishing along $\del Y$.
\end{defn}

\begin{lemma}[Lemma 2.34 \cite{gps}] For every cut-off Reeb vector field denoted by ${\bf R}$
on  a contact manifold $(Y,\xi)$, there exists a compact subset of the interior of
$Y$ which intersects all periodic obits of ${\bf R}$.
\end{lemma}

Next choosing a vector field $V$ outward pointing along $\del Y$, we fix a splitting
$$
\nbhd(\del Y) \cong \del Y \times \R_{\geq 0}
$$
so that $V = - \frac{\del}{\del t}$. Then $\del Y$ has decomposition
 $$
 \del Y = (\del Y)_+ \cup (\del Y)_-
 $$
where  $V$ is positively (resp. negatively) on $(\del Y)_+$ (resp. on $(\del Y)_-$)
transverse to $\xi$. The two sets meet on the dividing set
$$
\Gamma_{\del Y}: = \{x \in \del Y \mid V(x) \in \xi|_x\} = (\del Y)_+ \cap (\del Y)_-
$$
of $\del Y$ along $V$.

\begin{defn} A function $M: \R_{\geq 0} \to \R_{\geq 0}$ defined on a germ of intervals inside
$\R_{\geq 0}$ containing $0$ is called \emph{admissible} if $M(0) = 0 = M'(0)$, $M''(0) > 0$
and $M'(t) > 0$ for $t > 0$.
\end{defn}

Then the following is proved in \cite{gps}.

\begin{prop}[Proposition 2.35 \cite{gps}] \label{prop:gps2.35}
Let $Y$ be a contact manifold with convex boundary.
There exists a canonically defined contractible family of pairs $(t,\alpha)$ consisting of
a choice of coordinates $\nbhd(\del Y) \cong \del Y \times \R_{\geq 0}$ and a $\frac{\del}{\del t}$
invariant contact form $\alpha$, such that for any admissible function $M: \R_{\geq 0} \to \R_{\geq 0}$,
the dynamics of the associated cut-off Reeb vector field defined over $\nbhd(\del Y)$ satisfy the following
property: For any trajectory $\gamma : \R \to \del Y \times \R_{\geq 0}$ with $\gamma = \gamma(\tau)$,
we have
\begin{itemize}
\item If $dt(\gamma'(\tau_0) \geq 0$, then $dt(\gamma'(\tau)) > 0$ for all $\tau < \tau_0$.
\item If $dt(\gamma'(\tau_0) \leq 0$, then $dt(\gamma'(\tau)) < 0$ for all $\tau > \tau_0$.
\end{itemize}
In particular, for all sufficiently small $\delta > 0$, no trajectory enters the region
$\del Y \times \R_{ 0 \leq t \leq \delta}$ and then exits.
\end{prop}

We apply the above discussion to $Y = \del_\infty X$ with $\del Y = \del_\infty X \cap \del X$.
By taking a symplectization end $\del_\infty X \times [0,\infty) \hookrightarrow  X$, we can find
a one-parameter subgroup $\Phi: \R_{\geq 0} \to \Ham(X)$ which at infinity corresponds to the
Reeb flow of the contact form $\alpha$ chosen above. More specifically let $s$ be the radial coordinate
for the decomposition $\nbhd(\del_\infty X): = \del_\infty X \times [0,\infty)$. Then $\Phi$ is
nothing but the Hamiltonian flow of the vector field $X_s$.

We denote by $\CH(X)$  the space of functions $H: X \to \R$ and recall the definition of
\emph{dissipative} (single) Hamiltonian employed in \cite[Section 4.2]{gps}.

\begin{defn}[Dissipative Hamiltonian]
 Let $H = H(t,x)$ be a one-periodic Hamiltonian. We call $H$ \emph{dissipative}
if Its time-one map $\phi_H^1$ is nondegenerate and
$$
\inf_{x \in \nbhd(\del_\infty X)} d(x,\phi_H^1(x)) \geq \epsilon > 0
$$
on some neighborhood $\nbhd(\del_\infty)$ at infinity.
\end{defn}

\begin{defn} [Adapted pair; Definition 4.3 \cite{gps}]
A pair $(H,J)$ of one-periodic pair  $H = H(t,x)$ and $J = J(t,x)$
is said to be \emph{adapted to } $\del X$ when both $H = \Re \pi$ and $\pi$ is $J$-holomorphic
over $\pi^{-1}(\C_{|\Re| \leq \epsilon})$ for some $\epsilon > 0$.
\end{defn}

\cite{gps} also uses the $n$-simplex versions of these definitions. Unfortunately these
two requirements are incompatible along the neighborhood $\nbhd(\del_\infty X \cap \nbhd X)$
which forces \cite{gps} to introduce the rather complex definitions of \emph{admissible} Hamiltonians
\cite[Definition 4.4]{gps} and \emph{dissipative $n$-simplex} of Hamiltonians \cite[Definition 4.5]{gps}.
All these notions are introduced to ensure $C^0$ bounds and the relevant
confinement results \emph{via the monotonicity arguments}.

Our sectorial package $(H,J)$ replaces the adapted pair $(H,J)$ appearing in
\cite[Definition 4.6]{gps} which is possible because the sectorial pair $(H,J)$
is amenable to the maximum principle and to the strong maximum principle when
paired with the $Z$-invariant Lagrangians at infinity.

With the above reparation, we now state the following

\begin{prop} Let $(t,\alpha)$ be the pair as given in Proposition \ref{prop:gps2.35}
and ${\bf R}$ be the associated Reeb vector field. We fix
the symplectization radial coordinate $s$ so that its Hamiltonian vector field $X_s$
corresponds to the Reeb vector field ${\bf R}$ on $\del_\infty X$ as described above.
For all sufficiently small $\delta > 0$, no Hamiltonian trajectory of sectorial
Hamiltonian $H = \rho(\frak{s}_{\varphi,\kappa})$ enters the region
$$
\nbhd_{\leq \delta}(\del X) \cap \nbhd(\del_\infty X)
$$
and then exits.
\end{prop}
\begin{proof} We recall that near a corner of codimension $k+1$
in $\nbhd(\del_\infty X \cap \del X)$ the end-profile
function $\frak{s}_{\varphi, \kappa} = \frak{s}_{\varphi} \circ \phi_\kappa$ where
\begin{enumerate}[(a)]
\item  $\frak{s}_{\varphi}$ has the form
$$
s_{k+1, \varphi} = -\log \varphi(R_1, \ldots, R_k, e^{-s}).
$$
\item $\phi_\kappa$ is the flow map of the vector field $X_t$
given in \eqref{eq:Xs}  for the splitting
$\nbhd(\del X) = \del X \times \C_{\Re \leq \delta}$.
\end{enumerate}

For the simplicity of exposition, we assume $k = 1$ so that
$$
s_{1+1,\varphi} = -\log \varphi(R, e^{-s}).
$$
By the property (b), we have
$$
\frak{s}_{\varphi,\kappa} = s_{1+1,\varphi}  = -\log \varphi\left(R, e^{-s}\right).
$$
Therefore, we compute
$$
d\frak{s}_{\varphi,\kappa}
 = \frac{1}{\varphi} \left(- \frac{\del \varphi}{\del x}(R,e^{-s})dR
 + \frac{\del \varphi}{\del y}(R,e^{-s})e^{-s}\, ds \right),
$$
and hence
$$
X_{\frak{s}_{\varphi,\kappa}} =
 \frac{1}{\varphi} \left(\frac{\del \varphi}{\del x}(R,e^{-s})\frac{\del}{\del I}
 + \frac{\del \varphi}{\del y}(R,e^{-s})e^{-s}\, X_s \right),
$$
From this expression where both $\frac{\del \varphi}{\del x_1}(R,e^{-s})$ and
$\frac{\del \varphi}{\del x_2}(R,e^{-s})e^{-s}$ are nonnegative,
and the property of the Hamiltonian vector field $X_s$ of
 $s$ spelled out in Proposition \ref{prop:gps2.35},  we conclude the proposition.
\end{proof}

Once this proposition in our disposal, the amenability of the maximum principle stated in
Section \ref{sec:structure-maps} and \ref{sec:continuation}
conclude construction of wrapped Floer cohomology by the same arguments
provided in \cite[Section 3.4]{gps} or by the more geometric argument in setting of
Liouville bundles \cite[Section 4]{oh-tanaka:actions}.

\section{Construction of covariant inclusion functor}
\label{sec:inclusion-functor}

In this section, we explain how the necessary confinement results
entering in the construction of the covariantly functorial morphisms defined in
\cite{gps} can be established utilizing the present sectorial Floer package.
We will focus on the geometric part of the construction since the algebraic part of
the construction is almost the same as that of \cite[Section 3.4-3.6]{gps}.

Let $X \hookrightarrow X'$ be an inclusion of Liouville sectors in the sense of \cite{gps}.
We will only consider the case
$$
X \cap \del X' = \emptyset
$$
for the simplicity of exposition which can be easily accommodate the case with $X \cap \del X' \neq \emptyset$.
In particular, there is some neighborhood $\nbhd(\del X')$ of $\del X'$ such that
$$
X \cap \nbhd(\del X') = \emptyset
$$
where
$$
\nbhd(\del X')
$$
is a neighborhood of $\del X'$ that is $Z$-invariant at infinity
We would like to mention that this choice of neighborhood depends only on the given pair
$X \subset X'$ with $X \cap \del X' = \emptyset$.

With this geometric preparation, our construction of the covariantly functorial inclusion functor follows the following
steps. We let $\varepsilon_0$ be the same constant that is used before.

\subsection{Step 1: Choice of smoothing profiles}

We take smoothing profiles and other Floer data of $X$ and $X'$ to define the wrapped Fukaya categories
$$
\Fuk(X), \quad \Fuk(X')
$$
respectively.

We may take the smoothing profiles and associated end-profile functions of $X$ and $X'$
as follows by suitably decreasing $\varepsilon_0> 0$ and $\epsilon> 0$ if necessary.

We first have a splitting data and the decomposition
$$
X = F \times \C_{[0, 4\varepsilon_0]} \bigcup X_{-4\varepsilon_0}
$$
where $ X_{-4\varepsilon_0}: = X \setminus F \times \C_{[0, 4\varepsilon_0]}$,
so that
$$
X \cap \nbhd_{4\varepsilon_0}(\del X') = \emptyset
$$
and $\nbhd_{4\varepsilon_0}(\del X') \subset \nbhd(\del X')$  for the neighborhood
$\nbhd(\del X')$ chosen in the beginning.

We take a contact type hypersurface $S_0'$ of $X'$ so that
$$
S_0' \pitchfork \del X
$$
and the associated radial function $s$ will also have the form
$$
s = s(I)
$$
on $F \times \C_{[0, 2\varepsilon_0]} \bigcup X_{-4\varepsilon_0} \subset X$. Recall by definition that $s$
also has the form $s(I')$ on $\nbhd_{\frac{3\varepsilon_0}{2}}(\del X')$.
Then we define the end-profile functions of ${\mathfrak s}_{\varphi,\kappa}^X$ and
${\mathfrak s}_{\varphi,\kappa}^{X'}$ simultaneously using this radial function $s$ defined on $X'$.

\subsection{Step 2: Choice of generating Lagrangian branes}

We take a collection of sectorial Lagrangians in $X$ indexed by a
countable set $\mathcal I$ containing all isotopy classes of Lagrangians.
(This collection corresponds to the indexing set $I$
used in \cite{gps}. Since we use the letter $I$ for the Hamiltonian $I$ already, we avoid
using  the same letter by changing $I$ to the calligraphic $I$ here.)

The following definition is given in \cite[Section 3.4]{gps}.

\begin{defn}[Positive wrapping category, \cite{gps}] The \emph{positive wrapping category}
denoted by
$$
(L \rightsquigarrow -)^+
$$
is defined as follows:
\begin{itemize}
\item {[{\bf Objects}]} isotopies of exact Lagrangians $\phi: L \rightsquigarrow L^w$,
\item {[{\bf Morphisms}]} The morphisms
$$
(\phi: L \rightsquigarrow L^w) \to (\phi': L \rightsquigarrow L^{w'})
$$
are the homotopy classes of positive isotopies of exact Lagrangians $\psi: L^w \to L^{w'}$ such that
$\phi \# \psi = \phi'$.
\end{itemize}
\end{defn}
Motivated by this, we also consider the following.

\begin{defn} We say two sectorial Lagrangians are in the same isotopy class if
there is a positive isotopy from one to the other.
\end{defn}

By applying a suitable Hamiltonian isotopy, we may assume that all elements of
$\mathcal I$ are contained in
$$
X \setminus \nbhd^Z_{\sqrt{\epsilon_1}}(\del X)
$$
for some sufficiently small $\epsilon_1> 0$.
\begin{remark}
\begin{enumerate}
\item
For example, we may use the Hamiltonian isotopy generated by $-X_I$ where
$I$ is the function that naturally appears when a splitting data is chosen later.
But we prefer to make this operation independent of such a choice.
\item Here one should regard $\epsilon_1$ as the $\epsilon$ that
appears in the construction of convex smoothing functions $\varphi$.
\end{enumerate}
\end{remark}

Then we consider the poset
$$
\cO = \ZZ_{\geq 0} \times {\mathcal I}
$$
as in \cite[Section 3.6]{gps}, except the changing the letter $I$ to $\mathcal I$ as mentioned above.

\subsection{Step 3: Defining a $\lambda$-sectorial almost complex structure}

We first take an associated $\lambda$-sectorial almost complex structure on
$X \setminus \nbhd_{[2\varepsilon_0, \infty)}(\del X)$ that satisfies the defining
properties $\lambda$-sectorial almost complex structures thereon by regarding
$$
X \setminus \nbhd_{[2\varepsilon_0, \infty)}(\del X) =: X_{2\varepsilon_0}
$$
itself as a Liouville sector with the decomposition
$$
X_{2\varepsilon_0}\setminus X_{4\varepsilon_0} = F \times \C_{[4 \varepsilon_0,\infty)}.
$$
Then we extend this arbitrarily to a $\lambda$-sectorial $J$ for the whole $X'$.
We take
$$
2\varepsilon_0 < \sqrt{\epsilon_1}.
$$

\subsection{Step 4: Enlarging the generating collection to $X'$}
Next we consider the corresponding data for the sector $X'$. We similarly choose ${\mathcal I}'$ so that
all elements thereof are contained both in
$$
X' \setminus \nbhd^Z_{\sqrt{\epsilon_1}}(\del X) \subset X' \setminus \nbhd_{2\varepsilon_0}(\del X').
$$
We set
$$
\cO' = \ZZ_{\geq 0} \times ({\mathcal I}' \cup \cO).
$$
\subsection{Step 5: Study of pseudoholomorphic curves on $X'$}

Now using the above chosen Floer data $(J,H)$ on $X'$,
we consider the $J$-holomorphic equations on the marked discs $\dot \Sigma$
$$
\dot \Sigma = D^2 \setminus \{z_0, \ldots, z_m\}, \quad z_0, \ldots, z_m \in \del D^2
$$
on $X'$ with the corresponding boundary conditions with the Lagrangians
$L_0, \cdots, L_m$ from ${\mathcal I}'$. Then thanks to the properties of the sectorial Floer package
$$
(J,H),  \quad \{L_0, \cdots, L_m\}
$$
established in the present paper, the maximum (and strong maximum principle)
give rise to the covariantly functorial inclusion functor in the same was as
\cite{gps} does.

\appendix

\section{Discussion}

Here we would like to attract attention of readers to the implication of our sectorial Floer package
to the geometric framework needed to perform the analysis of pseudoholomorphic curves on noncompact
symplectic manifolds which have been commonly adopted in relation to the study of wrapped Fukaya
category on the Liouville manifolds (or the sectors). We recall that on noncompact manifold setting, the $C^0$ estimate must
precede for the study of compactness properties of pseudoholomorphic curves.

\subsection{Monotonicity argument relies on the energy estimates}

A Liouville equivalence between two Liouville manifolds $(M_1, \lambda_1)$ and $(M_2,\lambda_2)$
is a diffeomorphism $\phi: M_1 \to M_2$ satisfying $\phi^*\lambda = \lambda + df$
for a function $f: M_1 \to \RR$ \emph{with compact support}. (See \cite{abouzaid-seidel}.) This support condition enters
the analysis of pseudoholomorphic curves through the study of relevant action estimates
which in turn is used to get the energy bound.
When the maximum principle does not apply, using the monotonicity is a tool to establish
for the \emph{autonomous} (that is time-independent) Floer equation.
Similarly for the study of Floer homology
theory of \emph{exact Lagrangian submanifolds} as in the study of wrapped Fukaya category,
a choice of Liouville primitive $g$ of $i^*\lambda = dg$ crucially enters first to
equip the relevant Floer trajectory equation with a gradient structure. (See \cite{oh:jdg}, \cite{kasturirangan-oh}.)
This is the reason why the pair $(L,g)$, not just $L$ itself, has been a part of essential brane data for
the Floer theory of exact Lagrangian submanifolds (see \cite{seidel-long-exact-sequence}).
But this approach requires the uniform energy estimate to proceed and hence comes the requirement of
some fixed behavior of the primitives $f$ and $g$ at infinity. This is what has been adopted in the existing literature
and also in \cite{gps}. For the case of noncompact Lagrangian submanifolds such as the case of cylindrical Lagrangians
in Liouville manifolds researchers have been keen to the support behavior of the primitive $g$
so that the primitive remains to be compactly supported after the various operations such as taking the
product of two Lagrangian branes $(L_1,g_1)$ and $(L_2,g_2)$ in the product of
Liouville manifolds $M_1 \times M_2$. (See \cite{gps-2} for example.) Obviously
the product operation destroys this compact support property.

\subsection{Maximum principle does not rely on the energy estimates}

Of course, the aforementioned compact support requirement for the above functions $f$ or for $g$
do not enter the Hamiltonian flow equation or the Floer equation themselves other than through the action estimates
and hence through the energy estimates.
While the presence of the function $g$ is fundamental for the Floer theory to start, a choice thereof
does not affect the analysis of pseudoholomorphic curves for the compact
Lagrangian submanifolds, \emph{unless one is interested in a quantitative aspect of the Floer theory.}

When the maximum principle applies, it is the most efficient way of establishing the $C^0$ estimates
because it first requires only the local property of the second derivative estimates and are
independent of the energy estimates \emph{once the asymptotic limits are prescribed}. More importantly
it does not depend on the choice of primitives and even applies equally well to non-exact (cylindrical)
Lagrangian submanifolds.
Since we use the maximum principle to establish the $C^0$ estimates which does not depend on the
choice of primitives, we need no longer require this compact support condition for the primitives.
(See \cite{oh:jdg} for an early advocation of the \emph{asymptotically constant} framework.)
This enables us to widen the class of Hamiltonians in the literature
e.g., those used in \cite{groman}, \cite{gps} to that of the \emph{simply defined} sectorial
ones we introduce in the present paper once the $\lambda$-sectorial almost complex structure is paired with them
in the study of Floer equations. The Liouville diffeomorphisms that we used to deform
the Liouville form  in Section \ref{sec:deformation-lambda} implicitly plays a crucial role
along the way in the process of the construction of $\lambda$-sectorial almost complex structures.

%

\section{Giroux's ideal completion and the asymptotically linear framework}
\label{sec:ACI-framework}

Now we explain an elegant way of packaging all the definitions of sectorial objects
utilizing Giroux's notion of ideal Liouville domain which we generalize to
the case of Liouville sector, \emph{provided} we employ the analytical framework of
pseudoholomorphic curves on open symplectic manifolds developed by Bao
in \cite{bao1,bao2} (see also \cite{oh-wang2}) which extends the existing
cylindrical-at-infinity (CI) analytical framework to that of more flexible
asymptotically-cylindrical-at-infinity (ACI).
The point of view of Giroux's ideal completion as the completion of
the Liouville domain makes natural to extend the CI category to the ACI category.

\subsection{Giroux' notion of ideal Liouville domain}

We recall Giroux's notion of the \emph{ideal completion} of the Liouville domain $(M^0,\lambda^0)$.
We closely follow Giroux' exposition, in particular his notations in this section for the
comparison purpose.

\begin{defn}[Ideal Liouville Domains \cite{giroux}] An \emph{ideal
Liouville domain} $(F,\omega)$ is a domain $F$ endowed with an ideal Liouville structure $\omega$.
The \emph{ideal Liouville structure} is an exact symplectic form on $\Int F$ admitting a primitive
$\lambda$ such that: For some (and then any) function $u:F \to \RR_{\geq 0}$ with regular level set
$\del F = \{u=0\}$, the product $u\lambda$ extends to a smooth one-form on $F$ which induces a
contact form on $\del F$.
\end{defn}
\emph{We alert readers that the $F$ appearing in this definition has nothing to do with $F$ in the main text.}

We first recall the construction of ideal completion from \cite[Example 9]{giroux}.
Let $u:M^0 \to \RR_{\geq 0}$ be a function with the following properties:
\begin{itemize}
\item $u$ admits $C = \del M^0$ as its regular level set $\{u=0\}$,
\item $X_\lambda [\log u] < 1$ at every point in $\Int M^0$.
\end{itemize}
Then the form $\omega^0: = d(\lambda^0/u)$ defines a symplectic form on $\Int M^0$.
The pair $(M^0,\omega^0)$ is an \emph{ideal Liouville domain} which Giroux calls
the \emph{ideal completion} of the Liouville domain $(M^0,\lambda^0)$.

We write
\begin{itemize}
\item $\Aut(M^0,\omega^0)$ for the group of diffeomorphisms of $M^0$ preserving $\omega^0$
on $\Int M^0$,
\item  $\Aut_{\del}(M^0,\omega^0) \subset \Aut(M^0,\omega^0)$ for the subgroup of diffeomorphisms fixing $C = \del M^0$
pointwise, and
\item  $\Aut(C,\xi)$ for the group of contactomorphisms of $(C,\xi)$.
\end{itemize}

The following propositions are proved by Giroux.

\begin{prop}[Proposition 7 \cite{giroux}]\label{prop:giroux-Serre-fibraion} The restriction homomorphism
$$
\Aut(M^0,\omega^0) \to \Aut(C,\xi)
$$
is a Serre fibration, with associated long exact sequence of homotopy groups
\beastar
\ldots \pi_k(\Aut_{\del}(M^0,\omega^0)) \to \pi_k(\Aut(M^0,\omega^0)) \to \\
\to \pi_k(\Aut(C,\xi))
\to \pi_{k-1}(\Aut_{\del}(M^0,\omega^0)) \to \ldots.
\eeastar
\end{prop}

Let $u: M^0 \to \RR_{\geq 0}$ be a defining function of $C = \del M^0$ as above.

\begin{lemma}[Proposition 3 \cite{giroux}] Let $\beta^\flat$ be the vector field dual to
the form $\beta = \lambda/u$  defined by $\beta^\flat \rfloor \omega^0 = \beta$.
Then it is complete on $\Int M^0$ and defines a unique
embedding $\iota_\beta: \overline{SC} \to M^0$ such that
$\iota_\beta|_{C} = id$ and $\iota_\beta^*\beta = \lambda_\xi$ and its image is
an open collar neighborhood $U$ of $C = \del M^0$.
\end{lemma}

Therefore $\iota_\beta$ provides a natural decomposition
$$
\Int M^0 = (M^0 \setminus U) \cup U
$$
such that $\iota_\beta: (SC,d\lambda_\xi) \to (F,d\beta)$ is a
Liouville embedding.

In regard to the Liouville isomorphisms, an upshot of Giroux's ideal completion is that the group
$\Aut(M^0,\omega^0)$ corresponds to the group of asymptotically-linear-at-infinity symplectomorphisms, i.e.,
\emph{those converging to linear ones  exponentially fast at infinity}.
(See \cite{bao1} for the relevant definition of
asymptotic linearity or asymptotic cylindricality.)
This is an amplification of the standard group of symplectomorphisms
\emph{linear-at-infinity}.  The Giroux' framework of ideal
Liouville domains allows one to elegantly package this asymptotically linear framework.
Therefore it is aesthetically very satisfying to see that Bao's asymptotically-linear-at-infinity framework
enables one to make all the definitions of sectorial objects very simple. Explanation of this is now in order.

\subsection{ACI-sectorial Floer package in Liouville sectors}

The definition of \emph{ideal sectorial domains} is verbatim the same as that of
\cite{giroux} for the case of Liouville domains (without corners) (under the assumption
that the Liouville vector field is tangent to near infinity), and then
a Liouville sector is nothing but $M \setminus \del_\infty M$, an ideal completion of sectorial domain.

In this ACI point of view, it is very simple to
give the definitions of sectorial almost complex structures, Hamiltonians and Lagrangians.

We borrow the discussion from \cite[Example 9]{giroux} here:
Let $(W,\lambda)$ be a sectorial domain and
$u: W \to \RR_{\geq 0}$ be
a defining function of $\del_\infty W$ with the properties laid out as above.
Then $(F,\omega_{[u]})$ with $\omega_{[u]}: = d(\lambda/u)$ defines an ideal Liouville domain
in the sense of \cite[Definition 1]{giroux}.

\begin{defn}[Ceiling-profile function] Let $(W,\lambda)$ be a sectorial domain.
We call a function $u$ above the \emph{ceiling-profile function} the
(asymptotically-linear) open Liouville sector $(M, \omega)$ if
$$
M = W \setminus \del_\infty W, \quad \omega = d(\lambda/u).
$$
\end{defn}

\begin{defn}[ACI-sectorial Floer package] Let $(W,\lambda)$ be a sectorial domain and $(M, \omega)$ be
an open Liouville sector with the ceiling-profile function $u$, i.e., $\omega = d(\lambda/u)$.
\begin{enumerate}
\item {[ACI-sectorial Lagrangian]}
We call a Lagrangian submanifold $L \subset (M,\omega)$
a \emph{\text{\rm ACI}-sectorial} if there exists a smooth Lagrangian submanifold $\overline L$ in $W$ such that
\begin{itemize}
\item $\overline L$ intersects $\del_\infty W$ transversely,
\item $\overline L \cap \del W = \emptyset$.
\item $\overline L$ is $\beta^\flat$-invariant near $\del_\infty W$.
\end{itemize}
\item {[ACI-sectorial almost complex structure]} An $\omega$-tame almost complex structure $J$ on $M$ is
called \emph{\text{\rm ACI}-sectorial}, if $J$ makes the union $\del_\infty W \cup \del W$ $J$ pseudoconvex so that its defining exhaustion function is Liouville pseudoconvex.
\item {[ACI-sectorial Hamiltonian]} A function $H: M \to \RR$ is called \emph{\text{\rm ACI}-sectorial} if
it smoothly extends to $W$ and satisfies similar conditions as in Definition \ref{defn:sectorial-H}
in terms of the $J$-convex exhaustion function $\psi$.
\end{enumerate}
\end{defn}

\begin{remark} When $\alpha < 1$ with $I$ satisfying $Z[\log I] = \alpha$, the function $I - N$ for
a sufficiently large $N> 0$ can play the role of the profile function $u$ above on a neighborhood of the boundary of $\del M$,
provided the level set $I^{-1}(N)$ is compact in the neighborhood.
So if we want to take the ACI point of view, the case $\alpha = 1$ is special in that
the function $I$ does not belong to the set of ceiling-profile function: the two form $d(\lambda/I)$
has nontrivial kernel which is nothing but the span of the Liouville vector field.
Varying $\alpha$ will be also important, when the family story as in the bundle
of Liouville sectors considered in \cite{oh-tanaka:liouville-bundles} and it is good to have
the flexibility of choosing the eccentricity of $\alpha$ in the story, and also raises a question
what would be the implication when the family cannot avoid including the extreme eccentricity
$\alpha = 1$.
\end{remark}

\section{Proof of Proposition \ref{prop:tildefepsilon}}
\label{sec:tildefepsilon}

Let $\delta'>  \delta > 0$ and consider any smooth function
$f:[\delta, \delta') \to \R$ such that
$$
\begin{cases}
f(\delta) > 0, \quad f'(\delta) < 0, \quad  f''(\delta) > 0, \\
f(\delta) + \delta f'(\delta) = 0.
\end{cases}
$$
Let $\tau_{\delta} : = \delta -\frac{f(\delta)}{f'(\delta)}$,
i.e., the $x$-intersection of the tangent line
at $x = \delta$ of the graph of $f$.

Proposition \ref{prop:tildefepsilon} will be an immediate consequence of the following general lemma.

\begin{lemma}\label{lem:referee} Let $\delta > 0$, $f$ and $\tau_\delta> 0$ be as above.
Then for any $\tau > \tau_\delta$, there exists a smooth function $\widetilde f:[\delta, \infty) \to \R$ satisfying
the following conditions.
\begin{enumerate}
\item $f$ and $\widetilde f$ coincide in a neighborhood of $\delta$,
\item For each $x \in [0,\tau)$,
$$
\begin{cases}
\widetilde f(\delta) > 0, \quad \widetilde f'(\delta) < 0, \quad
\widetilde f''(\delta) > 0, \\
\widetilde f(\delta) + \delta \widetilde f'(\delta) = 0.
\end{cases}
$$
\item For $x \geq \tau$, $\widetilde f(x) = 0$.
\item For all $x \geq \delta $, $2\widetilde f'(x) + x \widetilde f''(x) \geq 0$.
\end{enumerate}
\end{lemma}
\begin{proof} By a suitable pre/post-composing $f$ with appropriate linear functions, it is
enough to consider the case
$$
f(\delta) = 1, \quad f'(\delta) = -1.
$$
Writing $-f' =: g$,  it is enough to find a smooth function
$\widetilde g: [\delta,\infty) \to [\delta,\infty)$ such that
\begin{enumerate}[(a)]
\item $g$ and $\widetilde g$ coincide in a neighborhood of $\delta$,
\item For $x \in [\delta,\tau)$, $0 < \widetilde g(x) < 1,
\, \widetilde g'(\delta) < 0$,
\item For $x \geq \tau$, $\widetilde g(x) = 0$.
\item $\int_\delta^\tau \widetilde g = 1$.
\item $\widetilde g'(x) < 0$ and $x (\log \widetilde g)' \geq - 2$ for $ \delta \leq x < \tau$.
\end{enumerate}
Then we just set
$$
\widetilde f(x): = 1 - \int_\delta^x \widetilde g(t)\, dt
$$
which will finish the proof.

Now we give construction of $\widetilde g$. Note that $g(\delta)
= -f'(\delta) = 1$ and $g'(\delta) = - f''(\delta) < 0$.  Furthermore by the normalization hypothesis
$f(\delta) = 1, \, f'(\delta) = -1$ imposed in the beginning, we have
$\tau_\delta = \delta -\frac{f(\delta)}{f'(\delta)} =1 + \delta $
and hence
$$
\tau +\delta > \tau_\delta + \delta = 1 +2\delta > \tau_\delta= 1+\delta.
$$
Note that the condition (e) can be rewritten as
\be\label{eq:2/x+}
\frac{2}{x} + \frac{\widetilde g'(x)}{\widetilde g(x)} \leq 0
\ee
which is equivalent to $(\log x^2 \widetilde g(x))' \leq 0$. Since $\log$ is an increasing
function, this is equivalent to $x^2 \widetilde g(x)$ is a decreasing function. 
Recalling the requirement 
$$
\begin{cases}
\widetilde g(x) = g(x)  \quad & \text{near } \, x = \delta \\
\widetilde g(x) \equiv 0 \quad & \text{for } \, x \geq \tau, 
\end{cases}
$$
we choose any smooth decreasing function $\varphi$ satisfying
$$
\begin{cases}
\varphi(x) = 0 \quad  & \text{for },  x \geq \tau\\
\varphi(x) = x^2 g(x)  \quad & \text{near }\, x = \delta.
\end{cases}
$$
Then we set $\widetilde g(x) = \frac{\varphi(x)}{x^2}$.

By the area consideration of the graph of $\widetilde g$, we can make a suitable choice of
$\varphi$ so that it satisfies all the requirements $(a) - (e)$
above, \emph{as long as we choose $\tau$ so that say $\tau > 3$.} We will set $\tau = 5$.
This finishes the proof.
\end{proof}

Now we apply this construction to the function
$$
f: = f_\e(x), \quad \delta: = \sqrt{\epsilon/2}
$$
and get $\widetilde f$
on $ x \geq \sqrt{\epsilon/2})$. This finishes the proof of Proposition
\ref{prop:tildefepsilon} (1) -- (3) immediately. To prove (4), we put
$$
h(x) = f_\e(x) + x f_\e'(x).
$$
Then recalling $f_\epsilon(x) = \frac{\epsilon}{2x}$
on $[\sqrt{\epsilon}/4,2\sqrt{\epsilon}]$, an
explicit calculation with the function $\frac{\epsilon}{2x}$
proves $h(0) = 0$. Therefore to prove (4) of Proposition \ref{prop:tildefepsilon}
it is enough to show $h'(x) \geq 0$ for $x \geq \sqrt{\epsilon/2}$.
But we compute
\beastar
h'(x) & = & 2f_\e'(x) + xf_\e''(x) = -2 g(x) - x g'(x) \\
& = & =
 - g(x)\left(2 + \frac{xg'(x)}{g(x)}\left(x\right)\right) 
\eeastar
which is positive by Lemma \ref{lem:referee} (e) (or more directly by \eqref{eq:2/x+})
since $g > 0$ on $[\sqrt{\epsilon/2}, \tau)$. This finishes the proof of the proposition.

\section{Review of the proof of Gray's stability theorem}
\label{sec:gray-theorem}

In this appendix, we recall the details of the proof of Gray's
stability theorem \cite{gray} via Moser's deformation method.
\begin{theorem}[Gray's stability theorem]
 Suppose $\xi_t$ be a homotopy of contact structures
for all $t \in [0,1]$. Then there is a diffeomorphism $\psi_t: Q \to Q$
such that $\psi_t^*(\xi_t) = \xi_0$ for all $t \in [0,1]$.
\end{theorem}
\begin{proof}  Let $\alpha_t$ be a smooth family of contact forms with
$\ker \alpha_t = \xi_t$. We would like to find a diffeomorphism
$\psi_t$ with
 $$
 \psi_t^* \alpha_t = e^{g_t}\, \alpha_0
 $$
 for some known function $g_t$ on $Q$.

We differentiate the equation $\psi_t^*(\alpha_t) = e^{g_t} \alpha$ in $t$ and obtain
By setting $h_t: = \frac{\del g_t}{\del t} \circ \psi_t^{-1}$, this equation is
transformed into
$$
\psi_t^*\left(d (X_t \rfloor \alpha_t) + X_t \rfloor d\alpha_t + \frac{\del \alpha_t}{\del t} \right) = \psi_t^*(h_t \alpha_t).
$$
Now \emph{if we can take $X_t$ so that $X_t \in \xi_t$}, this equation is reduced to
\be\label{eq:Xt-equation}
\frac{\del \alpha_t}{\del t} + X_t \rfloor d\alpha_t = h_t \alpha_t.
\ee
By applying the Reeb vector field $R_{\alpha_t}$, we obtain
\be\label{eq:mut-equation}
h_t = \frac{\del \alpha_t}{\del t}(R_{\alpha_t}).
\ee
Then we can indeed take $X_t \in \xi_t$ by the defining equation
solving \eqref{eq:Xt-equation} by
\be\label{eq:Xtdalphat}
X_t \rfloor d\alpha_t = h_t \alpha_t -
\frac{\del \alpha_t}{\del t}(R_{\alpha_t}).
\ee
We emphasize the fact that the right hand side is already solved
from \eqref{eq:mut-equation} for $h_t$ because
$\frac{\del \alpha_t}{\del t}(R_{\alpha_t})$ is determined from
the given isotopy $\alpha_t$ of one-forms such that $\ker \alpha_t = \xi_t$.
\end{proof}
\section{Strong maximum principle and $Z$-invariance of Lagrangians}
\label{sec:maximum-principle}

In this section, we provide the details on how we apply the maximum (and strong maximum) principle to
the circumstances of Theorems \ref{thm:nonautonomous-confinement} and \ref{thm:closed-open}.

In both cases, we have derived a differential inequality of
the following type.
\eqn\label{eq:Laplacian}
\begin{cases}
\Delta({\mathfrak s}_{\varphi,\kappa} \circ u) \geq c(\tau,t)
\frac{\del}{\del \tau}({\mathfrak s}_{\varphi,\kappa} \circ u)\\
u(\tau,i) \subset L_i, \quad i=0, \, 1
\end{cases}
\eqnd
for the function
$$
c(\tau,t) = - (\rho^{\chi(\tau)})'({\mathfrak s}_{\varphi,\kappa}(u(\tau,t))
$$
in the circumstances of Theorem \ref{thm:nonautonomous-confinement} and \ref{thm:closed-open}.
 Suppose to the contrary that there exists some
local maximum point $(\tau_0,i)$, say, $i = 0$ so that $u(\tau_0,0) \in L_0$ at which
the vector field $Z$ is tangent to $L_0$. Then we have
\eqn\label{eq:dudtau=0}
\frac{\del}{\del \tau}({\mathfrak s}_{\varphi,\kappa} \circ u)(\tau_0,0) = 0
\eqnd
which implies $\frac{\del u}{\del \tau}$ is tangent to
the intersection
$$
s^{-1}(r_0) \cap L_0, \quad r_0: = {\mathfrak s}_{\varphi,\kappa}(u(\tau_0,0)).
$$
On the other hand, $Z$ is tangent to $L_0$ at $u(\tau_0,0)$ by the $Z$-invariance of $L$ thereon.
In particular $T_{u(\tau_0,0)} ({\mathfrak s}_{\varphi,\kappa}^{-1}(r_0) \cap L_0) \subset \ker (\lambda|_{{\mathfrak s}_{\varphi,\kappa}^{-1}(r_0)})$.
Then we rewrite
$$
\frac{\del}{\del t}({\mathfrak s}_{\varphi,\kappa} \circ u)(\tau_0,0)
= d{\mathfrak s}_{\varphi,\kappa} \left(\frac{\del u}{\del t}\right)
= - (d{\mathfrak s}_{\varphi,\kappa} \circ J)J\left(\frac{\del u}{\del t}\right)
= \lambda\left(\frac{\del u}{\del \tau}\right) = 0
$$
where the last vanishing follows from the $Z$-invariance of $L$.
On the other hand, \eqref{eq:dudtau=0} and \eqref{eq:Laplacian} also imply
$\Delta({\mathfrak s}_{\varphi,\kappa} \circ u) \geq 0$. This is a contradiction to the strong
maximum principle which should imply
$$
\frac{\del}{\del t}({\mathfrak s}_{\varphi,\kappa} \circ u)(\tau_0,0) > 0
$$
unless ${\mathfrak s}_{\varphi,\kappa} \circ u$ is a constant function in which case there is
nothing to prove. This finishes the proof.

\bibliographystyle{amsalpha}

\bibliography{biblio-oh}

\end{document}